\documentclass[11pt,twoside]{article}
\usepackage{atmp,latexsym}
\copyrightnotice{2002}{6}{873}{978}	
\newtheorem{thm}{Theorem}[section]
\newtheorem{prop}[thm]{Proposition}
\newtheorem{claim}[thm]{Claim}

\newtheorem{thm-defi}[thm]{Theorem/Definition}
\newtheorem{example}[thm]{Example}
\newtheorem{cor}[thm]{Corollary}

\newtheorem{lemma}[thm]{Lemma}
\newtheorem{lemma-definition}[thm]{Lemma/Definition}
\newtheorem{defi}[thm]{Definition}
\newtheorem{rem}[thm]{Remark}

\newtheorem{condition}[thm]{Condition}
\newcommand{\A}{{\cal A}}
\newcommand{\B}{{\cal B}}
\newcommand{\C}{{\cal C}}
\newcommand{\D}{{\cal D}}
\newcommand{\E}{{\cal E}}
\newcommand{\F}{{\cal F}}

\newcommand{\U}{{\cal U}}

\renewcommand{\P}{{\cal P}}
\newcommand{\PP}{{\Bbb P}}
\newcommand{\linsys}[1]{{\mid}#1{\mid}}

\newcommand{\IsomRightArrow}{\stackrel{\cong}{\rightarrow}}
\newcommand{\LongIsomRightArrow}{\stackrel{\cong}{\longrightarrow}}
\newcommand{\RightArrowOf}[1]{\stackrel{#1}{\rightarrow}}

\newcommand{\LongLeftArrowOf}[1]{\stackrel{#1}{\longleftarrow}}

\newcommand{\LongRightArrowOf}[1]{\stackrel{#1}{\longrightarrow}}

\newcommand{\StructureSheaf}[1]{{\cal O}_{#1}}
\newcommand{\CompletedStructureSheaf}[1]{\hat{{\cal O}}_{#1}}
\newcommand{\CompletedFunctionField}[1]{\hat{{\cal K}}_{#1}}
\newcommand{\EndProof}{\hfill  $\Box$}
\newcommand{\restricted}[2]{#1_{\mid_{#2}}}

\newcommand{\tr}{{\rm tr}}
\newcommand{\Pic}{{\rm Pic}}
\newcommand{\Sym}{{\rm Sym}}
\newcommand{\Ext}{{\rm Ext}}

\newcommand{\Hom}{{\rm Hom}}
\newcommand{\Aut}{{\rm Aut}}
\newcommand{\End}{{\rm End}}
\newcommand{\Isom}{{\rm Isom}}
\newcommand{\Mod}{{\rm mod}}
\newcommand{\Ord}{{\rm ord}}
\newcommand{\Abs}[1]{{\mid\!#1\!\mid}}

\newcommand{\SheafHom}{{\cal H}om}
\newcommand{\SheafEnd}{{\cal E}nd}
\newcommand{\SheafExt}{{\cal E}xt}
\newcommand{\RelExt}{{\cal E}xt}

\newcommand{\Wedge}[1]{\stackrel{#1}{\wedge}}
\newcommand{\Span}{{\rm span}}

\newcommand{\Higgs}{{\rm Higgs}}
\newcommand{\Hecke}{{\rm Hecke}}
\newcommand{\Prym}{{\rm Prym}}
\newcommand{\orbit}{{\Bbb O}}
\newcommand{\CoadjointOrbit}{{\Bbb O}'}
\newcommand{\rank}{{\rm rank}}

\newcommand{\coker}{{\rm coker}}
\newcommand{\bba}{{\Bbb A}}

\newcommand{\bbc}{{\Bbb C}}

\newcommand{\bbh}{{\Bbb H}}
\newcommand{\bbi}{{\Bbb I}}

\newcommand{\bbp}{{\Bbb P}}
\newcommand{\bbq}{{\Bbb Q}}
\newcommand{\bbr}{{\Bbb R}}

\newcommand{\bbz}{{\Bbb Z}}
\newcommand{\Got}[1]{{\frak #1}}
\newcommand{\ga}{{\Got a}}
\renewcommand{\gg}{{\Got g}}
\newcommand{\ggl}{{\Got gl}}

\newcommand{\gsl}{{\Got sl}}
\newcommand{\gso}{{\Got so}}

\newcommand{\gh}{{\Got h}}

\newcommand{\gb}{{\Got b}}
\newcommand{\gc}{{\Got c}}
\newcommand{\gt}{{\Got t}}
\newcommand{\gz}{{\Got z}}

\input amssymb.sty

\begin{document}

\title{Elliptic Sklyanin integrable systems for arbitrary reductive groups}
\url{math.AG/0203031}		
\vspace*{-.3in}
 
\author{J.C. Hurtubise\footnote{
The first author of this article would like 
to thank NSERC and FCAR for their support}}
\vspace*{-.3in}

\address{Department of Mathematics and Statistics\\ 
McGill University
}

\addressemail{hurtubis@math.mcgill.ca}	
\vspace*{-.3in}

\author{E. Markman
\footnote{Partially supported by NSF grant number DMS-9802532}}
\vspace*{-.3in}

\address{Department of Mathematics and Statistics\\ 
University of Massachusetts at Amherst
}
\addressemail{markman@math.umass.edu}	
\markboth{\it Elliptic Sklyanin integrable systems \ldots}
{\it J.C. Hurtubise and E. Markman}

\begin{abstract}
We present the analogue, for an arbitrary complex reductive group $G$,
of the elliptic integrable systems of Sklyanin. The Sklyanin integrable
systems were originally constructed on symplectic leaves, 
of a quadratic Poisson structure, on a loop group of type A. 
The phase space, of our integrable systems, is a group-like analogue of the 
Hitchin system over an elliptic curve $\Sigma$. The phase space
is the moduli space of pairs $(P,\varphi)$, 
where $P$ is a principal $G$-bundle on $\Sigma$, 
and $\varphi$ is a meromorphic section of the adjoint group bundle. 
The Poisson structure, on the moduli space, 
is related to the Poisson structures on loop groupoids,
constructed by Etingof and Varchenko, 
using Felder's elliptic solutions of the 
Classical Dynamical Yang-Baxter Equation. 
\end{abstract}

\cutpage


\tableofcontents 


\section{Introduction}
\label{sec-introduction}

We present in this paper the analogue, for an arbitrary complex 
reductive group $G$, of the elliptic integrable systems of Sklyanin.
The Sklyanin integrable systems were constructed on symplectic leaves of a 
quadratic Poisson structure on a loop group of type A
\cite{belavin-drinfeld,sklyanin,reiman-semenov}. 

Let $\Sigma$ be a complex elliptic curve. 
We study the moduli spaces $M_\Sigma(G,c)$ of pairs $(P,\varphi)$, 
where $P$ is a principal $G$-bundle over $\Sigma$, of topological type $c$, 
and $\varphi$ is a meromorphic section of its adjoint {\em group} bundle.
$M_\Sigma(G,c)$ are infinite dimensional ind-varieties.
These moduli spaces may be regarded as algebro-geometric analogues of a 
meromorphic loop group. The main results of this paper are:

\begin{enumerate}
\item
$M_\Sigma(G,c)$ admits an algebraic Poisson structure 
(Theorem \ref{thm-infinite-dimensional-Poisson-moduli}). 
The algebraic Poisson structure coincides with the reduction 
of the one constructed by Etingof and Varchenko on a loop groupoid, 
using solutions of the Classical Dynamical Yang-Baxter Equations 
(Theorem \ref{thm-comparison-with-cdybe-simply-connected-case}).
\item
The symplectic leaves $M_\Sigma(G,c,\orbit)$ 
of $M_\Sigma(G,c)$ are finite dimensional 
(Section \ref{sec-symplectic-leaves-folliation}). 
Their dimension is calculated in Corollary 
\ref{cor-dimension-of-a-symplectic-leaf}. 
The data $\orbit$, parametrizing the symplectic leaves, is a decorated divisor.
It consists of 
a labeling, of finitely many points of the elliptic curve $\Sigma$, 
by finite-dimensional irreducible representations of the Langland's 
dual group (i.e., by orbits of the co-character lattice of $G$, 
under the action of the Weyl group).
\item
The symplectic leaves $M_\Sigma(G,c,\orbit)$ are algebraically completely 
integrable hamiltonian systems. The generic Lagrangian leaf is an
abelian variety, which admits a precise modular description in terms of 
$T$-bundles on a branched $W$-cover of the elliptic curve $\Sigma$
(Theorem \ref{thm-general-complete-integrability}). 
\end{enumerate}

En route, we prove that all smooth moduli spaces of simple sheaves, on
an algebraic Poisson surface, admit a Poisson structure (Theorem 
\ref{thm-Poisson-structure-on-moduli-of-sheaves-on-a-Poisson-surface}). 
This extends results of 1) S. Mukai and A. Tyurin, 
who constructed an anti-symmetric
tensor on these moduli spaces \cite{mukai-symplectic,tyurin}, 
2) S. Mukai, and Z. Ran,
who proved the Jacobi identity for the Mukai-Tyurin tensor in the case of 
locally free sheaves on a symplectic surface \cite{mukai-survey,ziv-ran},
and 3) F. Bottacin, who proved the Jacobi identity in the case of 
locally free sheaves on a Poisson surface \cite{Bo2}. 
The moduli spaces of sheaves, arising in the theory of integrable systems, 
consist of torsion sheaves, with pure one-dimensional 
support, on a Poisson surface (typically, a line bundle on a spectral curve). 

Sections \ref{sec-orbits} through 
\ref{sec-infinitesimal-study-of-M-infinity} 
provide an algebro-geometric treatment of the 
Poisson geometry and complete integrability 
of the moduli spaces $M_\Sigma(G,c)$.
In section \ref{sec-Hecke-correspondences} we construct a Poisson structure on
the Hecke correspondences over an elliptic curve. These are 
algebro-geometric analogues of orbits in the loop Grassmannian. 
In section \ref{sec-dynamical-yang-baxter}
we relate our constructions to the elliptic solutions of the 
Classical Dynamical Yang-Baxter equations. 
In section \ref{sec-products-of-partial-flag-varieties} we discuss the 
spectral curves arising from 
particular examples of symplectic leaves in $M_\Sigma(G,c)$. 
These examples are related to compact orbits of the loop Grassmannian, 
which are isomorphic to flag varieties of maximal parabolic subgroups. 
A particular case leads to a deformation of the Elliptic 
Calogero-Moser system of type A.

We expect most of the results of the paper to go through in the degenerate 
case of a nodal or cuspidal rational curve (a singular elliptic curve). 
The resulting systems will correspond to trigonometric, or rational, 
solutions of the dynamical Yang-Baxter equation.
Some aspects of 
the construction admit a generalization for curves of genus $\geq 2$,
but the moduli spaces are neither symplectic, nor Poisson.

\medskip
{\it Acknowledgments:} 
We  would like to thank Ron Donagi,
Pavel Etingof, Robert Friedman, John Harnad, Jim Humphreys, Ivan Mirkovic, 
Tony Pantev and Alexander Polishchuk for helpful communications.

\subsection{Poisson surfaces and the type A Sklyanin systems}
\label{sec-introduction-loop-group}
The starting point for our investigation was the relationship, 
established in \cite{HM-sklyanin}, between the 
Poisson geometry of three objects:
\begin{enumerate}
\item
The type A elliptic Sklyanin systems on the meromorphic loop group. 
\item
Poisson structures on moduli spaces $M_\Sigma(GL(n),c)$ of pairs $(E,\varphi)$,
consisting of a rank $n$ vector bundle $E$ of degree $c$ and a meromorphic,
generically invertible, endomorphism $\varphi:E\rightarrow E$. 
The rank $n$ and degree $c$ are assumed to be relatively prime. 
\item
The Poisson geometry 
of certain moduli spaces of sheaves on a ruled Poisson surface $(S,\psi_S)$. 
The surface $S$ is the compactification of a line bundle on the elliptic curve.
The Poisson structure $\psi_S$ is $\bbc^\times$-{\em invariant}. 
\end{enumerate}
The relationship between the first two objects is most direct once we pass
to the projectivization $P=\bbp{E}$ and the group $PGL(n)$. 
A Zariski open subset of $M_\Sigma(PGL(n),\bar{c})$ consists of pairs with
a fixed stable and rigid bundle $P$. Trivialize $P$ over 
a small disc around the origin of 
$\Sigma$ and restrict $\varphi$ to the boundary to obtain a loop in 
$PGL(n)$. 
The section $\varphi$ is determined as 
the meromorphic continuation of that loop. 

The relation between the last two objects is provided by 
the spectral dictionary
commonly used in the theory of integrable systems
\cite{hitchin-integrable-system,hurtubise-surfaces,simpson}.
The pair $(E,\varphi)$, with polar divisor $D$, determines a 
complete spectral curve $C$ in 
the total space of the line bundle $\StructureSheaf{\Sigma}(D)$ and
an eigensheaf $F$ on $C$. We regard $F$ as a torsion sheaf on 
a projective surface; the compactification 
$S=\bbp[\StructureSheaf{\Sigma}(D)\oplus \StructureSheaf{\Sigma}]$ 
of the total space of $\StructureSheaf{\Sigma}(D)$. 
The $\bbc^\times$-action on $\StructureSheaf{\Sigma}(D)$ extends to $S$
and $S$ has a $\bbc^\times$-invariant algebraic Poisson structure, which
is unique up to a scalar factor.

We expressed the Sklyanin Poisson structure on the loop group, in terms
of the deformation theory of the pairs $(E,\varphi)$. 
In fact, the formula (\ref{eq-Psi-M-D}), for the Poisson structure 
in this paper, is a direct generalization of the formula in the $PGL(n)$ 
case of \cite{HM-sklyanin}. 
The insight for obtaining this formula came from the theory
of moduli spaces of sheaves. Starting with a Poisson structure 
$\psi_S$ on a projective surface $S$,
Mukai and Tyurin constructed a Poisson
structure on the moduli spaces of sheaves on $S$ 
\cite{mukai-symplectic,tyurin}.
The formula of \cite{HM-sklyanin} is a translation of the Mukai-Tyurin tensor,
via the infinitesimal version of the above dictionary.

We would like to stress the significance of the $\bbc^\times$-invariance
of the Poisson structure $\psi_S$ on the surface $S$. 
Let us first contrast it with the linear Poisson structure on $S$, 
coming from the cotangent bundle via the isomorphism 
$\StructureSheaf{\Sigma}(D)\cong T^*_\Sigma(D)$. 
The $\bbc^\times$-invariant Poisson structure $\psi_S$ induces on 
$M_\Sigma(GL(n),c)$ a $\bbc^\times$-invariant Poisson structure. 
In contrast, the linear Poisson structure
induces on $M_\Sigma(GL(n),c)$ another, linear, Poisson structure, which 
corresponds to the generalized Hitchin systems 
(see \cite{hitchin-integrable-system,hurtubise-surfaces,M,Bo1}
and Theorem \ref{thm-generalized-hitchin} below).
The $\bbc^\times$-invariance, of the Poisson structure induced by $\psi_S$ 
on the loop group, is equivalent to its {\em quadratic} nature.
This is an analogue of the following simple observation:
An algebraic Poisson structure on a vector space $V$ is
$\bbc^\times$-invariant, if and only if it has quadratic coefficients.
In other words, the weight $0$ subspace in
$\Sym^\bullet(V^*)\otimes \Wedge{2}V$ is $\Sym^2(V^*)\otimes \Wedge{2}V$. 
A $\bbc^\times$-invariant algebraic Poisson structure exists on the surface 
$T^*_\Sigma(D)$, with $D$ effective, if and only if the genus of 
$\Sigma$ is $\leq 1$. In contrast, a linear Poisson structure exists for 
curves of arbitrary genus. Consequently, the Hitchin systems were constructed
for a base curve of arbitrary genus, while the Sklyanin systems exist
only for a base curve of genus $0$ and $1$. 

\subsection{Overview of the relationship with the Dyanamical Yang-Baxter
Equations}
\label{sec-introduction-cdybe}

The work presented here is a generalization, 
to arbitrary complex reductive groups, 
of the picture presented above in Section \ref{sec-introduction-loop-group}. 
Half way through the project, we discovered (thanks to A. Polishchuk)
that our work is related to the Classical Dynamical Yang-Baxter Equations
(CDYBE), introduced by G. Felder \cite{felder}. 
Etingof and Varchenko interpreted solutions of 
the CDYBE, with spectral parameter, as Poisson structures on a loop 
groupoid \cite{etingof-varchenko}. 
We hope to convince the reader, that the algebro-geometric approach 
to the subject provides new insight to the geometry of the CDYBE. 
To the best of our knowledge, complete integrability of the Etingof-Varchenko 
Poisson structures was not known for a general reductive group. 
Even the symplectic leaves foliation was not understood, in the 
absence of an analogue of the Dressing Action. 

Following is an overview of the relationship with the CDYBE. 
Let $\gg$ be a complex simple Lie algebra. 
Choose an embedding $\gg\subset A$ as a Lie sub-algebra of an associative 
algebra ($A$ could be the universal enveloping algebra, or $\ggl(V)$,
for a faithful representation $V$).  
The Classical Yang-Baxter equation is the following algebraic equation 
\begin{equation}
\label{eq-cybe}
[r^{12}(z_1-z_2),r^{13}(z_1-z_3)]+
[r^{12}(z_1-z_2),r^{23}(z_2-z_3)]+
[r^{13}(z_1-z_3),r^{23}(z_2-z_3)] = 0,
\end{equation}
for a meromorphic function 
\[
r \ : \ \bbc \ \ \ \longrightarrow \ \ \ \gg\otimes \gg.
\]
Above, $z$ is the coordinate on $\bbc$.
The notation $r^{12}$ means $r\otimes 1$ as an element of 
$A\otimes A\otimes A$, while $r^{ij}$ indicates that $r$ acts on the $i$-th 
and $j$-th factors.  
Drinfeld discovered the geometric meaning of the CYBE. 
Solutions of (\ref{eq-cybe}), satisfying certain additional assumptions, 
define a natural Poisson-Lie structure on the loop group of $G$
\cite{drinfeld-icm,reiman-semenov}. 
Similarly, constant solutions define Poisson-Lie structures on $G$. 
Quantization of the CYBE led to the theory of Quantum groups. 

Belavin and Drinfeld classified 
solutions of (\ref{eq-cybe}), satisfying a non-degeneracy condition, 
for all simple Lie algebras \cite{belavin-drinfeld}. 
There are rational, trigonometric, and elliptic solutions. 
{\em Elliptic solutions exist only for type A}. 
These solutions are related to 
factorization of the loop group.
Such factorizations arise from a rigid 
principal bundle on an algebraic curve \cite{HM-sklyanin}. 
On a rational curve, every trivial principal bundle is rigid. 
When the curve is elliptic, the only example of a rigid principal
bundle, is the $PGL(N)$ bundle coming from 
a stable vector bundle with rank and degree co-prime
(see \cite{FM} Theorem 5.13).  

Let $\gt\subset \gg$ be a Cartan sub-algebra and $\{x_i\}$ a basis of $\gt$,
regarded as a set of coordinates on $\gt^*$. 
The Knizhnik-Zamolodchikov-Bernard (KZB) equations are the 
differential equations 
for flat sections of a connection over the moduli space of 
punctured elliptic curves. It was discovered by Bernard in the study of the 
Wess-Zumino-Witten conformal field theory on a torus \cite{bernard}. 
A coefficient in one of the two KZB equations, is a certain 
meromorphic elliptic function 
\[
r: \gt^* \times \bbc \ \ \ \longrightarrow \ \ \ 
(\gg\otimes\gg)^\gt.
\] 
The function $r$ satisfies the Classical Dynamical Yang-Baxter Equation 
\begin{eqnarray}
\nonumber
\sum_i x_i^{(1)}\frac{\partial r^{23}(x,z_2-z_3)}{\partial x_i} & + &
\sum_i x_i^{(2)}\frac{\partial r^{31}(x,z_3-z_1)}{\partial x_i} +
\\
\label{eq-cdybe}
\sum_i x_i^{(3)}\frac{\partial r^{12}(x,z_1-z_2)}{\partial x_i} & + &
{\rm CYBE} \ \ \ = \ \ \ 0,
\end{eqnarray}
where the term CYBE indicates the expression on the left hand side of 
(\ref{eq-cybe}). 
Moreover, let $r$ be an arbitrary meromorphic function as above, satisfying 
$r^{12}(x,z)=-r^{21}(x,-z)$. 
The KZB equation with coefficient $r$ is consistent, if and only if 
$r$ satisfies equation (\ref{eq-cdybe}) \cite{felder}.
We suppressed the dependence of $r$ 
on the modular parameter of the elliptic curve, as it may be fixed 
in (\ref{eq-cdybe}). 

Etingof and Varchenko found the geometric meaning of the CDYBE. 
A solution of (\ref{eq-cdybe}), satisfying certain additional conditions, 
corresponds to a natural Poisson structure on the groupoid 
\[
V \times LG \times V,
\]
where $V$ is an open subset of $\gt^*$ and $LG$ is the loop group of $G$. 
Naturality means, roughly,  compatibility with the group structure of $LG$. 
We will see that a point in $V$ should be considered as a principal
$G$-bundle, of trivial topological type,
on the elliptic curve $\Sigma$. More precisely, 
$\gt^*$ should be considered as a branched cover of the moduli of $G$-bundles,
with the affine Weyl group as Galois group \cite{looijenga}. 
The subset $V$ is the complement of the ramification locus. 
More general versions of the CDYBE are
obtained by replacing $\gt$ with certain subalgebras. 
Bundles of arbitrary topological types arise in this manner. 
Taking the zero sub-algebra, reduces equation (\ref{eq-cdybe}) to the CYBE.
(Solutions correspond to a rigid $G$-bundle). 

\medskip
A quotient of the groupoid $V \times LG \times V$, 
by the cartesian square of the affine Weyl group, 
may be regarded as a set of triples $(P_1,\varphi,P_2)$, 
where $P_1$ and $P_2$ are principal $G$-bundles, and $\varphi$ is a
generically invertible (formal or meromorphic) homomorphism between them. 
The automorphism group of a generic semi-stable $G$-bundle $P$, 
of trivial topological type, is isomorphic to the maximal torus 
$T\subset G$. Consequently, the product $T\times T$ acts on 
$V \times LG \times V$. The action on $\varphi\in LG$ is given by
\[
(t_1,t_2)\cdot \varphi \ \ \ = \ \ \ t_1\varphi t_2^{-1}. 
\]
We consider the subspace 
\begin{equation}
\label{eq-groupoid}
[V \times LG \times V](\Sigma), 
\end{equation}
where $\varphi$ is meromorphic over the elliptic curve $\Sigma$. 
The quotient of $[V \times LG \times V](\Sigma)$ by $T\times T$ 
admits a Poisson structure as well. It is a cover of an open set in 
the Hecke correspondence, 
with the cartesian square of the affine Weyl group as Galois group. 
The algebro-geometric construction of 
the Poisson structure on the Hecke correspondences is described in section 
\ref{sec-Hecke-correspondences}.

Our main objective, however, is the generalization of the 
Sklyanin integrable systems. These are obtained from 
the groupoid $[V \times LG \times V](\Sigma)$ as follows.
One considers the locus in $[V \times LG \times V](\Sigma)$, 
where the two principal 
bundles are equal, and takes the quotient by the diagonal of $T\times T$. 
The action of the diagonal corresponds to the conjugation action of the 
automorphism group of a bundle. The reduced Etingof-Varchenko Poisson structure
is invariant under the action of the affine Weyl group, and descends 
to the quotient. The quotient by the affine Weyl group is partially
compactified by the algebro-geometric moduli space $M_\Sigma(G,0)$.
This moduli space parametrizes isomorphism classes of pairs $(P,\varphi)$, 
consisting of a principal $G$-bundle $P$, of trivial topological type, 
and a meromorphic section of its adjoint group bundle $PG$. 
The Poisson structure extends ``across the walls'' to the whole of
$M_\Sigma(G,0)$. 
The moduli spaces $M_\Sigma(G,c)$, of arbitrary topological type $c$, 
are the main object of study of this paper. 

{\em Note: Luen-Chau Li has recently also obtained integrable
models associated with dynamical Poisson groupoids \cite{li}.}

\section{Statement of the main Theorem \ref{thm-main}}
\label{sec-main-results}
Let $G$ be a connected reductive group, $T$ a maximal torus,
$N(T)$ its normalizer and $W:=N(T)/T$ the Weyl group. 
Choose a non-degenerate symmetric invariant bilinear 
form on the Lie algebra $\gg$. 
The center of $G$ is denoted by $Z$ and its Lie algebra by $\gz$. 
Let $\Phi$ be the set of roots 
of $\gg$. We denote by $G((t))$ the formal loop group
and by $G[[t]]$ its natural subgroup of non-singular elements.
The group $G[[t]]\times G[[t]]$ acts on $G((t))$ via the left and right
action. 

Given a point $p\in \Sigma$, we denote by 
$\CompletedStructureSheaf{(p)}$ the formal
completion of the structure sheaf at $p$ and by 
$\CompletedFunctionField{(p)}$ its fraction field. 
Given a principal $G$-bundle $P$, 
we denote by $P(G)$ the adjoint group bundle and
by $P(G(\CompletedStructureSheaf{(p)}))$ and 
$P(G(\CompletedFunctionField{(p)}))$ the corresponding realizations of 
$G[[t]]$ and $G((t))$. 
A $G[[t]]\times G[[t]]$ orbit $\orbit$ in $G((t))$ determines a well-defined
$P(G(\CompletedStructureSheaf{(p)}))\times 
P(G(\CompletedStructureSheaf{(p)}))$ 
orbit $\orbit$ in $P(G(\CompletedFunctionField{(p)}))$. 

\begin{defi}
{\rm
Let $P$ be a principal $G$-bundle, and $\varphi$ a section of its adjoint 
bundle $P(G)$. The pair $(P,\varphi)$ is said to be 
{\em infinitesimally simple}, if 
$H^0(\gz(\varphi))=\gz$, where $\gz(\varphi)\subset P\gg$ is 
the subbundle  centralizing $\varphi$. Equivalently, the center $Z$
of $G$ has finite index in 
the subgroup $Z(\varphi)\subset H^0(P(G))$ centralizing $\varphi$. 
We say that the pair $(P,\varphi)$ is {\em simple} if $Z(\varphi)=Z$. 
}
\end{defi}

Given the following data:
\begin{enumerate}
\item 
an elliptic curve $\Sigma$, 
\item
a reductive group $G$, 
\item
a class $c\in \pi_1(G)$,
\item
a formal sum ${\displaystyle \orbit=\sum_{p\in\Sigma}\orbit_p}$ 
of $G[[t]]\times G[[t]]$-orbits $\orbit_p$ in  $G((t))$, all
but finitely many of which are equal to $G[[t]]$, 
\end{enumerate}
we denote by 
\[
M_\Sigma(G,c,\orbit)
\] 
the moduli space of isomorphism classes of simple pairs 
$(P,\varphi)$, consisting of a 
principal $G$-bundle $P$, of topological type $c$, and a meromorphic section 
$\varphi$ of $P(G)$ with singularities in $\orbit$. 
Two pairs $(P_1,\varphi_1)$ and $(P_2,\varphi_2)$ are isomorphic, if there is
an isomorphism $f: P_1 \IsomRightArrow P_2$ conjugating $\varphi_1$ to 
$\varphi_2$. Considering the case $P_1=P_2=P$, we see that the class of 
$(P,\varphi)$ determines $\varphi$ only up to 
conjugation by the group of global sections of $P(G)$. 
The existence of $M_\Sigma(G,c,\orbit)$ is briefly addressed in remark 
\ref{rem-existence}.

The quotient $G//G$, of $G$ by its adjoint action, is naturally
isomorphic to $T/W$. Consequently, $P(G)//P(G)$ is isomorphic to
the trivial $T/W$ bundle. 
Let $char: P(G) \rightarrow (T/W)_\Sigma$ be the quotient morphism. 
The data $\orbit$ determines a variety $(\overline{T/W})_\Sigma(\orbit)$ 
over $\Sigma$ 
(see (\ref{eq-overline-T-mod-W-of-orbit})). 
Over the open subset of $\Sigma$, where $\orbit_p=G[[t]]$, 
the variety $(\overline{T/W})_\Sigma(\orbit)$ is the trivial $T/W$-bundle. 
If a meromorphic map $s:\Sigma \rightarrow T/W$ comes from a pair
$(P,\varphi)$ with singularities in $\orbit$, 
then $s$ is a holomorphic section of $(\overline{T/W})_\Sigma(\orbit)$. 
Denote by $Char(\orbit)$ the space (Hilbert scheme) of global sections of 
$(\overline{T/W})_\Sigma(\orbit)$ and let 
\begin{equation}
\label{eq-char}
char \ : \ M_\Sigma(G,c,\orbit) \ \ \longrightarrow \ \ Char(\orbit)
\end{equation}
be the natural morphism. The main result of this paper is:  

\begin{thm}
\label{thm-main} 
The moduli space $M_\Sigma(G,c,\orbit)$ is smooth and its dimension 
is calculated in Corollary \ref{cor-dimension-of-a-symplectic-leaf}. 
$M_\Sigma(G,c,\orbit)$ admits a natural symplectic
structure. Under mild technical assumptions on the data $\orbit$, 
the invariant polynomial map (\ref{eq-char}) induces a
Lagrangian fibration.
The generic Lagrangian fiber is a smooth compact abelian variety. 
\end{thm}

\noindent
The Theorem is a summary of the results in 
Theorem \ref{thm-algebraic-2-form}, 
Corollary \ref{cor-dimension-of-a-symplectic-leaf}, 
Theorem \ref{thm-smoothness}, 
and
Theorem \ref{thm-general-complete-integrability}. 

\medskip
We will construct a Poisson structure on the infinite dimensional 
moduli space $M(G,c)$ of simple pairs $(P,\varphi)$, 
where $P$ is a principal $G$-bundle over $\Sigma$ of topological type $c$ 
and $\varphi$ is a meromorphic section of $P(G)$. 
We will show that $M(G,c)$ is filtered by finite dimensional 
Poisson subvarieties. The integrable systems in Theorem \ref{thm-main}
are the symplectic leaves of these Poisson subvarieties. 
There are various ways to bound the divisor in $\Sigma$ of singular points of 
$\varphi$. Different bounds give rise to different filtrations of $M(G,c)$. 
The filtration we use below will depend on a choice of a faithful 
representation $\rho:G\hookrightarrow GL(V)$. 
We denote by $E=P\times_\rho V$ the vector bundle associated to $P$. 
Let $D$ be a divisor on $\Sigma$. Denote by $M(G,c,\rho,D)$ the moduli space 
of simple
pairs $(P,\varphi)$, where $P$ is a principal $G$-bundle over $\Sigma$
and $\varphi$ is a meromorphic section of the adjoint group bundle $P(G)$,
which corresponds to a holomorphic section of $End(E)(D)$.

\begin{thm}
\label{thm-infinite-dimensional-Poisson-moduli}
$M(G,c)$ is endowed with a natural Poisson structure. 
The finite dimensional subvarieties 
$M(G,c,\rho,D)$ are Poisson subvarieties.
The connected components of $M(G,c,\orbit)$ in Theorem
\ref{thm-main} are symplectic leaves of $M(G,c,\rho,D)$,
where the singularity divisor $D:=D(\orbit,\rho)$ is given in
(\ref{eq-polar-divisor}).
\end{thm}

The Theorem is a summary of 
Theorem \ref{thm-Lambda-is-a-Poisson-structure} and
Lemma \ref{lemma-symplectic-leaves-are-loci-with-fixed-singularity-type}.

It is instructive to compare Theorem \ref{thm-main} to the 
analogous Theorem \ref{thm-generalized-hitchin}
for the generalized Hitchin system. 
We will relate the two theorems by a degeneration, 
in a particularly interesting case
(see Lemma \ref{lemma-example-deforms-to-the-calogero-moser-of-type-A-n}). 
We view $\gg[[t]]^*$ as an ind-variety 
$\cup_{d=1}^\infty \left\{\gg[t]/(t^d)\right\}^*$. 
Each of the finite dimensional subvarieties in the filtration is the dual of 
a Lie-algebra. Let $\CoadjointOrbit:=\sum_{p\in \Sigma}\CoadjointOrbit_p$
be a formal sum of (finite-dimensional) coadjoint orbits of $\gg[[t]]^*$, 
all but finitely many of which are the zero orbit. 
Given a point $p$ and a principal bundle $P$, let 
$P(\gg(\CompletedStructureSheaf{(p)}))$ denote the formal stalk at $p$ 
of the adjoint Lie algebra bundle. 
$P(\gg(\CompletedStructureSheaf{(p)}))$ is isomorphic to $\gg[[t]]$.
Composing the residue with the 
invariant bilinear form on $\gg$, we get an isomorphism 
\[
P(\gg(\CompletedStructureSheaf{(p)}))^* \ \ \ \ \cong \ \ \ \ 
[P\gg\otimes T^*\Sigma](\infty p) \ / \ [P\gg\otimes T^*\Sigma]
\]
between the dual of the Lie algebra and the space of polar tails of
meromorphic sections of $P\gg\otimes T^*\Sigma$. Denote by
\[
\Higgs_\Sigma(G,c,\CoadjointOrbit)
\]
the moduli space of stable pairs $(P,\varphi)$, consisting off
a principal $G$-bundle $P$ and a meromorphic section
$\varphi$ of $P\gg\otimes T^*\Sigma$ with poles in $\CoadjointOrbit$. 
The data $\CoadjointOrbit$ determines a variety 
$(\gt/W)(\CoadjointOrbit)$. Over the open subset of $\Sigma$, where
$\CoadjointOrbit_p=0$, the variety $(\gt/W)(\CoadjointOrbit)$ is the trivial
$\gt/W$-bundle. One has an invariant polynomial morphism 
\begin{equation}
\label{eq-hitchin-map}
h \ : \ \Higgs_\Sigma(G,c,\CoadjointOrbit) \ \rightarrow \ 
H^0(\Sigma,(\gt/W)(\CoadjointOrbit)).
\end{equation}

\begin{thm}
\label{thm-generalized-hitchin}
\cite{Bo1,hitchin-integrable-system,M,markman-sw}
Let $\Sigma$ be a smooth, compact, and connected algebraic curve
of arbitrary genus. 
The moduli space $\Higgs_\Sigma(G,c,\CoadjointOrbit)$ admits
an algebraic symplectic structure. Under mild technical assumptions on 
the orbits $\CoadjointOrbit_p$ in 
$\CoadjointOrbit$, the morphism (\ref{eq-hitchin-map}) 
is a Lagrangian fibration. 
\end{thm}

\noindent
The technical assumption is satisfied, in particular, 
if each of the orbits $\CoadjointOrbit_p$ in 
$\CoadjointOrbit$ is regular or closed \cite{markman-sw}. 

\begin{rem}
\label{rem-existence}
{\rm 
The existence of the moduli spaces $M(G,c,\orbit)$, of simple pairs, 
follows from that of $M(G,c,\rho,D)$. 
The existence of a coarse moduli space $M(G,c,\rho,D)$, 
as a (non-Hausdorff) complex analytic space, can be shown, 
for example, by a variant 
of the arguments proving Theorem 4.3 in \cite{ramanathan}.
We expect that $M(G,c,\rho,D)$ is an algebraic space, in analogy to the
case of simple sheaves \cite{altman-kleiman}. 
The similar moduli spaces of semi-stable meromorphic Higgs pairs were 
constructed by Simpson \cite{simpson}. 
Ramanathan's definition of stability and semi-stability for principal 
bundles \cite{ramanathan}, has a straightforward analogue for the pairs 
parametrized by $M(G,c,\rho,D)$. 
We expect that the coarse moduli space of 
semi-stable pairs can be constructed, as a separated algebraic scheme,
by means of Geometric Invariant Theory. Stable pairs are infinitesimally 
simple (an analogue of Proposition 3.2 in \cite{ramanathan}). 

The group $GL(n)$ embeds in its Lie-algebra. Consequently, 
a pair $(P,\varphi)$ in $M(G,c,\rho,D)$ can be considered 
as a Higgs bundles with poles along $D$, when $G=GL(n)$ and
$\rho$ is the standard representation. 
In that case, the moduli space $M(G,c,\rho,D)^{ss}$ of semi-stable pairs,
is a Zariski open subset of the moduli space of 
semi-stable Higgs bundles constructed by
Nitsure and Simpson \cite{nitsure,simpson}. 

The construction of the moduli space of semi-stable pairs, 
for a general reductive $G$, definitely merits  
further consideration. We do not address it here, as the study of the
Poisson geometry is infinitesimal in nature. 
The only exception is the complete integrability of $M(G,c,\orbit)$. 
However, the functorial results of \cite{donagi-gaitsgory} lead to an 
alternative proof of the existence of a coarse moduli space 
$M(G,c,\orbit)_{reg}$, as an algebraic space. 
$M(G,c,\orbit)_{reg}$ is the open subset of $M(G,c,\orbit)$ introduced in
Definition \ref{moduli-of-simple-and-regular-pairs}.  
It is a relative moduli space of a particular type of $T$-torsors, 
on curves on a fixed algebraic variety $X(\orbit)$
(see Section \ref{sec-invariant-polynomials}). 
The existence of $M(G,c,\orbit)_{reg}$ 
follows from the existence of Hilbert schemes and 
of relative Picards \cite{altman-kleiman}.
}
\end{rem}

\section{Orbits}
\label{sec-orbits}

In Section \ref{sec-Iwahori} we recall a few basic facts about the
orbits of $G[[t]]\times G[[t]]$ in $G((t))$. 
In section \ref{sec-topological-invariant} we consider a 
discrete topological invariant, which marks connected components of
$M_\Sigma(G,c,\orbit)$, when the group $G$ is not simply connected.
In Section 
\ref{sec-invariant-polynomials} we define the variety $(T/W)_\Sigma(\orbit)$.
The space $Char(\orbit)$, of global sections of 
$(T/W)_\Sigma(\orbit)\rightarrow\Sigma$,
is the target space for the invariant polynomial map (\ref{eq-char}). 
The symplectic geometry of $(T/W)_\Sigma(\orbit)$ is studied in Section
\ref{sec-2-forms}. Sections \ref{sec-invariant-polynomials} and 
\ref{sec-2-forms} will not be needed until Section 
\ref{sec-complete-integrability}.

\subsection{Iwahori's Theorem}
\label{sec-Iwahori}
Iwahori's Theorem, for reductive groups over local fields such as $\bbc((t))$,
asserts the equality:
\[
G((t)) \ \ = \ \ G[[t]]\cdot T((t))\cdot G[[t]].
\]
Consequently, 
$G[[t]]\times G[[t]]$ orbits in $G((t))$ are in one-to-one correspondence
with orbits of $N(T[[t]])\times T[[t]]$ in $T((t))$,
where $N(T[[t]])$ acts by conjugation and $T[[t]]$ by multiplication. 
The latter orbits are parametrized by 
$W$-orbits in the lattice $Ch(T)^*$ of co-characters of the maximal torus $T$. 
An element $\varphi$ of $T((t))$ determines the co-character
$\Ord(\bullet,\varphi)$, 
which associates to $\lambda\in Ch(T)$ the order 
of the formal function $\lambda(\varphi)\in\bbc((t))$ at $t=0$. 

Given a faithful representation $\rho:G\hookrightarrow GL(V)$, denote 
its set of weights by $A(\rho)$. $A(\rho)$ is a finite $W$-invariant
subset of $Ch(T)$. 
Given a $W$-orbit $\orbit$  in $Ch(T)^*$, set
\[
d_+(\orbit,\rho) \  \ := \ \ \max\{\varphi(w) \ \mid \ w\in A(\rho) \ 
\mbox{and} \ \varphi\in \orbit\}.
\]
Let $d_-(\orbit,\rho)$ be the minimum of the above set. 
If $G$ is a simple group, then we have the equality 
$d_-(\orbit,\rho)=-d_+(\orbit,\rho)$. 
Define the ``polar'' divisor of the 
singularity data ${\displaystyle \orbit=\sum_{p\in\Sigma}\orbit_p}$ 
by
\begin{equation}
\label{eq-polar-divisor}
D := D(\orbit,\rho) \ := -\ \sum_{p\in \Sigma} d_-(\orbit_p,\rho)\cdot p
\end{equation}
If $(P,\varphi)$ is a pair in $M(G,c,\orbit)$ and $E$ is the vector bundle 
associated to $P$ via $\rho$, then $\rho(\varphi)$ is a holomorphic section of
$\End(E)(D)$. 


\subsection{A topological invariant}
\label{sec-topological-invariant}

We introduce next a discrete topological invariant of a pair 
$(P,\varphi)$ in $M_\Sigma(G,c,\orbit)$. 
Let $\Sigma^0$ be the complement in $\Sigma$ of the support of the
singularity data $\orbit$. Consider the restriction 
$P(G\restricted{)}{\Sigma^0}$ of the adjoint bundle to $\Sigma^0$. 
Fix a point $p_0\in \Sigma^0$. 
The identity section of $P(G\restricted{)}{\Sigma^0}$
maps $\pi_1(\Sigma^0,p_0)$ injectively into 
$\pi_1(p_0,P(G\restricted{)}{\Sigma^0})$. 
The image of $\pi_1(\Sigma^0,p_0)$ is a normal subgroup.
This can be seen by the fact, that 
$P\times G$ is a principal $G$-bundle over $P(G)$ and 
$P\times \{1\}$ maps onto the identity section. 
Consequently, 
the fundamental group of $P(G\restricted{)}{\Sigma^0}$ is naturally 
isomorphic to the product $\pi_1(\Sigma^0,p_0)\times \pi_1(G)$. 
Using this decomposition, the section $\varphi$ induces a homomorphism
\begin{equation}
\label{eq-topological-invariant}
\tau_\varphi \ \ : \ \ \pi_1(\Sigma^0,p_0) \ \ \ \longrightarrow \ \ \ 
\pi_1(G). 
\end{equation}

The fundamental group $\pi_1(G)$ is isomorphic to a quotient of 
$Ch(T)^*$ by a $W$-invariant sublattice
\begin{equation}
\label{eq-homomorphism-from-co-characters-to-pi-1-G}
Ch(T)^* \ \ \longrightarrow \ \ \pi_1(G).
\end{equation}
The kernel is the co-character lattice of the simply connected cover
of the semi-simple derived group $[G,G]$. The kernel is also identified as the 
co-weight lattice of $[\gg,\gg]$. 
The homomorphism $\tau_\varphi$ maps a loop around a point $p$, 
to the image via (\ref{eq-homomorphism-from-co-characters-to-pi-1-G})
of a co-character $a_p$ in the $W$-orbit $\orbit_p$. 
The homomorphism $\tau_\varphi$ corresponds to a cohomology class in
$H^1(\Sigma^0,\pi_1(G))$, since the group $\pi_1(G)$ is abelian. 
If $(P_1,\varphi_1)$ and $(P_2,\varphi_2)$ are two pairs in 
$M_\Sigma(G,c,\orbit)$, then the difference
$\tau_{\varphi_1}-\tau_{\varphi_2}\in H^1(\Sigma^0,\pi_1(G))$
is the restriction of a global class in $H^1(\Sigma,\pi_1(G))$.

If $G$ is semisimple, then 
the class $\tau_\varphi$ determines a  $\widetilde{T}$-local system
over $\Sigma^0$, where $\widetilde{T}$ is a maximal torus of the universal 
cover of $G$. Partial compactifications, of the corresponding
principal $\widetilde{T}$-bundle, will play 
an important role in the study of the fibers of the
invariant polynomial map (\ref{eq-char}) (see section
\ref{sec-spectral-curves-and-group-isogenies}).

We will need, in addition, the identification of $\tau$ as a connecting
homomorphism. 
Let $G^{sc}$ be the universal cover of $G$. We get the short exact 
sequence of groups
\[
0 \rightarrow \pi_1(G) \rightarrow G^{sc} \rightarrow G \rightarrow 0.
\]
The group $G$ acts on $G^{sc}$ by conjugation. We denote by 
$P(G^{sc})$ the associated bundle. It fits in a short exact sequence 
\[
0 \rightarrow \pi_1(G)_\Sigma \rightarrow P(G^{sc}) \rightarrow P(G) 
\rightarrow 0.
\]
The corresponding sequence, of sheaves of analytic sections, is exact 
in the analytic topology. 
We get a long exact sequence of sheaf cohomologies.
\begin{eqnarray}
\nonumber
0 & \longrightarrow & H^0_{an}(\Sigma^0,\pi_1(G)) \ \rightarrow \ 
H^0_{an}(\Sigma^0,P(G^{sc}))  \ \rightarrow \ 
H^0_{an}(\Sigma^0,P(G)) 
\\
\label{eq-long-exact-seq-of-analytic-sheaf-coho}
& \LongRightArrowOf{\delta_{an}} & H^1_{an}(\Sigma^0,\pi_1(G)) \ 
\rightarrow \ H^1_{an}(\Sigma^0,P(G^{sc})).
\end{eqnarray}
The class $\tau_\varphi$ is equal to $\delta_{an}(\varphi)$. When
$\Sigma^0$ is affine, 
the connecting homomorphism $\delta_{an}$ is surjective, 
as the last cohomology set is trivial.

\begin{rem}
\label{rem-topological-non-emptiness-condition} 
{\rm
The singularity data $\orbit$ should satisfy the following 
necessary topological condition, for the moduli space 
$M_\Sigma(G,c,\orbit)$ not to be empty. 
For each singularity point $p\in \Sigma$, choose a co-character $a_p$ 
in the $W$-orbit determined by $\orbit_p$. Then the sum 
$\sum_{p\in \Sigma}a_p$ projects to the identity 
in $\pi_1(G)$. The image is independent of the choices, 
due to the $W$-invariance of the projection 
(\ref{eq-homomorphism-from-co-characters-to-pi-1-G}). 
See Remark \ref{rem-non-emptiness-condition-for-Hecke} 
for a proof of a more general condition.

The above topological condition admits a more refined algebraic version. 
Compare with condition 
(\ref{eq-linear-equivalence-of-zero-of-det-to-poles-of-det}). 

}
\end{rem}

\subsection{Invariant polynomials}
\label{sec-invariant-polynomials}

Given a pair $(P,\varphi)$ with singularities in $\orbit$, we get a 
meromorphic section of $P(G)//P(G)$. The latter is isomorphic to the trivial 
$(T/W)$-bundle $(T/W)_\Sigma$ over $\Sigma$. 
We define below a birational model $(\overline{T/W})(\orbit)$ of 
$(T/W)_\Sigma$, 
together with a morphism 
\[
\pi \ : \ (\overline{T/W})(\orbit) \ \longrightarrow \ \Sigma.
\] 
If a meromorphic map $s:\Sigma \rightarrow T/W$ comes from 
a pair $(P,\varphi)$ with singularities in $\orbit$, then
$s$ is a holomorphic section of $(\overline{T/W})(\orbit)$. 

%

We construct first a Zariski open subset $(T/W)(\orbit)$ of 
$(\overline{T/W})(\orbit)$. Its complement has codimension $\ge 2$. 
If a pair $(P,\varphi)$ has singularities in $\orbit$
and $\varphi$ is strictly-semi-simple and split in the sense below,
at each of its singularity points, 
then $char(P,\varphi)$ is a section of $(T/W)(\orbit)$. 
(Call an element of $G((t))$ strictly-semi-simple and split, 
if it can be conjugated into $T((t))$ via an element of $G[[t]]$).
We define $(T/W)(\orbit)$ first as an orbifold (a scheme 
with quotient singularities). We will show that it is Hausdorff
(a separated quasi-projective variety) by considering the local toric model of
its $W$-cover. The construction of $(T/W)(\orbit)$, as an orbifold, 
is local over the curve $\Sigma$, 
so we may assume that $\orbit_p=G[[t]]$ for all $p\neq p_0$. Let
$A\subset Ch(T)^*$ be the $W$-orbit corresponding to 
$\orbit_{p_0}$. Given $a\in A$, define $X(a)$ to be the
principal $T$-bundle characterized by
\[
X(a)\times_T \lambda \ \ \ = \ \ \ 
\StructureSheaf{\Sigma}(a(\lambda)\cdot p_0)
\setminus \mbox{(zero section)}, \ \ \ \mbox{for all} \ 
\lambda\in Ch(T).
\]
We have a natural trivialization of $X(a)$ over $\Sigma\setminus \{p_0\}$.
We define $X(\orbit)$ to be a scheme with an open cover
by the $X(a)$, $a\in A$. Given $a$, $b$ in $A$, 
we glue $X(a)$ and $X(b)$ by the trivializations (leaving their fibers
over $p_0$ disjoint). The Weyl group acts on the disjoint union
of the $X(a)$ permuting the connected components. The $W$-action descends to
$X(\orbit)$. We set
\[
(T/W)(\orbit) \ \ := \ \ X(\orbit)/W.
\]

The smooth variety $X(\orbit)$ has a local toric model. 
We replace the base curve by the affine line $\Sigma:=\bbc$ and assume that 
$\orbit_p=G[[t]]$ for all $p\neq 0$. 
Set $\widetilde{T}:= \bbc^\times \times T$. 
Let $\widetilde{A}\subset Ch(\widetilde{T})^*$ be the set of
graphs of elements of $A$. Given $a:\bbc^\times \rightarrow T$, its
graph is
\[
\tilde{a}:= (id,a) : \bbc^\times \ \rightarrow \ [\bbc^\times \times T].
\]
The vectors in $\widetilde{A}$ generate a strongly convex rational polyhedral 
cone $\sigma$ in $Ch(\widetilde{T})^*\otimes_\bbz\bbr$. 
The fan $\overline{F}$, consisting of all
the faces of $\sigma$, corresponds to a singular affine toric variety
$X(\overline{F})$. Our toric variety $X(\orbit)$ is the complement in
$X(\overline{F})$ of all the $\widetilde{T}$-orbits of codimension $\geq 2$. 
In other words, we work with the fan $F$, 
consisting of the faces of $\sigma$ of dimension zero and $1$. 
Each element of $\widetilde{A}$ generates a one-dimensional face of $\sigma$.
Then $X(\orbit)$ is the smooth toric variety $X(F)$ and 
\[ 
(T/W)(\orbit) \ \ \cong \ \ X(F)/W. 
\]
The fiber of $(T/W)(\orbit)$ over $0$ is isomorphic to $T/W_a$, where
$W_a$ is the stabilizer of a co-character $a\in A$.

Let $\overline{X}(\orbit)$ be the partial toroidal compactification of
$X(\orbit)$, given locally by the toric model $X(\overline{F})$. 
The values of the invariant polynomial morphism (\ref{eq-char}) 
are global sections of the $W$-quotient 
\begin{equation}
\label{eq-overline-T-mod-W-of-orbit}
(\overline{T/W})(\orbit) \ \ := \ \ \overline{X}(\orbit)/W.
\end{equation}

Further information on the structure of $(T/W)(\orbit)$ 
and $(\overline{T/W})(\orbit)$ is provided in Lemma 
\ref{lemma-chi-bar-is-an-imersion-away-from-the-ramification-locus}.

\begin{example}
{\rm
Let $G=SL(2)$ and let $\orbit$ correspond to the $W$-orbit $\{ka,-ka\}$, 
where the co-character $a$ is the generator of $Ch(T)^*$ and $k$ is
a positive integer. Then $Ch(\widetilde{T})^*$
is spanned by the graph $e$ of the trivial co-character and by $a$.
The fan $F$ consists of faces of dimension $\leq 1$ of the cone $\sigma$ 
generated by $e+ka$ and $e-ka$. 
The dual cone $\sigma^\vee$ is generated by $ke^*+a^*$ and
$ke^*-a^*$. The semigroup of lattice points in $\sigma^\vee$ 
is generated by the three points $e^*$, $ke^*+a^*$, and $ke^*-a^*$.
$X(F)$ is the complement of the isolated singular point 
$(x,w,z)=(0,0,0)$, in the affine surface 
\[
{\rm Spec}(\bbc[x,x^ky,x^ky^{-1}]) \ \ \cong \ \ 
{\rm Spec}(\bbc[x,w,z]/(x^{2k}=wz). 
\]
It maps to $\bba^1={\rm Spec}(\bbc[x])$. The fiber over $0$ is $wz=0$, 
which is the disjoint union of two $\bbc^\times$'s. 
The fiber over $c$ is the smooth affine conic $wz=c^{2k}$, which is 
isomorphic to $\bbc^\times$. $W$ acts on $X(F)$ by interchanging 
$z$ and $w$. The quotient $X(F)/W$ is the complement of the origin in the 
surface
\[
{\rm Spec}(\bbc[x,w+z,wz]/(x^{2k}=wz) \ \ \cong \ \ 
{\rm Spec}(\bbc[x,w+z]) \ \ \cong \ \ \bba^2.
\]
}
\end{example}

\subsection{$2$-forms}
\label{sec-2-forms}

We investigate in this section 
the symplectic geometry of the variety $X(\orbit)$. Lemma 
\ref{lemma-t-valued-2-form} shows that 
the space of global sections of $\pi : (T/W)(\orbit)\rightarrow \Sigma$
is the base of a prym integrable system of the type 
studied in \cite{HM-pryms}. 
The results of this section will not be used in the rest of the paper.
We include them nevertheless, since they lead to a very short construction
of the integrable systems of the current paper, 
using the results of \cite{HM-pryms}. 
The results of \cite{HM-pryms} require the additional smoothness assumption
for the $W$-covers in $X(\orbit)$ of sections of $(T/W)(\orbit)$.

Given a character $\lambda\in Ch(T)$, we get the bundle homomorphism
\[
(\lambda\times \bbi) \ : \ T \times \Sigma \ \ \ \longrightarrow \ \ \ 
\bbc^\times \times \Sigma.
\]
Let $\omega$ be the $\bbc^\times$-invariant 
two-form $\frac{dt}{t}\wedge dz$ on the surface $\bbc^\times \times \Sigma$
and $\omega_\lambda$ its pullback to $X:=T \times \Sigma$. 
We get a $W$-equivariant {\em linear} homomorphism 
\begin{eqnarray*}
\Omega \ : \ Ch(T) & \longrightarrow &
H^0(X,\Wedge{2}T^*_X).
\\
\hspace{5ex} \lambda & \mapsto & \omega_\lambda.
\end{eqnarray*}
Consequently, we can regard $\Omega$ as a $W$-invariant section of 
$H^0(X,(\Wedge{2}T^*_X)\otimes_\bbc \gt)$. It is, in fact, a section
of the subbundle $\pi^*T_\Sigma\otimes T^*_{X/\Sigma}\otimes \gt$, where
$T^*_{X/\Sigma}$ is the relative cotangent bundle.

\begin{lemma}
\label{lemma-t-valued-2-form}
\begin{enumerate}
\item
The $\gt$-valued $2$-form $\Omega$ extends to an algebraic and $W$-invariant 
$\gt$-valued $2$-form $\Omega$ on $X(\orbit)$. 
\item
The contraction with $\Omega$ induces a surjective homomorphism
from $\Wedge{2}T_{X(\orbit)}$ to the trivial $\gt$-bundle. 
\item
The wedge square $\Omega\wedge\Omega$ vanishes as a section of
$\Wedge{4}T^*_X\otimes_\bbc \gt\otimes\gt$. 
\end{enumerate}
\end{lemma}

\noindent
{\bf Proof:}
The lemma is clear in case $X(\orbit)= T\times \Sigma$. 
In general, we do have a birational isomorphism between $X(\orbit)$
and $T\times \Sigma$, so the question is local in the neighborhood
of each fiber of $X(\orbit)$ over a singularity point $p\in \Sigma$. 
Given a divisor $D:=\sum_{i=1}^d m_i p_i$, $p_i \in \Sigma$,
the $2$-form $\omega$, being $\bbc^\times$-invariant, 
extends to a non-degenerate $2$-form on the 
complement of the zero section in $\StructureSheaf{\Sigma}(D)$. 
Consequently,  we get a {\em linear} homomorphism
\[
\Omega \ : \ Ch(T) \ \ \ \longrightarrow \ \ \
H^0(X(a),\Wedge{2}T^*_{X(a)}),
\]
for each $a\in Ch(T)^*$.
It is easy to see that contraction with $\Omega$ induces a surjective 
homomorphism from $\Wedge{2}T_{X(a)}$ to the trivial $\gt$-bundle over $X(a)$. 
Moreover, $\Omega$ glues to a global $\gt$-valued $2$-form on
$X(\orbit)$, which is $W$-invariant. 
\EndProof

Let $s:\Sigma \hookrightarrow (T/W)(\orbit)$ be a section, $C$ 
the inverse image of $s(\Sigma)$ in $X(\orbit)$, and  
$q:C\rightarrow s(\Sigma)$ the natural morphism.
Assume that $C$ is smooth. 
$W$ acts on the vector bundle $[q_*(T^*C)\otimes_\bbc \gt]$ via the 
diagonal action 
(via automorphisms of $C$ on the first factor and via the standard 
reflection representation on $\gt$). Denote by $N_{s(\Sigma)}$ the normal 
bundle of $s(\Sigma)$ in $(T/W)(\orbit)$. Contraction with $\Omega$ 
induces an injective homomorphism from
$N_{C/X(\orbit)}$ into $T^*_C\otimes_{\bbc}\gt$. The homomorphism is not
surjective along the ramification points of $q$.
Nevertheless, Lemma \ref{lemma-t-valued-2-form} implies:

\begin{cor}
\label{cor-normal-bundle-in-terms-of-cotangent-bundle-of-cameral-cover}
The following three vector bundles are isomorphic:
\[
[q_*(T^*C)\otimes_\bbc \gt]^W \ \ \cong \ \
[q_*N_C]^W \ \ \cong \ \ N_{s(\Sigma)}.
\]
\end{cor}

\noindent
{\bf Proof:}
This was proven, in a slightly more general setting, in
Lemma 5.6 in \cite{HM-pryms}.
\EndProof

\medskip
Note that the global sections of $[q_*(T^*C)\otimes_\bbc \gt]^W$ 
are the global one-forms on the generalized prym
(\ref{eq-generalized-prym}).
Hence, the dimension of the generalized prym 
is equal to the dimension of (the component in) 
the space of global sections of $(T/W)(\orbit)$. 
It was shown in \cite{HM-pryms} that the union $M$, 
of the generalized pryms parametrized by sections of $(T/W)(\orbit)$, 
admits a canonical symplectic structure. Moreover, the 
generalized pryms are Lagrangian subvarieties of $M$. 

\section{Infinitesimal study of $M(G,c,\orbit)$}

\subsection{Infinitesimal deformations}
\label{sec-infinitesimal-deformations}
We identify first the infinitesimal deformations of a meromorphic section 
$\varphi$ in the adjoint bundle $P(G)$ of a fixed principal bundle $P$. 
We further restrict attention to deformations 
with singularities is a fixed formal sum 
$\orbit=\sum_{p\in \Sigma} \orbit_p$. These infinitesimal deformations 
are determined by the global sections of a
vector bundle $ad(P,\varphi)$ defined below 
(Lemma \ref{lemma-sections-of-ad-P-varphi-are-infinitesimal-def}).

Let $\Delta_\varphi\subset P\gg\oplus P\gg$ be the subbundle
\begin{equation}
\label{eq-Delta-phi}
\Delta_\varphi \ := \ \{(a,b) \ : \ a+Ad_\varphi(b)=0\}.
\end{equation}
Denote the quotient by
\[
ad(P,\varphi) \ := \ [P\gg\oplus P\gg]/\Delta_\varphi.
\]
We get the short exact sequence
\begin{equation}
\label{eq-short-exact-seq-defining-ad-P-phi}
0\rightarrow \Delta_\varphi \rightarrow [P\gg\oplus P\gg] \LongRightArrowOf{q} 
ad(P,\varphi) \rightarrow 0.
\end{equation}
Notice that the intersection of $\Delta_\varphi$ with the sum of the centers
$P\gz\oplus P\gz$ is given by
\begin{equation}
\label{eq-intersection-of-Delta-phi-with-sum-of-centers}
\Delta_\varphi\cap [P\gz\oplus P\gz] \ \ \ = \ \ \ 
\{(a,b)\in P\gz\oplus P\gz \ : \ a+b=0\}
\end{equation}
and is isomorphic to $P\gz$. In particular, the image of
$P\gz\oplus P\gz$ in $ad(P,\varphi)$ is isomorphic to $P\gz$. 
Let $e_i:P\gg\hookrightarrow P\gg\oplus P\gg$, $i=1,2$, be the two inclusions
and denote by 
\[
R_\varphi, \ L_\varphi : P\gg \ \longrightarrow \ ad(P,\varphi)
\]
the two compositions $R_\varphi:=q\circ e_1$ and $L_\varphi:=q\circ e_2$. 
We set
\[
ad_\varphi \ := \  L_\varphi-R_\varphi.
\]

\begin{example}
{\rm
If $\varphi=\bbi$ is the identity section of $P(G)$, then $ad(P,\bbi)=P\gg$.
If $\varphi$ is a holomorphic section of $P(G)$, then 
$R_\varphi$ and $L_\varphi$ are isomorphisms. If $\varphi$ is a holomorphic
section of the center of $P(G)$, then $R_\varphi=L_\varphi$. 
}
\end{example}

The following Lemma explains the geometric significance of $ad(P,\varphi)$
in terms of a faithful representation $\rho$ and the associated vector bundle 
$E$. Assume that $\varphi$ is a holomorphic section of $\End(E)(D)$
and let 
\[
R_{\rho(\varphi)}, L_{\rho(\varphi)} \ : \ \End(E) \ \rightarrow \ 
\End(E)(D)
\]
be the homomorphisms given by left and right multiplication by 
$\rho(\varphi)$. 

\begin{lemma}
\label{lemma-sections-of-ad-P-varphi-are-infinitesimal-def}
\begin{enumerate}
\item
The image of the composition homomorphism
\[
[P\gg\oplus P\gg] \ \LongRightArrowOf{(\rho,\rho)} \ 
\End(E)\oplus \End(E) \ 
\LongRightArrowOf{(R_{\rho(\varphi)},L_{\rho(\varphi)})} \ 
\End(E)(D)
\]
is isomorphic to $ad(P,\varphi)$. Moreover, under this identification,
$R_\varphi$ and $L_\varphi$ correspond to 
$R_{\rho(\varphi)}$ and $L_{\rho(\varphi)}$. 
\item
The global sections of $ad(P,\varphi)$ are in one-to-one correspondence
with the infinitesimal deformations of $\varphi$ as a section of 
the sheaf of double orbits
\[
\rho(P(G))\cdot \varphi \cdot \rho(P(G)) \ \ \subset \ \ \End(E)(D).
\]
\end{enumerate}
\end{lemma}

\noindent
{\bf Proof:} The kernel of the above composition is clearly 
isomorphic to $\Delta_\varphi$.
\EndProof

\begin{cor}
The infinitesimal deformations of the pair 
$(P,\varphi)$ in $M(G,c,\orbit)$ 
are naturally identified by the first hyper-cohomology
of the complex (in degrees $0$ and $1$)
\begin{equation}
\label{eq-tangent-complex-for-symplectic-leaf}
P\gg \ \ \LongRightArrowOf{ad_\varphi} \ \  ad(P,\varphi). 
\end{equation}
\end{cor}

\noindent
{\bf Proof:}
A repetition of the arguments in
\cite{Bo1,biswas-ramanan,M}.
\EndProof

\medskip
We define next the notion of a family of pairs with singularities in
$\orbit$. Let $X$  a scheme of finite type over $\bbc$. 
Assume given a pair $(\P,\Phi)$, of a principal bundle $\P$ over 
$\Sigma\times X$ and a meromorphic section $\Phi$ of $\P(G)$.

\begin{defi}
\label{def-constant-singularity-data}
{\rm
We say that $\Phi$ has a {\em locally constant singularity type}, 
if each point $(p,x)\in \Sigma\times X$ admits
1) open neighborhoods (analytic, \'{e}tale, or even formal)
$U\subset \Sigma$ of $p$ and $V\subset X$ of $x$, 
2) a regular 
section (trivialization) $\eta : U\times V \rightarrow \P$, 
3) a regular section $f$ of $\P(G)$ over $U\times V$, and
4) a rational morphism $\varphi : U \rightarrow G$, satisfying
\[
f\cdot\Phi\cdot\eta \ \ \ = \ \ \ \eta \cdot (\varphi\circ \pi_\Sigma). 
\]
In other words, the trivialization $\eta$ takes the meromorphic section
$f\cdot \Phi$ to a rational morphism from $U\times V$ to $G$, 
which is the pullback of a rational morphism from $U$ to $G$. 
We say that $\Phi$ has {\em singularities in $\orbit$}, 
if, furthermore, $\varphi$ has singularities in $\orbit$.  
}
\end{defi}

If $\Phi$ has singularities in $\orbit$, then 
the quotient sheaf $ad(\P,\Phi)$ of $\P\gg\oplus \P\gg$,
by the relative analogue $\Delta_\Phi$ of (\ref{eq-Delta-phi}), 
is easily seen to be a locally free sheaf on $\Sigma\times X$. 

\begin{rem}
{\rm
In a relative setting, the Kodaira-Spencer map can be 
identified cohomologically by the homomorphism
$\kappa^0$ defined below in (\ref{eq-kappa-i}). Let $X$ be a scheme, $P$ 
a principal $G$-bundle over $X\times \Sigma$, and $\varphi$ 
a meromorphic section of $P(G)$ with a constant singularity data $\orbit$. 
Let $\pi_X:X\times \Sigma\rightarrow X$ be the projection. 
Denote by $\A(P)$ the Atiyah algebra of $G$-invariant vector fields on $P$
with symbols in $\pi_X^*TX$. $\A(P)$ is a sheaf on $X\times \Sigma$, 
which fits in the short exact sequence
of the symbol map $\sigma$
\[
0 \rightarrow \ P\gg \ \longrightarrow \
\A(P) \ \LongRightArrowOf{\sigma} \ \pi^*TX \ \rightarrow 0.
\]

Away from its singularities, we regard $\varphi$ as an automorphism of $P$. 
Given a vector field $\xi$ on $P$,
we get the new vector field $ad_\varphi(\xi)$, which takes a 
local function $f$ on $P$ to $\varphi_*(\xi(f))-\xi\varphi_*(f)$
(and $\varphi_*$ is pullback by $\varphi^{-1}$). This is simply the 
infinitesimal version of conjugating by $\varphi$ a one-parameter family of 
automorphisms of $P$. 
The $G$ action commutes with the $P(G)$-action. Hence, 
away from its singularities, $ad_\varphi$ 
takes sections of $\A(P)$ to sections of $\A(P)$. Moreover, its image 
consists of vertical tangent vectors, i.e., its image is in $P\gg$. 
Globally, we get a homomorphism 
\[
ad_\varphi \ : \ \A(P) \ \ \ \longrightarrow \ \ \ ad(P,\varphi),
\]
which restricts to the homomorphism 
(\ref{eq-tangent-complex-for-symplectic-leaf}) on $P\gg$.
We get the following diagram of complexes on $X\times \Sigma$.
\[
[0\rightarrow \pi_X^*TX] \ \LongLeftArrowOf{(0,\sigma)} \ 
[P\gg \rightarrow \A(P)] \ \LongRightArrowOf{(id,ad_\varphi)} \ 
[P\gg \LongRightArrowOf{ad_\varphi}ad(P,\varphi)].
\]
The left homomorphism of complexes $(0,\sigma)$ is a quasi-isomorphism. 
Hence, we get on the level of higher direct images a homomorphism
\begin{equation}
\label{eq-kappa-i}
\kappa^i \ : \ TX\otimes H^i(\Sigma,\StructureSheaf{\Sigma}) 
\ \ \ \longrightarrow \ \ \ 
R^{i+1}_{\pi_{X,*}}[P\gg \LongRightArrowOf{ad_\varphi}ad(P,\varphi)].
\end{equation}
The relative Kodaira-Spencer map is $\kappa^0$. 
}
\end{rem}

The universal properties of $M(G,c,\orbit)$ imply, that 
the Zariski tangent space at a point $(P,\varphi)$ is given by 
the first hypercohomology of 
(\ref{eq-tangent-complex-for-symplectic-leaf})
\[
T_{(P,\varphi)}M(G,c,\orbit) \ \ \cong \ \ H^1(P\gg \ \ 
\LongRightArrowOf{ad_\varphi} \ \  ad(P,\varphi)).
\]
More generally, if the pair $(P,\varphi)$ is infinitesimally simple,
then the first hypercohomology of 
(\ref{eq-tangent-complex-for-symplectic-leaf})
computes the tangent space of the moduli stack. 

\subsection{The symplectic structure}
\label{sec-symplectic-structure}
The cotangent space at a point $(P,\varphi)$ in $M(G,c,\orbit)$ 
is given by the first hypercohomology of the complex 
\begin{equation}
\label{eq-cotangent-complex-of-a-symplectic-leaf}
ad(P,\varphi)^* \ \ \LongRightArrowOf{ad_\varphi^*} \ \ P\gg^*.
\end{equation}
(here we used a trivialization of the canonical line bundle of $\Sigma$). 
We have a natural homomorphism $\Psi$ from the cotangent 
complex to the tangent complex: 
\[
\begin{array}{ccc}
P\gg & \LongRightArrowOf{ad_\varphi} & ad(P,\varphi)
\\
\Psi_0 \ \uparrow  \hspace{3ex} & & \hspace{3ex} \uparrow \ \Psi_1
\\
ad(P,\varphi)^* & \LongRightArrowOf{ad^*_{\varphi}} & P\gg^*.
\end{array}
\]
In degree $1$ the homomorphism $\Psi_1$
is the composition of the isomorphism $P\gg^*\cong P\gg$ with 
the homomorphism $L_\varphi$ from $P\gg$ to $ad(P,\varphi)$. 
It is easier to define the dual of the homomorphism $\Psi_0$ in degree
zero. The negative $-\Psi_0^*$ of the dual is the composition of the 
isomorphism $P\gg^*\cong P\gg$ with the {\it right} multiplication 
homomorphism $R_\varphi$ from $P\gg$ to $ad(P,\varphi)$. 

We need to check that $\Psi$ is a homomorphism of complexes.
It suffices to check it away from the singular divisor of $\varphi$.
Over this open subset of $\Sigma$, the invariant bilinear form on $\gg$ 
induces an isomorphism
$$ad(P,\varphi)^*\ \ \cong \ \  ad(P,\varphi^{-1}).$$
The composition, of $\Psi_0$ with the generic isomorphism
$ad(P,\varphi)^* \cong  ad(P,\varphi^{-1})$, is equal to 
the negative of left 
multiplication by $\varphi$. Similarly, $ad_\varphi^*$
becomes $-ad_\varphi$. The commutativity
$ad_\varphi\circ \Psi_0=\Psi_1\circ ad_\varphi^*$ follows. 

\begin{rem}
\label{rem-interchanging-roles-of-left-and-right}
{\rm
If we interchange the roles of
left and right multiplication in the above construction, we get  
another homomorphism of complexes.
The two homomorphism of complexes are homotopic. The homotopy is given 
by the complex homomorphism $h$, of degree $-1$, from 
(\ref{eq-cotangent-complex-of-a-symplectic-leaf}) to 
(\ref{eq-tangent-complex-for-symplectic-leaf}), which sends $P\gg^*$
isomorphically to $P\gg$. Indeed, 
$(\psi_0,\psi_1)-(h\circ ad^*_\varphi,ad_\varphi\circ h)$
is obtained from $(\psi_0,\psi_1)$ by interchanging the roles of $L_\varphi$
and $R_\varphi$.
}
\end{rem}

We get an induced homomorphism 
\[
\Psi \ : \ T^*_{(P,\varphi)}M(G,c,\orbit) \ \ \longrightarrow \ \ 
T_{(P,\varphi)}M(G,c,\orbit)
\]
on the level of first hypercohomology. We will see in 
Theorem \ref{thm-smoothness}, that the moduli space 
$M(G,c,\orbit)$ is smooth. 
The above construction can be carried out in the relative setting
of families of pairs. 
The relative construction gives rise to a global homomorphism from the 
cotangent to the tangent bundles:
\begin{equation}
\label{eq-algebraic-2-form}
\Psi \ : \ T^*M(G,c,\orbit) \ \ \longrightarrow \ \ TM(G,c,\orbit).
\end{equation}

\begin{thm}
\label{thm-algebraic-2-form}
The homomorphism $\Psi$, given in 
(\ref{eq-algebraic-2-form}), is an anti-self-dual isomorphism. 
Its inverse defines a non-degenerate closed algebraic $2$-form on 
$M(G,c,\orbit)$.
\end{thm}


\noindent
{\bf Proof:} (of Theorem \ref{thm-algebraic-2-form})
Note that $-\Psi^*$ is the
complex homomorphism obtained from $\Psi$ by interchanging the roles
of $R_\varphi$ and $L_\varphi$. Consequently, the complex
homomorphisms $\Psi$ and $-\Psi^*$ are homotopic 
(Remark \ref{rem-interchanging-roles-of-left-and-right}). 
It follows that the 
homomorphism  $\Psi$ is anti-self-dual, on the level of first cohomology. 
We prove next that $\Psi$ is a quasi-isomorphism. 
We have the short exact sequence of complexes
\[
\begin{array}{ccccc}
\left[
\begin{array}{c}
ad(P,\varphi)^* \\ \downarrow \ ad_\varphi^* \\
P\gg^* 
\end{array}
\right]
& 
\begin{array}{c}
{\buildrel{-R_\varphi^*}\over {\longrightarrow}} \\ \\
{\buildrel{L_\varphi}\over {\longrightarrow}} 
\end{array}
& 
\left[
\begin{array}{c}
P\gg \\ \downarrow \ ad_\varphi \\ ad(P,\varphi) 
\end{array}
\right]
& 
\begin{array}{c}
{\buildrel{q_0}\over {\longrightarrow}} \\ \\
{\buildrel{q_1}\over {\longrightarrow}} 
\end{array}
& 
\left[
\begin{array}{c}
P\gg/Im(R_\varphi^*) \cr \downarrow \ \overline{ad}_\varphi \cr 
ad(P,\varphi)/Im(L_\varphi)
\end{array} 
\right]
\end{array}
\]
where the first column is (\ref{eq-cotangent-complex-of-a-symplectic-leaf}) 
and the second is (\ref{eq-tangent-complex-for-symplectic-leaf}).
The image of the homomorphism 
$(L_\varphi,ad_\varphi) : [P\gg\oplus P\gg] \rightarrow ad(P,\varphi)$ 
is equal to the image of $(L_\varphi,R_\varphi)$. The latter is surjective,
by definition of $ad(P,\varphi)$. It follows that 
the vertical right hand homomorphism $\overline{ad}_\varphi$ is surjective. 
It is clearly a homomorphism
between two sheaves with zero-dimensional support
and the two sheaves have the same length. Hence, $\overline{ad}_\varphi$
is an isomorphism. 

The closedness of the $2$-form is proven in Section \ref{sec-Jacobi-identity}.
\EndProof

\subsection{The dimension formula}
We will need to calculate the degree of $ad(P,\varphi)$. 
We first observe that the sequence
(\ref{eq-short-exact-seq-defining-ad-P-phi})
is self-dual. Denote by $\kappa : P\gg \rightarrow P\gg^*$ the isomorphism
induced by the invariant bilinear pairing. 

\begin{lemma}
\label{lemma-ad-P-phi-is-dual-to-Delta-phi}
\begin{enumerate}
\item
\label{lemma-item-ad-P-phi-is-dual-to-Delta-phi}
We have the commutative diagram, whose vertical homomorphisms are
isomorphisms:
\[
\begin{array}{ccccccccc}
0 & \rightarrow & \Delta_\varphi & \rightarrow & P\gg\oplus P\gg & 
\rightarrow  & ad(P,\varphi)  & \rightarrow  & 0
\\
& & \downarrow & & \hspace{2ex}\cong \ \downarrow 
\left[
\begin{array}{cc}
0       & \kappa \\
-\kappa & 0
\end{array}
\right] 
& & 
\downarrow
\\
0  & \rightarrow  & ad(P,\varphi)^*  & \rightarrow  & [P\gg\oplus P\gg]^*  & 
\rightarrow & \Delta_\varphi^*  & \rightarrow  & 0.
\end{array}
\]
\item
\label{lemma-item-bundles-isomorphic-to-ad-P-phi}
The following vector bundles are canonically isomorphic: 
\end{enumerate}
\[
ad(P,\varphi)
 \cong 
\Delta_\varphi^*
 \cong  
[P\gg\cap Ad_\varphi(P\gg)]^*
 \cong 
[Ad_{\varphi^{-1}}(P\gg)\cap P\gg]^* 
 \cong  
\Delta_{\varphi^{-1}}^*
 \cong  
ad(P,\varphi^{-1}). 
\]
\end{lemma}

\noindent
{\bf Proof:}
\ref{lemma-item-ad-P-phi-is-dual-to-Delta-phi}) The kernel $K$ of
\[
0\rightarrow K \rightarrow [P\gg\oplus P\gg]^* \rightarrow \Delta_\varphi^*
\rightarrow 0
\]
is equal to
\[
K \ = \ \{(\alpha,\beta)\in  P\gg^*\oplus P\gg^* \ \mid \
\beta=Ad_\varphi^*(\alpha)\}. 
\]
We see, that 
$
{\displaystyle
\left[
\begin{array}{cc}
0       & \kappa \\
-\kappa & 0
\end{array}
\right] 
}$ 
maps $\Delta_\varphi$ isomorphically 
onto $K$. By its definition, the vector bundle  $ad(P,\varphi)^*$
is equal to $K$. In particular, 
$ad(P,\varphi)$ is isomorphic to $\Delta_\varphi^*$. 

\ref{lemma-item-bundles-isomorphic-to-ad-P-phi}) 
The first isomorphism was proven in part
\ref{lemma-item-ad-P-phi-is-dual-to-Delta-phi}.
The image of the homomorphism 
\begin{eqnarray*}
[Ad_\varphi(P\gg)\cap P\gg] & \longrightarrow & 
P\gg\oplus P\gg 
\\
a & \mapsto & (a, -Ad_{\varphi^{-1}}(a))
\end{eqnarray*}
is equal to $\Delta_\varphi$. 
The second isomorphism follows. 
The isomorphism
\[
Ad_{\varphi^{-1}} \ : \ [Ad_\varphi(P\gg)\cap P\gg] \ \ \LongIsomRightArrow
[P\gg \cap Ad_{\varphi^{-1}}(P\gg)] 
\]
implies the third isomorphism.
\EndProof

Given a  $W$-invariant set of characters $B\subset Ch(T)$ and a 
$G[[t]]\times G[[t]]$ orbit $\orbit$ in $G((t))$, we define the sum
\[
\sum_{\alpha\in B} \Ord(\alpha,\orbit)
\]
by choosing a representative of $\orbit$ lying in $T((t))\cap \orbit$. The sum
is independent off the choice because  $B$ is $W$-invariant. 

\begin{lemma}
\label{lemma-degree-of-ad-P-phi}
The degree of $ad(P,\varphi)$ is given by the following formula:
\[
\chi(ad(P,\varphi)) \ \ = \ \
\sum_{p\in \Sigma} \sum_{\alpha\in \Phi} \max\{0,-\Ord_p(\alpha,\orbit_p)\},
\]
where $\orbit_p$ is the 
$P(G(\CompletedStructureSheaf{(p)}))\times 
P(G(\CompletedStructureSheaf{(p)}))$
orbit of $\varphi$.
\end{lemma}

\noindent
{\bf Proof:} 
It suffices to calculate the degree of $[Ad_{\varphi^{-1}}(P\gg)\cap P\gg]$,
which is isomorphic to $ad(P,\varphi)^*$ by Lemma
\ref{lemma-ad-P-phi-is-dual-to-Delta-phi}.
The calculation of the subsheaf $[Ad_{\varphi^{-1}}(P\gg)\cap P\gg]$
of $P\gg$ is local and 
we may assume that $\varphi$ is an element of $T((t))$. 
Then the sheaf $[Ad_\varphi(P\gg)\cap P\gg]$ decomposes as a direct sum 
of line-bundles 
\begin{equation}
\label{eq-stabilizer-of-co-character}
[Ad_{\varphi^{-1}}(P\gg)\cap P\gg]   =   
\{a\in P\gg \ : \ Ad_\varphi(a)\in P\gg \}   =   
\oplus_{\alpha\in \Phi} t^{\max\{0,-\Ord(\alpha,\varphi)\}} P\gg_\alpha.
\end{equation}
\EndProof

We set
\[
\gamma(\orbit_p) \ \ := \ \ 
\sum_{\alpha\in \Phi} \max\{0,-\Ord_p(\alpha,\orbit_p)\}.
\]
The formula for $\gamma(\orbit_p)$ can be simplified once we choose a dominant 
co-character representing $\orbit_p$ 
(see Equation (\ref{eq-gamma-is-dot-product-with-delta})).
Let $\Lambda_r$ be the root lattice, 
$\Delta:=\{\alpha_1, \dots, \alpha_r\}$ its basis of simple roots, and 
$\Delta^\vee:=\{\alpha_1^\vee, \dots, \alpha_r^\vee\}$ the basis of co-roots, 
where $\alpha_i^\vee := \frac{2\alpha_i}{(\alpha_i,\alpha_i)}$. 
Denote by $\{\lambda_1, \dots, \lambda_r\}$ the basis 
of fundamental dominant weights of the weight lattice $\Lambda$. 
Recall that $\{\lambda_1, \dots, \lambda_r\}$ is the basis dual to 
$\Delta^\vee$ with respect to the inner product: 
$(\lambda_i,\alpha_j^\vee)=\delta_{i,j}$. 
We have the natural inclusions $\Lambda_r\subset \Lambda$,
$\Lambda_r\subset Ch(T)$, and a homomorphism $Ch(T)\rightarrow \Lambda$. 
The latter homomorphism is injective if $G$ is semi-simple 
and surjective if $G$ is simply connected. Let $\Phi_+$ be the set of 
positive roots. Set
\[
\delta \ \ := \ \ \frac{1}{2}\sum_{\alpha\in \Phi_+}\alpha. 
\]
Then the following relation holds (\cite{humphreys} Section 13.3 Lemma A
page 70):
\[
\delta \ \ = \ \  \sum_{i=1}^r \lambda_i.
\]
The fundamental Weyl chamber $C(\Delta)\subset \Lambda\otimes_\bbz \bbr$
is the cone generated by the fundamental weights. 
Given a co-character $a\in Ch(T)^*$, we denote by $\bar{a}$ its image
in $\Lambda_r^*$. We say that $a$ is {\em dominant}, 
if  $\bar{a}=(b,\bullet)$ for a dominant weight $b$.
Every $W$-orbit in $Ch(T)^*$ contains a dominant element.

Choose a co-character $a$ of $\orbit_p$ with $-a$ dominant. Then
$\alpha(\bar{a})\leq 0$, for every positive root 
$\alpha$. We get that
\begin{equation}
\label{eq-gamma-is-dot-product-with-delta}
\gamma(\orbit_p) \ \ = \ \ -2\bar{a}(\delta).
\end{equation}
Note that the integer $\gamma(\orbit_p)$ is the dimension of the orbit,
parametrized by $\orbit_p$, in the loop grassmannian $G((t))/G[[t]]$
(see \cite{lusztig}).

\begin{cor}
\label{cor-dimension-of-a-symplectic-leaf}
Let $(P,\varphi)$ be a simple pair. The dimension of the Zariski tangent 
space of $M(G,c,\orbit)$ at $(P,\varphi)$ 
is given by 
\[
\dim(M(G,c,\orbit)) \ \ = \ \ 2\dim(\gz)
+ \sum_{p\in \Sigma} \gamma(\orbit_p).
\]
\end{cor}

Smoothness $M(G,c,\orbit)$ is established in 
Theorem \ref{thm-smoothness}. 
Consequently, the corollary determines the dimension of $M(G,c,\orbit)$. 
The similarity, between the dimension of $M(G,c,\orbit)$ and the dimension of 
products of orbits in the loop Grassmannian, is explained in Section
\ref{sec-dimensions-of-Hecke-and-M}. 
The dimension $\gamma(\orbit_p)$, of such an orbit,
is even if the group $G$ is simply connected. 
In general, $\gamma(\orbit_p)$ may not be even (see Equation 
(\ref{eq-dimension-in-GL-case})).
Theorem \ref{thm-algebraic-2-form}
implies, that the sum $\sum_{p\in \Sigma} \gamma(\orbit_p)$ is even.
This follows independently from the topological non-emptiness 
condition in Remark \ref{rem-topological-non-emptiness-condition} 

\medskip
\noindent
{\bf Proof of the Corollary:}
We follow a straightforward calculation. 
The pair $(P,\varphi)$ is infinitesimally simple. 
It follows that $H^0$ of the complex
(\ref{eq-tangent-complex-for-symplectic-leaf})
is isomorphic to $\gz$. The second cohomology of 
(\ref{eq-tangent-complex-for-symplectic-leaf}) 
is dual to $H^0$ of the complex 
(\ref{eq-cotangent-complex-of-a-symplectic-leaf}).
The latter is isomorphic to $\gz^*$ (identify 
$ad(P,\varphi)^*$ with $\Delta_\varphi$ via Lemma
\ref{lemma-ad-P-phi-is-dual-to-Delta-phi} and use the equality
(\ref{eq-intersection-of-Delta-phi-with-sum-of-centers})).

The dimension of the Zariski tangent space of 
$M(G,c,\orbit)$ is equal to the 
dimension of $H^1$ of the complex
(\ref{eq-tangent-complex-for-symplectic-leaf}). It is given by
\[
\dim H^1 \ = \ \dim H^0+\dim H^2 - [\chi(P\gg)-\chi(ad(P,\varphi)]
\ = \ 2\dim(\gz) + \chi(ad(P,\varphi)). 
\]
\EndProof

\begin{example}
\label{example-dimension-of-leaf-in-GL-case}
{\rm
Let $G=GL(n)$, $D$  an effective divisor, $E$ a rank $n$ vector bundle and
$\varphi,\psi : E \rightarrow E(D)$ two homomorphisms.
$G[[t]]\times G[[t]]$ orbits in $G((t))$ are in one-to-one correspondence with 
$G[[t]]$-orbits in the loop Grassmannian $G((t))/G[[t]]$. 
Hence, $\varphi$ and $\psi$ belong to the same 
orbit $\orbit=\sum_{p\in \Sigma}\orbit_p$, if and only if their co-kernels
are isomorphic as $\StructureSheaf{\Sigma}$-modules. 
Denote by $\epsilon_i$, $i=1, \dots, n$, the weights of the standard 
representation. The roots of $\ggl_n$ and $\gsl_n$ are 
$\{\epsilon_i-\epsilon_j\}_{i,j=1, \ i\neq j}^n$,
the fundamental weights are 
$\lambda_i=\sum_{j=1}^i\epsilon_j-\frac{i}{n}(\epsilon_1+\cdots +\epsilon_n)$, 
$1\leq i\leq n-1$, and 
$
2\delta=\sum_{i=1}^n(n+1-2i)\cdot\epsilon_i.
$
Assume that the co-character $\Ord_p(\bullet,\varphi)$ is in the $W$-orbit
of $-(b_p,\bullet)$, for some dominant weight 
$b_p=\sum_{i=1}^{n-1}\beta_i(p)\lambda_i$,
$\beta_i(p)\geq 0$. Then
\begin{equation}
\label{eq-dimension-in-GL-case}
\dim M(GL(n),c,\orbit) \ = \ 
2+ \sum_{p\in \Sigma} (b_p,2\delta) \ = \ 
2+\sum_{p\in \Sigma}\sum_{i=1}^{n-1}\beta_i(p)[i\cdot(n-i)]. 
\end{equation}
See equation 
(\ref{eq-dimension-of-symplectic-leaf-in-GL-case-and-simple-zeros-poles})
for a version of this formula in a special case.

}
\end{example}

\begin{example}
\label{example-partial-flag-varieties-of-maximal-parabolics}
{\rm
Certain partial flag varieties are realized as orbits in the loop
Grassmannian $G((t))/G[[t]]$. If $\orbit_p$ corresponds
to such an orbit, then the contribution 
$\gamma(\orbit_p)$, of the point $p$ to the dimension of
$M(G,c,\orbit)$, is equal to the dimension of the flag variety
(see the formula in Corollary \ref{cor-dimension-of-a-symplectic-leaf}). 
We compile here a list of such cases. The geometry
of the corresponding moduli spaces $M(G,c,\orbit)$ is studied in
Section \ref{sec-products-of-partial-flag-varieties}, for groups $G$
of type $A$ and $D$. 
Let $P$ be a parabolic subgroup of $G$ and suppose that 
$a\in Ch(T)^*$ satisfies 
\begin{eqnarray*}
a(\alpha)& \geq 0 & \mbox{for every root} \ \alpha \ {\mbox of} \ P \ \ 
\mbox{and}
\\
a(\alpha) & =-1 & 
\mbox{for every root} \ \alpha \ {\mbox of}  \ G, \ 
\mbox{which is not a root of} \ P
\end{eqnarray*}
(it follows that $a$ vanishes on the roots of the Levi factor of $P$).
Let $\orbit_p$ be an orbit with co-characters in the $W$-orbit of $a$.
Then  $\gamma(\orbit_p)$ is equal to the dimension of $G/P$. 
Examples of such pairs $(P,a)$ can be obtained as follows. 
Assume that $G$ is of adjoint type. Then the root lattice coincides with the
character lattice. 
Choose a basis $\Delta$ of simple roots with dual basis
$\Delta^*\subset Ch(T)^*$. Let $\alpha\in \Delta$ be a simple root
and $a:=\alpha^*\in \Delta^*$ its dual. Let $P_\alpha$ be the maximal 
parabolic corresponding to $\alpha$. The set of roots of $P_\alpha$ is 
$\Phi_{P_\alpha}:=\{\beta\in\Phi \ : \ \alpha^*(\beta)\geq 0\}$. 
The set of roots $\Phi_{U_\alpha}$, of the unipotent radical of $P_\alpha$, 
is defined using a strict inequality. We get the equality
\begin{equation}
\label{eq-sum-over-roots-in-opposite-unipotent-radical}
\gamma(\orbit_p) \ \ := \ \ 
\sum_{\beta\in \Phi}\max\{0,-\alpha^*(\beta)\}
\ \ = \ 
\sum_{\beta\in -\Phi_{U_\alpha}}\alpha^*(\beta)
\end{equation}
(the sum of the coefficients of $\alpha$ in the roots of the 
unipotent radical of the opposite parabolic). 
Assume that $\alpha$ satisfies the following condition:
\begin{equation}
\label{eq-condition-for-a-maximal-parabolic}
\alpha^*(\beta) \ \ \in \ \ \{0,1,-1\}, \ \ \mbox{for every root} \ \beta.
\end{equation}
Then the sum (\ref{eq-sum-over-roots-in-opposite-unipotent-radical})
is precisely the dimension of $G/P_\alpha$.

The co-character $\alpha^*$ determines the orbit 
$\orbit:=G[[t]]\alpha^*G[[t]]/G[[t]]$ in the loop grassmannian 
$G((t))/G[[t]]$. 
The subgroup of $G((t))$, stabilizing the image of $G[[t]]$ in 
$G((t))/G[[t]]$,
is $G[[t]]$ itself. Hence, the subgroup of $G((t))$, stabilizing 
the image of $\alpha^*$ in $G((t))/G[[t]]$, is $Ad^{-1}_{\alpha^*}(G[[t]])$.
Consequently, the subgroup $G[[t]]_{\alpha^*}$ of $G[[t]]$, stabilizing 
the image of $\alpha^*$ in $G((t))/G[[t]]$, is 
$G[[t]]\cap Ad^{-1}_{\alpha^*}(G[[t]])$. We see, that evaluation at $t=0$
induces a surjective homomorphism from $G[[t]]_{\alpha^*}$ onto $P_\alpha$.
We have a natural morphism from the orbit $\orbit$
to the flag variety $G/P_\alpha$. The orbit $\orbit$ 
is a $tG[[t]]/(tG[[t]]\cap G[[t]]_{\alpha^*})$-bundle over $G/P_\alpha$. 
The Lie-algebra of  
$G[[t]]_{\alpha^*}$ is given by (\ref{eq-stabilizer-of-co-character}),
with $\varphi$ replaced by $\alpha^*$. 
Condition (\ref{eq-condition-for-a-maximal-parabolic}) is equivalent to the
condition, that the morphism $\orbit\rightarrow G/P_\alpha$ is an 
isomorphism. 

When $G$ is of type $A_n$, then all maximal parabolics satisfy 
condition (\ref{eq-condition-for-a-maximal-parabolic}). 
The corresponding flag varieties  
are Grassmannians of subspaces in the standard representation of $SL(n+1)$. 
When $G$ is of type $D_n$, then only the three roots at the ends of the 
Dynkin diagram satisfy condition 
(\ref{eq-condition-for-a-maximal-parabolic}). 
The corresponding flag varieties  of $SO(2n)$ parametrize 
isotropic subspaces 
in the standard $2n$-dimensional representation. 
The subspaces are of dimensions $1$ (the quadric in $\bbp^{2n-1}$)
or $n$ (two distinct components). 
When $G$ is of type $B_n$, then the only root satisfying condition
(\ref{eq-condition-for-a-maximal-parabolic}) is the long root 
at the edge of the Dynkin diagram. The corresponding flag variety is 
the $(2n-1)$-dimensional  quadric in $\bbp^{2n}$ of isotropic lines
in the standard representation of $SO(2n+1)$. 
When $G$ is of type $C_n$, then the only root satisfying condition
(\ref{eq-condition-for-a-maximal-parabolic}) is the long root. 
The corresponding flag variety is the Lagrangian grassmannian of $Sp(2n)$. 
When $G$ is of type $E_6$, the two simple roots $\alpha_1$ and $\alpha_6$
give rise to compact orbits in the loop Grassmannian.
When $G$ is $E_7$, then only one root $\alpha_7$ leads to a compact orbit. 
For $E_8$, $F_4$, and $G_2$, there isn't any maximal parabolic satisfying  
condition (\ref{eq-condition-for-a-maximal-parabolic})
(see the table of highest roots in \cite{humphreys} page 66 Table 2). 

Note, that all partial flag varieties of maximal parabolics, 
admit realizations as $G[[t]]$ orbits in $G((t^{1/d}))/G[[t]]$,
for $d$-th roots of $t$ of sufficiently high order.
}
\end{example}

\bigskip
{\bf A non-emptiness condition:}
The data $\orbit=\sum_{p\in \Sigma} \orbit_p$ needs to satisfy 
a general necessary condition, in order for  $M(G,c,\orbit)$ to be non-empty. 
Let $(P,\varphi)$ be a pair in $M(G,c,\orbit)$, 
$\rho:G\rightarrow GL(n)$  a representation, $E$ the associated vector bundle 
and $\rho(\varphi):E\rightarrow E$ the associated meromorphic endomorphism.
Then the zero and polar divisors of $\det(\rho(\varphi))$
must be linearly equivalent. In terms of the weights $A\subset Ch(T)$
of the representation, the condition becomes:
\begin{equation}
\label{eq-linear-equivalence-of-zero-of-det-to-poles-of-det}
\sum_{p\in \Sigma}\left[
\sum_{w\in A} \max\{0,-\Ord_p(w,\orbit_p)\}
\right]\cdot p
\ \ \ \sim \ \ \ 
\sum_{p\in \Sigma}\left[
\sum_{w\in A} \max\{0,\Ord_p(w,\orbit_p)\}
\right]\cdot p.
\end{equation}
(See Remark \ref{rem-topological-non-emptiness-condition} for a topological 
version of this condition and equation 
(\ref{eq-down-to-earth-linear-equivalence-of-zero-and-poles-of-det}) 
for a down to earth version of this condition in a special case). 


\subsection{Smoothness}
\label{sec-smoothness}
We use the $T^1$-lifting property of Ziv Ran to prove that 
$M(G,c,\orbit)$ is smooth \cite{kawamata-T1-lifting,ziv-ran-T1-lifting}. 

\begin{thm}
\label{thm-smoothness}
The moduli space $M_\Sigma(G,c,\orbit)$ is smooth (or empty). 
\end{thm}

\noindent
{\bf Proof:}
We follow Kawamata's proof of the
Mukai-Artamkin criterion (\cite{kawamata-T1-lifting} Theorem 3). 
Let $(Art)$ be the category of Artinian local $\bbc$ algebras with residue
field $\bbc$. Let $(Set)$ be the category of sets. 
Fix a simple pair $(P,\varphi)$ with singularities in $\orbit$. 
Consider the following deformation functor $D: (Art) \rightarrow (Set)$. 
$D(A)$ is the set of isomorphism classes of triples $(\P,\Phi,f_0)$,
where $\P$ is a principal $G$-bundle on $\Sigma\times {\rm Spec}(A)$, 
$\Phi$ is a
section of $\P(G)$ with singularities in $\orbit$
(Definition \ref{def-constant-singularity-data}), and $f_0$ is an isomorphism 
of $(\P,\Phi)\times_{{\rm Spec}(A)}{\rm Spec}(k)$ with $(P,\varphi)$. 

Set $A_n := \bbc[t]/(t^{n+1})$ and $B_n := A_n\otimes_\bbc A_1$. Let
$\beta_n : B_n \rightarrow A_n$ be the natural homomorphism.  
Set $\Sigma_n:=\Sigma\times {\rm Spec}(A_n)$. 
Given $(\P_n,\Phi_n,f_0)\in D(A_n)$, we denote 
by $(\P_{n-1},\Phi_{n-1},f_0)$ its restriction to ${\rm Spec}(A_{n-1})$.
Define
$
T^1((\P_n,\Phi_n)/A_n)
$
to be the set of isomorphism classes of quadruples $(Y_n,\Psi_n,h_0,h_n)$,
where $(Y_n,\Psi_n,h_0)$ is an object in $D(B_n)$ and $h_n$ is an
isomorphism between $(\P_n,\Phi_n)$ and 
the restriction of $(Y_n,\Psi_n)$ to $A_n$. We require also that $h_0$ is the 
composition of $f_0$ with the restriction of $h_n$. 
Then we have a natural identification 
\[
T^1((\P_n,\Phi_n)/A_n) \ \ \cong \ \ 
\bbh^1\left(
\Sigma_n,\left[\P_n\gg \LongRightArrowOf{ad\Phi_n} ad(\P_n,\Phi_n)\right]
\right).
\]
Denote the complex $[P\gg\rightarrow ad(P,\varphi)]$
by $C_\bullet$ and let $C_\bullet^n$ be the complex 
$[\P_n\gg \LongRightArrowOf{ad\Phi_n} ad(\P_n,\Phi_n)]$. 
We get a short exact sequence of complexes on $\Sigma_n$
\begin{equation}
\label{eq-short-exact-sequence-of-complexes-of-lifts}
0 \rightarrow C_\bullet \rightarrow C_\bullet^n 
\rightarrow C_\bullet^{n-1} \rightarrow 0. 
\end{equation}
Considering the long exact sequence of hyper-cohomologies, 
we see that the vector space $\bbh^2(\Sigma,C_\bullet)$
is the obstruction space $T^2$ to lifting an infinitesimal deformation of 
$(\P_{n-1},\Phi_{n-1})$ to an infinitesimal deformation of 
$(\P_n,\Phi_n)$. 

Theorem 1 in \cite{kawamata-T1-lifting} and the universal properties of 
$M(G,c,\orbit)$ imply, that smoothness of 
$M(G,c,\orbit)$ at $(P,\varphi)$ follows from the
$T^1$-lifting property established in the Lemma \ref{lemma-T1-lifting}.
This completes the proof of Theorem \ref{thm-smoothness}.
\EndProof

\medskip
Key to the proof of the $T^1$-lifting property for simple sheaves, 
is Mukai's observation, that obstruction classes must live in
the kernel of the trace map \cite{mukai-symplectic}. 
Next, we construct the analogous trace map for our pairs.
Let $\gz$ be the center of $\gg$ and $(\gz)_{\Sigma_n}$ be the 
trivial $\gz$-bundle over $\Sigma_n$. 
Let $\tr_0: \P_n\gg \rightarrow (\gz)_{\Sigma_n}$ 
be the projection with respect to the 
orthogonal decomposition, arising from the invariant bilinear form on $\gg$. 
Using the equality (\ref{eq-intersection-of-Delta-phi-with-sum-of-centers}), 
we get a commutative diagram defining the
homomorphism $\tr_1$
\[
\begin{array}{ccccc}
\Delta_{\Phi_n} & \rightarrow & \P_n\gg\oplus \P_n\gg & \rightarrow 
& ad(\P_n,\Phi_n)
\\
\downarrow & & \downarrow & & \hspace{3ex} \ \downarrow \ \tr_1
\\
(\gz)_{\Sigma_n} & \rightarrow & (\gz)_{\Sigma_n} \oplus (\gz)_{\Sigma_n}
& \rightarrow & (\gz)_{\Sigma_n},
\end{array}
\]
where the bottom right arrow is summation. 
We get the morphism of complexes $\tr := (\tr_0,\tr_1)$
\[
\tr \ : \ C^n_\bullet \ \ \ \longrightarrow \ \ \ 
\left[(\gz)_{\Sigma_n} \ \LongRightArrowOf{0} (\gz)_{\Sigma_n}\right].
\]
Denote the complex $[(\gz)_{\Sigma_n}\LongRightArrowOf{0}(\gz)_{\Sigma_n}]$
by $Z_\bullet^n$. Then $\bbh^2(Z_\bullet^n)$ is isomorphic to 
$\bbh^1((\gz)_{\Sigma_n})$ and hence to 
$\gz\otimes H^1(\Sigma_n,\StructureSheaf{\Sigma_n})$. 
If $(P,\varphi)$ is simple, then the morphism 
$\tr^i:\bbh^i(C^n_\bullet)\rightarrow \bbh^i(Z^n_\bullet)$
is an isomorphism for $i=0,2$. We have shown that, when $n=0$, 
in the proof of Corollary \ref{cor-dimension-of-a-symplectic-leaf}. 
It follows for every $n$ by Nakayama's Lemma. We get the isomorphism 
\[
\tr^2 \ : \ \bbh^2(C_\bullet^n) \ \ \ \longrightarrow \ \ \ 
\bbh^2(Z_\bullet^n) \ \ \ \cong \ \ \ 
\gz\otimes H^1(\StructureSheaf{\Sigma_n}).
\]

\begin{lemma}
\label{lemma-T1-lifting}
For any positive integer $n$ and any object
$(\P_n,\Phi_n,f_0)$ in $D(A_n)$, the restriction homomorphism
\[
r \ \ : \ \ T^1((\P_n,\Phi_n)/A_n) \ \ \ \longrightarrow \ \ \ 
T^1((\P_{n-1},\Phi_{n-1})/A_{n-1})
\]
is surjective. 
\end{lemma}

\noindent
{\bf Proof:}
Consider the commutative diagram with exact rows, where the upper row is
a segment of the long exact sequence arising from 
(\ref{eq-short-exact-sequence-of-complexes-of-lifts}). 
\[
\begin{array}{ccccc}
\bbh^1(C_\bullet^n) & \LongRightArrowOf{r} & \bbh^1(C_\bullet^{n-1}) &
\LongRightArrowOf{\alpha} & \bbh^2(C_\bullet)
\\
\tr^1 \ \downarrow \ \hspace{3ex} & & 
\hspace{3ex} \ \downarrow \ \tr^1 & & 
\hspace{3ex} \ \downarrow \ \tr^2
\\
\bbh^1(Z_\bullet^n) & \LongRightArrowOf{\beta} & \bbh^1(Z_\bullet^{n-1}) &
\LongRightArrowOf{\gamma} & \bbh^2(Z_\bullet)
\end{array}
\]
Now, $\beta$ is surjective, being equal to $\bbi\otimes \beta'$, where
$\bbi$ is the identity of $\gz$ and 
$\beta': H^1(\Sigma_n,\StructureSheaf{\Sigma_n})\rightarrow 
 H^1(\Sigma_{n-1},\StructureSheaf{\Sigma_{n-1}})$ 
is the natural homomorphism. Hence, $\gamma=0$. 
The injectivity of $\tr^2$ implies that $\alpha=0$. 
This establishes the surjectivity of $r$.
\EndProof

\section{The Jacobi Identity}
\label{sec-Jacobi-identity}

We complete the proof of Theorem \ref{thm-algebraic-2-form}
in this section. 
Consider the pairing
\begin{eqnarray}
\nonumber
\theta \ : \  
\Wedge{2}T^*M(G,c,\orbit) & \rightarrow & H^2(S,\omega_S)\cong \bbc
\\
\label{eq-pairing-theta}
\theta(\alpha,\beta)& = & \langle\Psi_G(\alpha),\beta\rangle,
\end{eqnarray}
where $\Psi_G$ is given in (\ref{eq-algebraic-2-form})
and $\langle\bullet,\bullet\rangle$ is the Serre Duality pairing. 
The bracket $\{f,g\}$ of two local functions on $M(G,c,\orbit)$ 
is defined by 
\[
\{f,g\} \ \ \ := \ \ \ \theta(df,dg).
\]
The statement of Theorem 
\ref{thm-algebraic-2-form}, 
that $\Psi_G$ is a Poisson structure, 
means that the structure sheaf of $M(G,c,\orbit)$ is a Lie algebra with 
respect to that bracket. One has to prove the Jacobi identity
\[
\{f,\{g,h\}\} + \{g,\{h,f\}\}+\{h,\{f,g\}\} \ \ \ = \ \ \ 0,
\]
for every three local functions on $M(G,c,\orbit)$. 
We provide an algebraic proof of the Jacobi identity
in this section.  
Note that an alternative proof follows from the relation with
the Dynamical Yang-Baxter equation
discussed in section \ref{sec-dynamical-yang-baxter}.

We reduce the proof of the Jacobi identity to the $GL(n)$ case, 
by the choice of the faithful representation $\rho: G \hookrightarrow GL(n)$. 
We would have liked to argue as follows. 
The moduli space $M(G,c,\orbit)$ is mapped 
into a Poisson moduli space $M(GL(n),\rho(c),D)$. The 
Jacobi identity for the Poisson structure on $M(GL(n),\rho(c),D)$
would be verified by realizing this moduli space as a moduli space of
sheaves on a Poisson surface. Moreover, 
$M(G,c,\orbit)$ is mapped into a single symplectic leaf 
$M(GL(n),\rho(c),\rho(\orbit))$ of $M(GL(n),\rho(c),D)$. 
Furthermore, the symplectic form on $M(GL(n),\rho(c),\rho(\orbit))$
restricts to the $2$-form defined on $M(G,c,\orbit)$  via the tensor $\Psi_G$ 
in Theorem \ref{thm-algebraic-2-form}. 

Essentially, we will carry out the above argument.
There are, however, two major technical difficulties. First, the $GL(n)$ 
pair $(E,\rho(\varphi))$, associated to an infinitesimally simple
$G$-pair $(P,\varphi)$, may no longer be simple. Second, the Jacobi identity
for moduli spaces of sheaves on a Poisson surface was proven, in general, 
only for locally free sheaves \cite{Bo2} (see also
\cite{ziv-ran,mukai-survey} for proofs in the case of locally free sheaves on 
a symplectic surface). 
However, the spectral construction associates 
to a $GL(n)$-pair a sheaf supported on a spectral curve in the Poisson surface
\[
S:=\PP[\StructureSheaf{\Sigma}(D)\oplus \StructureSheaf{\Sigma}].
\] 
Thus, we need to deal with torsion sheaves on $S$ with pure one-dimensional 
support. Moreover, the spectral curves arising from a general $G$ could be 
all singular and the eigensheaves on the spectral curves need not be 
line-bundles. 
In section \ref{sec-fourier-mukai}
we prove the Jacobi identity for {\em all smooth moduli spaces 
of simple sheaves on any Poisson surface}
(Theorem 
\ref{thm-Poisson-structure-on-moduli-of-sheaves-on-a-Poisson-surface}).
This is done by means of a Fourier-Mukai transform, which maps the moduli of
simple sheaves into a moduli space of locally free and simple sheaves. 

Let us summarize the various constructions that we apply. 
We start with a pair 
$(P,\varphi)$ in $M(G,c,\orbit)$. We use the representation 
$\rho:G\hookrightarrow GL(V)$
to get the $GL(V)$-pair $(P\times_\rho V,\rho(\varphi))$.
The spectral dictionary associates to the $GL(V)$-pair a sheaf $F$ on the
surface $S$ with pure one-dimensional support \cite{simpson}. 
We then apply a Fourier-Mukai functor $\Phi$ to obtain 
a locally free sheaf $\Phi(F)$. 
\begin{equation}
\label{eq-four-moduli-spaces}
\begin{array}{ccccccc}
M(G,c,\orbit) & \LongRightArrowOf{e} & M(GL(V),\rho(c),D) & \subset & M_\rho 
& \LongRightArrowOf{\Phi} & X,
\\
(P,\varphi) & \mapsto & (P\times_\rho V,\rho(\varphi)) & \mapsto &
F & \mapsto & \Phi(F).
\end{array}
\end{equation}

The proof of the Jacobi identity, for a general $G$, is based on a reduction 
to the $GL(n)$ case. This requires a separate proof of the $GL(n)$ 
case first. We chose to avoid duplications. 
In order to extract the proof of Theorem \ref{thm-algebraic-2-form} in 
the $GL(n)$ case first, the reader need only read Theorem 
\ref{thm-Poisson-structure-on-moduli-of-sheaves-on-a-Poisson-surface} 
in Section \ref{sec-fourier-mukai} and proceed to 
Lemma \ref{lemma-symplectic-leaves-are-loci-with-fixed-singularity-type}. 

\subsection{Restriction of Poisson structures}
\label{sec-induced-poisson-structures}

We will first explain the proof of Theorem 
\ref{thm-algebraic-2-form}, 
ignoring the two difficulties 
mentioned above. The changes required
will be addressed in the following sections starting with 
\ref{sec-fourier-mukai}. 
When $G=GL(n)$ and $\rho$ is the standard representation, 
the moduli space $M(G,c,\rho,D)$ is naturally
a Zariski open subset of a moduli space $M_\rho$ of simple sheaves on
the Poisson surface $S$. 
The surface $S$ has a 
$\bbc^\times$-invariant Poisson structure \cite{HM-sklyanin}. 
The first Chern class, of sheaves in $M_\rho$, 
is equal to 
$c_1\left(\StructureSheaf{S}(n)\otimes \pi^*\StructureSheaf{\Sigma}(nD)
\right)$, where $\StructureSheaf{S}(-1)$ is the tautological line 
subbundle, of the pull back of 
$\StructureSheaf{\Sigma}(D)\oplus\StructureSheaf{\Sigma}$ 
to the ruled surface $S$. 
This first chern class is represented by 
the $n$-th multiple, of the image in $S$, of a section of the line bundle 
$\StructureSheaf{\Sigma}(D)$. The Euler characteristic of 
sheaves in $M_\rho$ is the degree $c$, of the vector bundle $E$, 
of a pair $(E,\varphi)$ in $M(GL(n),c,\rho,D)$. 
A generic pair $(E,\varphi)$ corresponds to a line bundle $F$ on 
the spectral curve of $(E,\varphi)$. 
The spectral curve belongs to the linear system of the line bundle 
$\StructureSheaf{S}(n)\otimes \pi^*\StructureSheaf{\Sigma}(nD)$. 
A proof of the Jacobi identity, applicable for the $GL(n)$ case, 
was given in \cite{DM,hurtubise-surfaces} under certain regularity 
assumptions.

Let $\sigma_\infty\subset S$ be the section at infinity,
so that $S\setminus \sigma_\infty$ is the total space of 
$\StructureSheaf{\Sigma}(D)$. The spectral curves of pairs 
$(E,\varphi)$ in $M(GL(n),c,\rho,D)$ do not meet $\sigma_\infty$. 
We note, that $M_\rho$ admits additional components consisting of sheaves,
whose one-dimensional support does meet $\sigma_\infty$. 
The support of such sheaves necessarily contain $\sigma_\infty$,
as well as some fibers of $\pi:S\rightarrow\Sigma$. 
The determinant line-bundle of such sheaves belongs to 
the component of $\Pic(S)$ containing the line-bundle 
$\StructureSheaf{S}(n)\otimes \pi^*\StructureSheaf{\Sigma}(nD)$, 
but it need not be equal to 
$\StructureSheaf{S}(n)\otimes \pi^*\StructureSheaf{\Sigma}(nD)$. 

For a general $G$, we adopt throughout section 
\ref{sec-induced-poisson-structures}
the ideal assumptions: 1) the representation 
$\rho$ induces a morphism
\begin{equation}
\label{eq-idealized-morphism}
e \ : \ M(G,c,\orbit) \ \ \ \longrightarrow \ \ \ M(GL(n),\rho(c),D)
\end{equation}
into a smooth moduli space of simple sheaves.
2) The tensor $\Psi_\rho$, defined on $M(GL(n),\rho(c),D)$ 
in \cite{HM-sklyanin} (see also (\ref{eq-Psi-M-D}) below), 
is a Poison structure. 
3) The moduli space $M(GL(n),\rho(c),\rho(\orbit))$ is a symplectic
symplectic leaf of $M(GL(n),\rho(c),D)$, where $\rho(\orbit)$ is the 
double orbit in $GL(n)((t))$ containing the image of $\orbit$ via $\rho$. 
Assumption 3 is proven in Lemma 
\ref{lemma-symplectic-leaves-are-loci-with-fixed-singularity-type}. 
Assumption 2 will be dropped in section \ref{sec-fourier-mukai}
and all three assumptions will be dropped in section 
\ref{sec-non-simple-sheaves}.

We describe in Lemma \ref{lemma-induced-poisson-structure} 
a restriction operation for Poisson structures. The essence of the
Lemma is that, once a Poisson tensor can be restricted, then there 
aren't any differential obstructions for its restriction to be a Poisson 
structure. This parallels the symplectic case, in which the pull back 
of a closed form is automatically closed. 
Let $(X,\Psi_X)$ be a smooth Poisson algebraic variety, 
$e:M\rightarrow X$ a morphism from a smooth variety $M$, and 
$\Psi_M\in H^0(M,TM\otimes TM)$ a section. 
Let $\restricted{T^*X}{M}$ be the pull back of the cotangent 
bundle by $e$. Denote by 
$\Psi_X^{-1}(TM)$ the subsheaf of $\restricted{T^*X}{M}$, which is mapped 
via $\Psi_X$ to the image of $de:TM\rightarrow \restricted{TX}{M}$. 

\begin{lemma-definition}
\label{lemma-induced-poisson-structure}
If all the conditions below are satisfied, 
then $\Psi_M$ is a Poisson structure. 
We refer to $\Psi_M$  as the 
{\em restriction of the Poisson structure $\Psi_X$}.
\begin{enumerate}
\item 
\label{lemma-item-injective-differential}
The differential $de$ is everywhere injective.
\item
\label{lemma-item-image-of-psi-M-in-image-of-psi-X}
The image of $de\circ \Psi_M: T^*M\rightarrow \restricted{TX}{M}$ 
is contained in the image of $\Psi_X$.
\item
\label{lemma-outer-is-commutative}
The outer square of the following diagram is commutative.
\[
\begin{array}{ccccc}
\Psi_X^{-1}(TM) & \rightarrow & \restricted{T^*X}{M} & 
\LongRightArrowOf{\Psi_X} & \restricted{TX}{M}
\\
e^* \ \downarrow & & e^* \ \downarrow & & \uparrow \ de
\\
T^*M & = & T^*M & \LongRightArrowOf{\Psi_M} & TM.
\end{array}
\]
\end{enumerate}
\end{lemma-definition}

\begin{rem}
{\rm
\begin{enumerate}
\item
Lemma \ref{lemma-induced-poisson-structure} can be formulated
in terms of the Poisson variety $(X,\Psi_X)$ and an embedding 
$e: M \hookrightarrow X$, without giving the tensor $\Psi_M$. 
The Poisson structure $\Psi_X$ can be restricted to $M$, if and only if
the restriction of $\Psi_X$ to $\Psi_X^{-1}(TM)$ factors through $e^*$. 
\item
In the symplectic case, where $\Psi_X$ is non-degenerate, 
commutativity of the outer square means that $\Psi_M$ is 
non-degenerate as well and the $2$-form $\Psi_X^{-1}$ restricts to 
$\Psi_M^{-1}$. 
\item
Note that the right square in the above diagram does not commute in
general. 
\end{enumerate}
}
\end{rem}

\medskip
\noindent
{\bf Proof of Lemma \ref{lemma-induced-poisson-structure}:}
The Lemma is local and even formal, so we may assume, for simplicity of
notation, that $e$ is an embedding. 
The anti-symmetry of $\Psi_M$ follows from that of $\Psi_X$. 
Let $\theta_X$ be pairing obtained from the Poisson
structure as in (\ref{eq-pairing-theta}). 
The Jacobi identity on $X$ 
is known to be equivalent to the vanishing of the section 
$\tilde{d}\theta_X$ of $\Wedge{3}TX$ defined by:
\[
\tilde{d}\theta_X(\alpha,\beta,\gamma) \ \ \ := \ \ 
\Psi_X(\alpha)\theta_X(\beta,\gamma) - 
\langle[\Psi_X(\alpha),\Psi_X(\beta)], \gamma\rangle + 
{\rm cp}(\alpha,\beta,\gamma),
\]
for $1$-forms $\alpha$, $\beta$, and $\gamma$ on $X$
(see \cite{Bo2} Proposition 1.1 or the equalities 
(\ref{eq-cyclically-shifted-Jacobi-identity}) below). 
Above, ${\rm cp}(\alpha,\beta,\gamma)$ denotes the terms obtained from the 
previous two by cyclic permutations of $\alpha$, $\beta$, 
and $\gamma$. The first term on the right hand side is Lie derivative of the 
function by the vector field $\Psi_X(\alpha)$. 
The second term involves the Lie bracket of vector fields on $X$. 

It is easy to check that $\tilde{d}\theta_M$ is a tensor 
(i.e., a section of $\Wedge{3}TM$). It is, in fact, a section of the image of
$\Wedge{3}\Psi_M : \ \Wedge{3}T^*M\rightarrow \Wedge{3}TM$. 
It suffices to check the last statement locally for exact $1$-forms
$\alpha_i=df_i$. Expanding the first term using Lie derivatives, we get
$\Psi_M(\alpha_1)\langle\Psi_M(\alpha_2),\alpha_3\rangle  = 
\langle[\Psi_M(\alpha_1),\Psi_M(\alpha_2)],\alpha_3\rangle +
\langle\Psi_M(\alpha_2),L_{\Psi_M(\alpha_1)}\alpha_3\rangle$, and thus
\begin{eqnarray}
\nonumber
\tilde{d}\theta_M(\alpha_1,\alpha_2,\alpha_3) & = & 
\langle\Psi_M(\alpha_2),L_{\Psi_M(\alpha_1)}\alpha_3\rangle+ 
{\rm cp}(1,2,3)  \ \ = \ \mbox{by exactness}
\\
\nonumber
& = &  
\langle\Psi_M(\alpha_2),d\langle\Psi_M(\alpha_1),\alpha_3\rangle\rangle+
{\rm cp}(1,2,3)
\ \ = 
\\
\label{eq-cyclically-shifted-Jacobi-identity}
& = &  
\{f_2,\{f_1,f_3\}\} + {\rm cp}(1,2,3).
\end{eqnarray}
The last term depends only on $\Psi_M(\alpha_i)$, $i=1,2,3$. 


Denote by $T(X,M)$ the subsheaf of $TX$ of vector fields tangent to $M$. 
$T(X,M)$ is a sheaf of Lie sub-algebras. Moreover, the restriction 
$T(X,M)\rightarrow TM$ is a Lie algebra homomorphism. 
Assumptions 
\ref{lemma-item-image-of-psi-M-in-image-of-psi-X} and 
\ref{lemma-outer-is-commutative} imply 
that the restriction of $\tilde{d}\theta_X$ to 
$\Psi_X^{-1}T(X,M)$ is the pull-back of $\tilde{d}\theta_M$. In other words,
if the one forms $\alpha$, $\beta$, and $\gamma$ are sections of the subsheaf
$\Psi_X^{-1}T(X,M)$ of $T^*X$, then we have the equality 
\[
\left(\tilde{d}\theta_X(\alpha,\beta,\gamma)\restricted{\right)}{M} 
\ \ \ = \ \ \ 
\tilde{d}\theta_M(\restricted{\alpha}{M},\restricted{\beta}{M},
\restricted{\gamma}{M}).
\]
Assumptions 
\ref{lemma-item-image-of-psi-M-in-image-of-psi-X} and 
\ref{lemma-outer-is-commutative} imply that the restriction
$\Psi_X^{-1}(TM) \rightarrow \Psi_M(T^*M)$ is surjective. 
Thus, the image $\Psi_M(\alpha)$, of every $1$-form $\alpha$ on $M$, 
can be lifted, locally, 
to a $1$-form on $X$ in $\Psi_X^{-1}(TM)$. 
Consequently, the vanishing of $\tilde{d}\theta_X$ implies the vanishing of
$\tilde{d}\theta_M$. 
\EndProof

\medskip
Let us check the assumptions of Lemma 
\ref{lemma-induced-poisson-structure} in case $e$ is the morphism
(\ref{eq-idealized-morphism}). 
We denote by $\Psi_G$ the tensor on 
$M(G,c,\orbit)$ provided by (\ref{eq-algebraic-2-form}). 
The tensor on $M(GL(n),\rho(c),\rho(\orbit))$
is denoted by $\Psi_\rho$. 

Given a pair $(P,\varphi)$ in $M(G,c,\orbit)$, with associated pair 
$(E,\rho(\varphi))$ in $M(GL(n),\rho(c),\rho(\orbit))$, we denote by 
$P\gg^\perp$ the subbundle of $\End(E)$ orthogonal to $P\gg$ with respect to
the trace pairing
\[
\End(E) \ \ \ = \ \ \ P\gg \ \oplus \ P\gg^\perp.
\]
The sequence (\ref{eq-short-exact-seq-defining-ad-P-phi}) fits in a 
commutative diagram
\[
\begin{array}{ccccccccc}
0 & \rightarrow & \Delta_\varphi & \rightarrow &
[P\gg\oplus P\gg] & \LongRightArrowOf{q} & ad(P,\varphi) & \rightarrow & 0
\\
& & \downarrow & & \downarrow \ \rho & & \downarrow
\\
0 & \rightarrow & \Delta_{\rho(\varphi)} & \rightarrow & 
[\End(E)\oplus \End(E)] &
\longrightarrow & ad(E,\rho(\varphi)) & \rightarrow 0.
\end{array}
\]
Since $P\gg$ and $P\gg^\perp$ are $Ad_{\rho(\varphi)}$-invariant, 
the vector bundle $\Delta_{\rho(\varphi)}$ decomposes as a direct sum
of $\Delta_\varphi$ with a subbundle $\Delta_\varphi^\perp$ of
$P\gg^\perp\oplus P\gg^\perp$. Consequently, 
$ad(E,\rho(\varphi))$ decomposes as a direct sum
\begin{eqnarray*}
ad(E,\rho(\varphi)) & = & ad(P,\varphi) \ \oplus \ ad(P,\varphi)^\perp,
 \ \ \ \mbox{where} 
\\
ad(P,\varphi)^\perp & = & [P\gg^\perp\oplus P\gg^\perp]/\Delta_\varphi^\perp.
\end{eqnarray*}
The homomorphisms
$q$, $e_1$ and $e_2$, associated to $\rho(\varphi)$ as in section
\ref{sec-2-forms}, all preserve the decompositions. 
It follows from the equalities 
$L_{\rho(\varphi)}=q\circ e_2$ and $R_{\rho(\varphi)}=q\circ e_1$,
that both homomorphisms 
$L_{\rho(\varphi)}, R_{\rho(\varphi)}:\End(E) \rightarrow ad(E,\rho(\varphi))$
preserve the decompositions. 
Consequently, the complex (\ref{eq-tangent-complex-for-symplectic-leaf}) 
of $(E,\rho(\varphi))$ decomposes as well. 
The injectivity assumption \ref{lemma-item-injective-differential} follows. 
Pull back $\Psi_\rho$ to $M(G,c,\orbit)$ as a homomorphism
\[
e^*(\Psi_\rho) \ : \  e^*T^*M(GL(n),\rho(c),\rho(\orbit)) \ \ \ 
\longrightarrow \ \ \ e^*TM(GL(n),\rho(c),\rho(\orbit)).
\]
The pulled back Poisson structure $e^*(\Psi_\rho)$ preserves the 
decomposition $e^*(\Psi_\rho)=(\Psi_G,\Psi_G^\perp)$. 
Assumptions \ref{lemma-item-image-of-psi-M-in-image-of-psi-X}
and \ref{lemma-outer-is-commutative} follow.

\subsection{A Fourier-Mukai transform}
\label{sec-fourier-mukai}
 
The moduli space $M(GL(n),\rho(c),D)$ is contained, as a Zariski open subset, 
in a moduli space $M_\rho$ of simple sheaves on $S$ 
of pure one-dimensional support. 
The moduli space $M(GL(n),\rho(c),D)$ has an anti-symmetric tensor $\Psi_\rho$ 
given in (\ref{eq-Psi-M-D}). It was verified in 
\cite{HM-sklyanin}, that $\Psi_\rho$ coincides with the Mukai-Tyurin
tensor on $M_\rho$. The Mukai-Tyurin tensor is induced
by the $\bbc^\times$-invariant Poisson structure $\psi_S$ on the surface $S$. 
One uses push-forward via $\psi_S:\omega_S\rightarrow \StructureSheaf{S}$ 
to get the homomorphism:

\begin{equation}
\label{eq-Mukai-poisson-structure}
\Ext^1(F,F\otimes\omega_S) \ \ \ 
\LongRightArrowOf{\psi_{S,*}} \ \ \ 
\Ext^1(F,F).
\end{equation}
We will now prove that $\Psi_\rho$ is a Poisson structure.
The argument in section \ref{sec-induced-poisson-structures} will imply, 
that $\Psi_G$ is the restriction of $\Psi_\rho$ via $e$,
provided we keep the ideal assumption, that we have the morphism 
(\ref{eq-idealized-morphism}) into the smooth moduli space  
$M_\rho$. 

Next, we extend Bottacin's proof of the Jacobi identity,
to the general case of all simple sheaves. 

\begin{thm}
\label{thm-Poisson-structure-on-moduli-of-sheaves-on-a-Poisson-surface}
Let $M$ be a smooth moduli space of simple sheaves on a Poisson surface
$(S,\psi_S)$. Then the Mukai-Tyurin tensor
$\Psi_M$, given in (\ref{eq-Mukai-poisson-structure}), is a Poisson structure.
\end{thm}

\begin{rem}
{\rm
Note that the smoothness assumption is necessary, as there are examples of 
singular moduli spaces of simple sheaves on a Poisson surface
(consider the structure sheaf of the curve in Example 8.6 in \cite{DM}). 
Over a symplectic surface, smoothness of any moduli space of simple 
sheaves was proven by Mukai \cite{mukai-symplectic}. 
Mukai's smoothness Theorem extends, 
over any Poisson surface $(S,\psi_S)$, to moduli spaces of simple and 
{\em locally free} sheaves. If $E$ is a locally free sheaf on 
$S$, then $\Ext^2_S(E,E)^*$ is isomorphic to $H^0(\End(E)\otimes \omega_S)$.
The latter is mapped injectively via $\psi_S$ to a subspace of 
$H^0(\End(E))$. If $E$ is simple, then $H^0(\End(E))$ is one-dimensional. 
If $S$ is not symplectic, then the image of $H^0(\End(E)\otimes \omega_S)$ 
consists of endomorphisms, which vanish along the zero divisor of $\psi_S$.
Hence, both $H^0(\End(E)\otimes \omega_S)$ and $\Ext^2_S(E,E)$ vanish. 
The vanishing of $\Ext^2_S(E,E)$ implies the smoothness of the moduli
space at $E$ \cite{artamkin,mukai-symplectic}. 
}
\end{rem}
 
The rest of section \ref{sec-fourier-mukai} is devoted to the 
proof of Theorem 
\ref{thm-Poisson-structure-on-moduli-of-sheaves-on-a-Poisson-surface}. 
The only simple sheaf with zero-dimensional support is the sky scraper sheaf 
of a point. The moduli space $M$ is then isomorphic to the surface $S$, 
$\Psi_M$ is just $\psi_S$, and the Theorem holds trivially. 

We assume from now on that the sheaves parametrized by $M$ have 
support of pure dimension $\geq 1$. We will construct a morphism 
\begin{equation}
\label{eq-morphism-Phi}
\Phi \ : \ M \ \ \ \longrightarrow \ \ \ X
\end{equation}
into a moduli space $X$ of simple and locally free sheaves on $S$. 
The tensor $\Psi_X$, given by (\ref{eq-Mukai-poisson-structure}) 
on the moduli space $X$, is known to be a Poisson structure \cite{Bo2}.
The morphism $\Phi$ will satisfy the assumptions of 
Lemma \ref{lemma-induced-poisson-structure} and
$\Psi_M$ is the restriction of $\Psi_X$. 

\medskip
Assume that all sheaves $F$, in some Zariski open subset
of $M$, satisfy the following assumptions:
\begin{equation}
\label{eq-vanishing-assumptions}
\left.
\begin{array}{rcl}
H^i(F) & = & 0, \ \ \ \mbox{for} \ i>0,
\\
H^i(F\otimes\omega_S) & = & 0, \ \ \ \mbox{for} \ i>0, \ \ \ \mbox{and}
\\
F & & \mbox{is generated by global sections}.
\end{array}
\right\}
\end{equation}
The assumptions will always hold in some neighborhood of a given 
sheaf, after we compose $e$ in 
(\ref{eq-idealized-morphism}) with tensorization by a sufficiently 
ample line bundle on $S$.

Let $S_1=S_2=S$, $\pi_i:S_1\times S_2\rightarrow S_i$
the projection on the $i$-th factor, $i=1,2$, 
$\Delta\subset S_1\times S_2$ the diagonal, and
$I_\Delta$ its ideal sheaf. 
We denote by $D(S)$ the bounded derived category of coherent sheaves on $S$.
Let 
\[
\Phi^{I_\Delta}_{S_2\rightarrow S_1} \ \ : \ \ D(S_2) \ \ \
\rightarrow \ \ \ D(S_1) 
\]
be the exact functor defined by the formula
\[
\Phi^{I_\Delta}_{S_2\rightarrow S_1}(\bullet) \ \ := \ \
R_{\pi_{1_*}}(I_\Delta\stackrel{L}{\otimes}\pi_2^*(\bullet)).
\]
We denote this functor by $\Phi$ for short.
The basic properties of such functors were studied in
\cite{mukai-fourier,orlov,bridgeland}.

The assumptions (\ref{eq-vanishing-assumptions})
imply that the object $\Phi(F)$ is represented by a sheaf, denoted
by $\Phi(F)$ as well, which is the kernel of the evaluation map
\begin{equation}
\label{eq-Phi-F}
0\rightarrow \Phi(F) \ \LongRightArrowOf{\alpha} \ 
H^0(F)\otimes\StructureSheaf{S} \
\LongRightArrowOf{ev} \ F \rightarrow 0.
\end{equation}
This can be seen by replacing $I_\Delta$ by the quasi-isomorphic complex
$[\StructureSheaf{S\times S} \rightarrow \StructureSheaf{\Delta}]$.
The transform of $F$ is seen to be the complex 
$[H^0(F)\otimes\StructureSheaf{S} \rightarrow F]$, 
which is quasi-isomorphic to $\Phi(F)$. 

\begin{lemma}
\label{lemma-properties-of-Phi-F}
Assume that $F$ is a sheaf on the surface $S$ satisfying the assumptions
 (\ref{eq-vanishing-assumptions}).
Then the sheaf $\Phi(F)$ satisfies the following properties:
\begin{enumerate}
\item
\label{lemma-item-Phi-F-is-locally-free}
If $F$ is a sheaf of pure $d$-dimensional support, $d>0$, then 
$\Phi(F)$ is locally free.
\item
\label{lemma-item-Phi-F-is-simple}
If $F$ is a simple sheaf, so is $\Phi(F)$.
\item
\label{lemma-item-Phi-is-injective-or-surjective}
If $G$ is a sheaf satisfying $H^k(G)=0$, for $k>0$, then the
homomorphism
\[
\Phi \ : \ \Ext^i_S(F,G) \ \ \ \longrightarrow \ \ \
\Ext^i_S(\Phi(F),\Phi(G))
\]
is injective for $i=0,1$ and surjective for $i=0,2$.
\end{enumerate}
\end{lemma}

\noindent
{\bf Proof:}
\ref{lemma-item-Phi-F-is-locally-free}) 
The sheaf $\Phi(F)$ is clearly torsion free. 
We have the short exact sequence defining $Q$
\[
0\rightarrow \Phi(F) \rightarrow (\Phi(F)^*)^* \rightarrow Q\rightarrow 0.
\]
$\Phi(F)$ is locally free if and only if it is reflexive, i.e., 
if and only if $Q$ vanishes. If a point $p$ is in the support of $Q$,
then there is an injective homomorphism from 
the sky-scraper sheaf $\bbc_p$ into $Q$. 
Local freeness of $\Phi(F)$ would follow, if we prove that 
the sheaf $\SheafExt^1(\bbc_p,\Phi(F))$ vanishes for
every point $p$ in $S$. Take $\SheafHom_S(\bbc_p,\bullet)$
of the short exact sequence (\ref{eq-Phi-F}) and observe that both of the
sheaves $\SheafHom_S(\bbc_p,F)$ and 
$\SheafExt^1_S(\bbc_p,H^0(F)\otimes\StructureSheaf{S})$ vanish.

\ref{lemma-item-Phi-F-is-simple}) follows from 
\ref{lemma-item-Phi-is-injective-or-surjective}. We sketch the proof of 
\ref{lemma-item-Phi-is-injective-or-surjective}. 
Let $\Gamma$ be the functor 
$\Phi^{I_\Delta^\vee\otimes \pi^*_1\omega_S[2]}_{S_1\rightarrow S_2}$.
It is known that $\Gamma$ is the left adjoint of $\Phi$. 
A straightforward calculation yields that 
$\Gamma(\Phi(F))$ fits in the distinguished triangle 
\[
A \ \longrightarrow \ \Gamma(\Phi(F)) \ \longrightarrow \
F \ \longrightarrow \ A[1]
\]
where $A=(H^1[\Phi(F)^*]^*\otimes\StructureSheaf{S})[1]$.
Given any sheaf $G$, we get 
the long exact sequence
\[
\cdots \rightarrow \Ext^i(F,G) \rightarrow \Ext^i(\Gamma\Phi(F),G)
\rightarrow \Ext^i(A,G) \rightarrow \Ext^{i+1}(F,G) \rightarrow \cdots
\]
Left adjointness implies the isomorphism $\Ext^i(\Gamma\Phi(F),G)\cong
\Ext^i(\Phi(F),\Phi(G))$. Now $\Ext^i(A,G)$ is isomorphic to
$H^1(\Phi(F)^*)\otimes H^{i-1}(G)$. Part
\ref{lemma-item-Phi-is-injective-or-surjective}
follows from the assumption that $H^k(G)=0$, for $k\neq 0$.
\EndProof

From now on, the notation $\Phi(F)$ will always denote the sheaf
associated to $F$ (rather than an object in $D(S)$). 
We let $X$ be the moduli space of simple locally free sheaves 
with the rank and Chern classes of $\Phi(F)$, $F\in M$. 
We denote by
\begin{equation}
\label{eq-d-Phi}
d\Phi \ : \ \Ext^1(F,F) \ \ \ \longrightarrow \ \ \ \Ext^1(\Phi(F),\Phi(F))
\end{equation}
the natural homomorphism. It is the differential of the morphism $\Phi$ in 
(\ref{eq-morphism-Phi}). The differential $d\Phi$ is injective, by part 
\ref{lemma-item-Phi-is-injective-or-surjective}
of Lemma \ref{lemma-properties-of-Phi-F}.
Let us check the assumptions of Lemma 
\ref{lemma-induced-poisson-structure} for the morphism $\Phi$. 
It remains to prove assumptions 
\ref{lemma-item-image-of-psi-M-in-image-of-psi-X}
and \ref{lemma-outer-is-commutative}. 
We will need an explicit description of both $d\Phi$ and its dual $d^*\Phi$.
We start with $d\Phi$. 

We will adopt the following convention. 
If $\eta\in \Ext^1(A,B)$ is an extension class of two sheaves,
then we will denote a representing extension by
\[
0\rightarrow \ B \ \LongRightArrowOf{\iota_\eta} \
V_\eta \ \LongRightArrowOf{j_\eta} \ A \rightarrow 0.
\]
Let $\eta$ be a class in $\Ext^1(F,F)$. 
Then $V_\eta$ satisfies assumptions (\ref{eq-vanishing-assumptions}) as well.
We get a commutative diagram with short exact rows and columns
\begin{equation}
\label{eq-diagram-of-d-Phi}
\begin{array}{ccccc}
\Phi(F) & \rightarrow & \Phi(V_\eta) & \rightarrow & \Phi(F)
\\
\downarrow & & \downarrow & & \downarrow
\\
H^0(F)\otimes\StructureSheaf{S} & \rightarrow & 
H^0(V_\eta)\otimes\StructureSheaf{S} & \rightarrow & 
H^0(F)\otimes\StructureSheaf{S}
\\
\downarrow & & \downarrow & & \downarrow
\\
F & \rightarrow & V_\eta & \rightarrow & F.
\end{array}
\end{equation}
The top row is an extension representing $d\Phi(\eta)$.

Next we identify explicitly the codifferential
\[
d^*\Phi \ : \ \Ext^1(\Phi(F),\Phi(F)\otimes\omega_S)
\ \ \ \longrightarrow \ \ \ \Ext^1(F,F\otimes\omega_S). 
\]
Taking $\Hom(\bullet,H^0(F)\otimes\omega_S)$ of the short exact sequence
(\ref{eq-Phi-F}), we get the long exact sequence
\begin{eqnarray*}
& \rightarrow & \Ext^i(F,H^0(F)\otimes\omega_S) \ 
\rightarrow \ \Ext^i(H^0(F)\otimes\StructureSheaf{S},H^0(F)\otimes\omega_S) \ 
\\
& \rightarrow & \Ext^i(\Phi(F),H^0(F)\otimes\omega_S) 
\ \rightarrow \ 
\Ext^{i+1}(F,H^0(F)\otimes\omega_S) \ \rightarrow \ \cdots
\end{eqnarray*}
$\Ext^i(F,H^0(F)\otimes\omega_S)$ vanishes for $i<2$, by assumptions 
(\ref{eq-vanishing-assumptions}). Consequently, the pull back homomorphism
\begin{equation}
\label{eq-pull-back-via-alpha}
\alpha^* \ : \ 
\Ext^1(H^0(F)\otimes\StructureSheaf{S},H^0(F)\otimes\omega_S) \ 
\longrightarrow \ \Ext^1(\Phi(F),H^0(F)\otimes\omega_S)
\end{equation}
is injective. We claim that $\alpha^*$ is surjective as well. 
Surjectivity follows from injectivity of 
$\Ext^2(F,H^0(F)\otimes\omega_S)\rightarrow 
\Ext^2(H^0(F)\otimes\StructureSheaf{S},H^0(F)\otimes\omega_S)$. 
Injectivity of the latter follows from Serre's Duality and the
fact that $F$ was assumed to be generated by global sections
(Assumption (\ref{eq-vanishing-assumptions})). 

Let $\epsilon$ be a class in $\Ext^1(\Phi(F),\Phi(F)\otimes\omega_S)$. 
The push forward $(\alpha\otimes 1)(\epsilon)$ is the class of an extension 
$V_{(\alpha\otimes 1)(\epsilon)}$ of $\Phi(F)$ by $H^0(F)\otimes\omega$.
Since the pull back (\ref{eq-pull-back-via-alpha}) is an isomorphism,
$(\alpha\otimes 1)(\epsilon)=\alpha^*(\delta)$
for some extension class 
$\delta\in \Ext^1(H^0(F)\otimes\StructureSheaf{S},H^0(F)\otimes\omega_s)$. 
Denote by $e:V_\epsilon \hookrightarrow V_{(\alpha\otimes 1)(\epsilon)}$
and $f:V_{(\alpha\otimes 1)(\epsilon)} \hookrightarrow V_\delta$ the natural 
homomorphisms. Both $e$ and $f$ are injective. We get the commutative diagram
with injective vertical homomorphisms:
\begin{equation}
\label{eq-codifferential-diagram}
\begin{array}{ccccccccc}
0 & \rightarrow & 
\Phi(F)\otimes\omega & \LongRightArrowOf{\iota_\epsilon} & 
V_\epsilon & \LongRightArrowOf{j_\epsilon} & \Phi(F) 
& \rightarrow & 0
\\
& & \alpha\otimes 1 \ \downarrow \ \hspace{6ex} & & 
f\circ e \ \downarrow \ \hspace{5ex} & & 
\alpha \ \downarrow \ \hspace{1ex} 
\\
0 & \rightarrow & 
H^0(F)\otimes\omega_S & \LongRightArrowOf{\iota_\delta} & V_\delta & \LongRightArrowOf{j_\delta} & H^0(F)\otimes\StructureSheaf{S}
& \rightarrow & 0.
\end{array}
\end{equation}

\noindent
We obtain the quotient short exact sequence
\[
0\rightarrow F\otimes\omega_S \ \longrightarrow \ V_{\lambda(\epsilon)} \
 \longrightarrow \ F 
\rightarrow 0.
\]
representing an extension class  
$\lambda(\epsilon)\in \Ext^1(F,F\otimes\omega_S)$. 

\begin{lemma}
$d^*\Phi(\epsilon) \ \ = \ \ -\lambda(\epsilon)$. 
\end{lemma}

\noindent
{\bf Proof:}
Let $\eta$ be a class in $\Ext^1(F,F)$ and 
$\epsilon\in \Ext^1(\Phi(F),\Phi(F)\otimes\omega_S)$. We need to prove the 
equality
\begin{equation}
\label{eq-sum-of-two-traces}
{\rm tr}(\epsilon\circ d\Phi(\eta)) 
\ + \ {\rm tr}(\lambda(\epsilon)\circ \eta) 
\ \ = \ \ 0,
\end{equation}
where the left trace is the natural homomorphism
\[
{\rm tr} \ : \ \Ext^2(\Phi(F),\Phi(F)\otimes\omega_S) \ \ \longrightarrow
\ \ H^2(\omega_S) \ \cong \ \bbc,
\]
and the right trace is from $\Ext^2(F,F\otimes\omega_S)$. 
Compose the diagrams (\ref{eq-diagram-of-d-Phi}) and 
(\ref{eq-codifferential-diagram}) to 
get the commutative diagram of $2$-extensions
\[
\begin{array}{ccccccc}
\Phi(F)\otimes\omega_S & \LongRightArrowOf{\iota_\epsilon} & V_\epsilon & 
\LongRightArrowOf{\Phi(\iota_\eta)\circ j_\epsilon} & \Phi(V_\eta) & 
\LongRightArrowOf{\Phi(j_\eta)} & \Phi(F)
\\
\alpha\otimes 1 \ \downarrow & & f\circ e \ \downarrow & & \downarrow & & 
\downarrow \ \alpha
\\
H^0(F)\otimes\omega_S & \LongRightArrowOf{\iota_\delta} & V_\delta & 
\LongRightArrowOf{H^0(\iota_\eta)\circ j_\delta} & 
H^0(V_\eta)\otimes \StructureSheaf{S} & \longrightarrow &
H^0(F)\otimes \StructureSheaf{S} 
\\
ev\otimes 1 \ \downarrow & & \downarrow & & \downarrow & & \downarrow
\\
F\otimes \omega_S & \longrightarrow & V_{\lambda(\epsilon)} & \longrightarrow 
& V_\eta & \longrightarrow & F
\end{array}
\]
with short exact columns. The left hand side of
(\ref{eq-sum-of-two-traces})
is the sum of the traces of the top and bottom $2$-extensions. 
This sum is equal to the 
trace of the $2$-extension in the middle row of the diagram. 
This middle $2$-extension is the cup product of the middle $1$-extensions
in diagrams (\ref{eq-diagram-of-d-Phi}) and 
(\ref{eq-codifferential-diagram}). The middle row of
(\ref{eq-diagram-of-d-Phi}) is a trivial extension. Hence, 
the middle $2$-extension class vanishes. 
\EndProof

The following Lemma is key to the proof of 
Assumptions 
\ref{lemma-item-image-of-psi-M-in-image-of-psi-X} and 
\ref{lemma-outer-is-commutative}
of Lemma \ref{lemma-induced-poisson-structure}.

\begin{lemma}
\label{lemma-tfae}
Suppose we are given a commutative diagram
\[
\begin{array}{ccccccccc}
0 & \rightarrow &
\Phi(F)\otimes\omega_S & \rightarrow & V_\epsilon & \rightarrow & \Phi(F) 
& \rightarrow & 0
\\
& & 
\alpha\otimes 1 \ \downarrow & & \downarrow & & \downarrow \ \alpha
\\
0 & \rightarrow &
H^0(F)\otimes\omega_S & \rightarrow & V_\gamma & \rightarrow & 
H^0(F)\otimes\StructureSheaf{S} & \rightarrow & 0
\\
& &
ev\otimes 1 \ \downarrow & & \downarrow & &  \downarrow \ ev
\\
0 & \rightarrow &
F\otimes\omega_S & \rightarrow & V_\beta & \rightarrow & F 
& \rightarrow & 0
\end{array}
\]
with short exact rows and columns. Denote the corresponding row extension 
classes by $\epsilon$, $\gamma$, and $\beta$. Then the following are 
equivalent.
\begin{enumerate}
\item
\label{eq-psi-epsilon-equal-differential-of-psi-beta}
$\psi_{S_*}(\epsilon) \ \ \ = \ \ \  d\Phi(\psi_{S_*}(\beta))$,
\item
\label{eq-psi-epsilon-in-image-of-differential}
$\psi_{S_*}(\epsilon) \ \ \  \in \ \ \  {\rm Im}(d\Phi)$,
\item
\label{eq-psi-kills-gamma}
$\psi_{S_*}(\gamma) \ \ \  = \ \ \ 0$,
\item
\label{eq-alpha-tensor-psi-kills-epsilon}
$(\alpha\otimes \psi_S)_*(\epsilon) \ \ \  = \ \ \ 0$.
\end{enumerate}
\end{lemma}

\noindent
{\bf Proof:}
(\ref{eq-alpha-tensor-psi-kills-epsilon}) $\Longleftrightarrow$ 
(\ref{eq-psi-kills-gamma}): \  
We have the equality $\alpha^*(\gamma) = (\alpha\otimes 1)_*(\epsilon)$. 
Thus, (\ref{eq-alpha-tensor-psi-kills-epsilon}) is equivalent to 
$\psi_{S_*}(\alpha^*(\gamma))=0$ as well as to 
$\alpha^*(\psi_{S_*}(\gamma))=0$.
The vanishing of $\Ext^1(F,H^0(F)\otimes\StructureSheaf{S})$ 
implies, that the pull back homomorphism
\[
\alpha^* \ : \ 
\Ext^1(H^0(F)\otimes \StructureSheaf{S},H^0(F)\otimes \StructureSheaf{S}) \ \ 
\longrightarrow \ \ 
\Ext^1(\Phi(F),H^0(F)\otimes \StructureSheaf{S})
\]
is injective. 
Thus, $\alpha^*(\psi_{S_*}(\gamma))=0$ is equivalent to 
$\psi_{S_*}(\gamma)=0$. 

\medskip
\noindent
(\ref{eq-psi-kills-gamma}) $\Longrightarrow$ 
(\ref{eq-psi-epsilon-equal-differential-of-psi-beta}): \  
The vanishing of $\psi_{S_*}(\gamma)$ implies, that the vector bundle 
$V_{\psi_{S_*}(\gamma)}$ is trivial (being the trivial extension of two trivial
vector bundles). It follows, that the homomorphism 
$V_{\psi_{S_*}(\gamma)}\rightarrow V_{\psi_{S_*}(\beta)}$ 
is conjugate to the evaluation homomorphism. Thus, 
$V_{\psi_{S_*}(\epsilon)}$ is $\Phi(V_{\psi_{S_*}(\beta)})$ 
and $\psi_{S_*}(\epsilon)=d\Phi(\psi_{S_*}(\beta))$. 

\medskip
\noindent
(\ref{eq-psi-epsilon-equal-differential-of-psi-beta}) $\Longrightarrow$ 
(\ref{eq-psi-epsilon-in-image-of-differential}): \  Clear.

\medskip
\noindent
(\ref{eq-psi-epsilon-in-image-of-differential}) $\Longrightarrow$ 
(\ref{eq-alpha-tensor-psi-kills-epsilon}): \
Let $\psi_{S_*}(\epsilon)=d\Phi(\eta)$, for some 
$\eta\in \Ext^1(F,F)$. 
Then $V_{\psi_{S_*}(\epsilon)}$ is isomorphic to 
$\Phi(V_\eta)$. Thus, by the explicit construction of $\Phi$,
we have the equality 
$(\alpha\otimes 1)_*(\psi_{S_*}(\epsilon)) = \alpha^*(\delta)$,
where $\delta$ is the extension class of the middle row in 
(\ref{eq-diagram-of-d-Phi}). But $\delta=0$. 
\EndProof

\medskip
\noindent
{\bf Proof of Assumption 
\ref{lemma-item-image-of-psi-M-in-image-of-psi-X}
of Lemma \ref{lemma-induced-poisson-structure}:
}
We need to prove that the image of
\[
\Ext^1(F,F\otimes\omega_S) \ \LongRightArrowOf{\psi_{S_*}} \ 
\Ext^1(F,F) \ \LongRightArrowOf{d\Phi} 
\Ext^1(\Phi(F),\Phi(F))
\]
is contained in the image of 
\[
\psi_{S_*} \ : \ \Ext^1(\Phi(F),\Phi(F)\otimes\omega_S) \ \ 
\longrightarrow \ \ 
\Ext^1(\Phi(F),\Phi(F)). 
\]
Start with a class $\eta\in \Ext^1(F,F\otimes\omega_S)$. 
The extension group
$\Ext^1(H^0(F)\otimes\StructureSheaf{S},F\otimes\omega_S)$
vanishes. Hence,
the pullback $ev^*(\eta)$ is the trivial extension of 
$H^0(F)\otimes\StructureSheaf{S}$ by $F\otimes\omega_S$
\[
V_{ev^*(\eta)} \ \  = \ \ 
(F\otimes\omega_S) \ \oplus \ (H^0(F)\otimes\StructureSheaf{S}).
\]
Consequently, we have a commutative diagram with surjective vertical 
homomorphisms:
\[
\begin{array}{ccccc}
H^0(F)\otimes\omega_S &  \rightarrow & 
(H^0(F)\otimes\omega_S) \ \oplus \ H^0(F)\otimes\StructureSheaf{S}
& \rightarrow & H^0(F)\otimes\StructureSheaf{S}
\\
ev\otimes 1 \ \downarrow \ \hspace{5ex}& & \downarrow & & 
\hspace{2ex} \ \downarrow \ ev
\\
F\otimes\omega_S & \rightarrow & V_\eta & \rightarrow & F 
\end{array}
\]
Complete it to a diagram as in Lemma \ref{lemma-tfae} by adding the
top row of kernels. Let $\epsilon$ be the extension class of the top row. 
The implication 
\ref{eq-psi-kills-gamma} $\Longrightarrow$ 
\ref{eq-psi-epsilon-equal-differential-of-psi-beta}
in Lemma \ref{lemma-tfae} implies that 
$\psi_{S_*}(\epsilon) = d\Phi(\psi_{S_*}(\eta))$.
\EndProof

\medskip
\noindent
{\bf Proof of Assumption \ref{lemma-outer-is-commutative}
of Lemma \ref{lemma-induced-poisson-structure}:
}
We will prove the commutativity of the square for the
Poisson structures $\Psi_M$ on $M$ and $-\Psi_X$ on $X$. 
Let $\epsilon$ be a class in $\Ext^1(\Phi(F),\Phi(F)\otimes\omega_S)$ 
and $\eta$ in $\Ext^1(F,F)$. 
Assume that 
$
\psi_{S_*}(\epsilon) = d\Phi(\eta).
$
We need to prove the equality
$
\psi_{S_*} d^*\Phi(\epsilon) = -\eta.
$
Apply the implication 
\ref{eq-psi-epsilon-in-image-of-differential} $\Longrightarrow$ 
\ref{eq-psi-epsilon-equal-differential-of-psi-beta} 
of Lemma \ref{lemma-tfae} to the 
diagram (\ref{eq-codifferential-diagram})
in the explicit construction of $d^*\Phi$. We get the equality
$\psi_{S_*}(\epsilon)=d\Phi(\psi_{S_*}(-d^*\Phi(\epsilon))).$
Assumption \ref{lemma-outer-is-commutative} follows 
from the injectivity of $d\Phi$.
\EndProof

\medskip
This completes the proof of Theorem 
\ref{thm-Poisson-structure-on-moduli-of-sheaves-on-a-Poisson-surface}.

\subsection{Restriction from a first order infinitesimal neighborhood}
\label{sec-restriction-from-first-order-neiborhood}

We will drop the assumption that the $GL(V)$-pair in diagram
(\ref{eq-four-moduli-spaces}) is simple. 
Consequently, the sheaves $F$ and $\Phi(F)$ need not be simple. 
The deformations of non-simple
sheaves may be obstructed. Thus, we work only with the sheaves of 
infinitesimal deformations of $F$ and $\Phi(F)$ over 
$M(G,c,\orbit)$. Our main point is that we can apply, both Lemma
\ref{lemma-induced-poisson-structure} and the main result of
Bottacin in \cite{Bo2}, using only the first order germ along
$M(G,c,\orbit)$ of the moduli space $X$ in diagram
(\ref{eq-four-moduli-spaces}). 
This section explains the reformulation of Bottacin's Theorem 
suitable for our purpose.

Let $M$ be a smooth algebraic or complex analytic space. 
Let $\E$ be a $0$-cochain, in the complex or \'{e}tale topology of
$M$, of families of locally free sheaves over $S\times M$.
Assume that $\SheafEnd(\E)$ glues to a global vector bundle over $S\times M$. 
Denote by $\RelExt^i_{\pi_M}(\E,\E)$ the relative extension sheaf over $M$. 
Let
\begin{equation}
\label{eq-Psi-E}
\Psi_{\E} \ : \ \RelExt^1_{\pi_M}(\E,\E\otimes\omega_S) \ \ \ 
\rightarrow \ \ \ 
\RelExt^1_{\pi_M}(\E,\E)
\end{equation}
be push-forward via the Poisson structure $\psi_S$ of $S$.

\begin{thm} \cite{Bo2}
\label{thm-bottacin}
Let $\Psi_M: T^*M\rightarrow TM$ be an anti-symmetric homomorphism.
Assume that the Kodaira-Spenser homomorphism
\[
e \ : \ TM \ \ \ \hookrightarrow \ \ \ \RelExt^1_{\pi_M}(\E,\E)
\]
is injective.
Assume furthermore that $e$ and $\Psi_M$ satisfy assumptions
\ref{lemma-item-image-of-psi-M-in-image-of-psi-X},
and \ref{lemma-outer-is-commutative} of Lemma 
\ref{lemma-induced-poisson-structure}, where 
$\Psi_X:T^*\restricted{X}{M}\rightarrow T\restricted{X}{M}$ is replaced
by (\ref{eq-Psi-E}).
Then $\Psi_M$ is a Poisson structure. 
\end{thm}

\noindent
{\bf Proof:}
We explain how the argument, given  in \cite{Bo2}, 
actually proves the Theorem. 
Let $\Psi_\E^{-1}(TM)$ be the subsheaf of 
$\RelExt^1_{\pi_M}(\E,\E\otimes\omega_S)$, which is taken by $\Psi_\E$ to 
$e(TM)$. 
Consider the pairing 
\begin{eqnarray*}
\theta_\E \ : \ \RelExt^1_{\pi_M}(\E,\E\otimes\omega_S) \ \otimes \ 
\RelExt^1_{\pi_M}(\E,\E\otimes\omega_S) & \longrightarrow &
H^2(S,\omega_S)
\\
\alpha \otimes \beta \hspace{16ex} & \mapsto & {\rm tr}(\psi_{S_*}(\alpha)\circ \beta).
\end{eqnarray*}
Note that $\theta_\E$ is anti-symmetric on $\Psi_\E^{-1}(e(TM))$,
because of the anti-symmetry of $\Psi_M$ and the assumptions of Lemma
\ref{lemma-induced-poisson-structure}.
Define the $\bbc$-linear homomorphism 
\[
\tilde{d}\theta_\E(\alpha,\beta,\gamma) \ := \ 
\Psi_\E(\alpha)\theta_\E(\beta,\gamma) - 
\langle[\Psi_\E(\alpha),\Psi_\E(\beta)],\gamma\rangle 
+ {\rm cp}(\alpha,\beta,\gamma), 
\]
for sections $\alpha$, $\beta$, $\gamma$ in $\Psi_\E^{-1}(TM)$. 
The bracket $[\Psi_\E(\alpha),\Psi_\E(\beta)]$ is the Lie bracket
in $TM$. The derivation by $\Psi_\E(\alpha)$ is the Lie derivation of
the function $\theta_\E(\beta,\gamma)$ with respect to the
vector field $\Psi_\E(\alpha)$ on $M$.
The argument in \cite{Bo2} proves that 
$\tilde{d}\theta_\E$ vanishes identically. 
The stability assumption in \cite{Bo2} was used only to deduce simplicity,
which was used to deduce that 
1) the universal endomorphism sheaf $\SheafEnd(\E)$ exists and 
2) $\theta_\E$ is anti-symmetric. 
We assumed the existence of $\SheafEnd(\E)$ and 
we need the anti-symmetry only on 
$\Psi_\E^{-1}(TM)$, where it known. 
Assumptions 
\ref{lemma-item-image-of-psi-M-in-image-of-psi-X} 
and \ref{lemma-outer-is-commutative} of Lemma 
\ref{lemma-induced-poisson-structure} imply that $\tilde{d}\theta_\E$
factors through the analogous tensor
$\tilde{d}\theta_M$ on $T^*M$. Hence, the latter vanishes. 
The vanishing of $\tilde{d}\theta_M$ is equivalent to the Jacobi identity.
\EndProof

\subsection{Non-simple sheaves}
\label{sec-non-simple-sheaves}

We can now prove Theorem 
\ref{thm-algebraic-2-form} without any additional assumptions.
The simplicity of the pairs $(P,\varphi)$ parametrized by 
$M(G,c,\orbit)$ implies that, locally in the \'{e}tale or complex topology
over moduli, we have 
a universal pair $(\P,u)$ over $\Sigma\times M(G,c,\orbit)$.
Moreover, the adjoint bundle $\P\gg$ is defined globally over
the moduli. Consequently, we get the 
$\ggl(N)$-bundle 
$\SheafEnd(\Phi(\F))$, associated to $\P\gg$ via the composition of 
$\rho$ and the functor $\Phi$. Here we use a relative version of the 
constructive description of $\Phi$ in (\ref{eq-Phi-F}), under the assumptions 
(\ref{eq-vanishing-assumptions}), 
so that $\Phi$ involves first tensorization by a sufficiently ample
line bundle. 

The triple $\{M,\E,e\}$ in Theorem \ref{thm-bottacin} is taken to be
$\{M(G,c,\orbit),\Phi(\F),d\Phi\circ e\}$, where $d\Phi$ is given in 
(\ref{eq-d-Phi}) and the latter $e$ is the infinitesimal version of 
(\ref{eq-idealized-morphism}). 
The injectivity of the latter $e$ was proven in section
\ref{sec-induced-poisson-structures}.
The injectivity of $d\Phi$, as well as the other assumptions of
Theorem \ref{thm-bottacin}, are proven by the relative version of the
arguments in section \ref{sec-fourier-mukai}.
%
Relative Fourier-Mukai transform, 
for a family of sheaves, are discussed, for example, in 
\cite{mukai-fourier,orlov,bridgeland}. 
One applies the relative version, 
to a (local) universal family of sheaves $\F$ over $S\times M$, 
in order to relate the two moduli spaces. The arguments in  
section \ref{sec-fourier-mukai} are presented for a single sheaf. 
The relative analogue is straightforward. 
One uses the cohomology and base change theorem and 
the following fact, about sections of relative extension sheaves, such 
as the sheaf $\RelExt^1_{\pi_M}(\E,\E)$ over $M$: 
Given a local section $\eta$ of 
 $\RelExt^1_{\pi_M}(\E,\E)$, over an open subset $U\subset M$,
we can find an open covering $U=\cup_{i\in I} U_i$ and extensions
\[
0\rightarrow \E \rightarrow V_{\eta,i}\rightarrow \E\rightarrow 0
\]
over $S\times U_i$ representing $\eta$ (see \cite{lange}). 
\section{Complete Integrability}
\label{sec-complete-integrability}

The invariant polynomial morphism (\ref{eq-char}) maps 
a pair $(P,\varphi)$ in $M(G,c,\orbit)$ to a section 
$char(P,\varphi)$ of $(\overline{T/W})(\orbit)$ 
(see (\ref{eq-overline-T-mod-W-of-orbit})). The quotient morphism 
$q:\overline{X}(\orbit)\rightarrow (\overline{T/W})(\orbit)$ 
is $W$-invariant. Hence, the curve 
\begin{equation}
\label{eq-spectral-curve}
C:=q^{-1}(char(P,\varphi))
\end{equation}
is a $W$-Galois cover of $\Sigma$. 
We will refer to $C$ as the {\em spectral cover}. 
Sections \ref{sec-normal-bundle}, \ref{sec-spectral-and-cameral-covers}, 
and \ref{sec-fibers-of-invariant-polynomial-map} deal with 
simply connected groups. 
In section \ref{sec-normal-bundle} we calculate the dimension of the space 
of sections of $(T/W)(\orbit)$. 
Sections \ref{sec-spectral-and-cameral-covers} and 
\ref{sec-fibers-of-invariant-polynomial-map} describe the 
fibers of the invariant polynomial morphism (\ref{eq-char}). In particular,
we prove that the fibers are Lagrangian subvarieties
(Theorem \ref{thm-complete-integrability}) and describe each fiber in terms 
of its spectral curve.
In section \ref{sec-spectral-curves-and-group-isogenies}
we introduce the new techniques needed, in order to extend 
the complete integrability Theorem \ref{thm-complete-integrability} 
to the case of non-simply connected groups. 
In section 
\ref{sec-dim-space-of-invariant-polynomials} we
prove the complete integrability of $M(G,c,\orbit)$, for any 
semi-simple group $G$ (Theorem \ref{thm-general-complete-integrability}). 
We sketch the proof for a general reductive group
in remark \ref{rem-comlete-integrability-for-reductive-groups}.
In section
\ref{sec-symplectic-surfaces} we study the fibers of the 
invariant polynomial map
(\ref{eq-char}), in terms of curves on a symplectic surface.
The topological invariant (\ref{eq-topological-invariant})
plays a crucial role, starting in section
\ref{sec-spectral-curves-and-group-isogenies}.

\subsection{The variety $(T/W)(\orbit)$ for a simply connected $G$}
\label{sec-normal-bundle}

Assume that $G$ is simply connected. 
Let $\chi_i$ be the trace of the representation $\rho_i$ of $G$, 
corresponding to the fundamental weight $\lambda_i$. 
The morphism
\[
\chi:= (\chi_1,\dots, \chi_r) \ : \ T \ \ \rightarrow \ \ \bbc^r
\] 
factors through an isomorphism $\bar{\chi}:T/W\rightarrow \bbc^r$. 

A local section of $X(\orbit)\rightarrow \Sigma$
around $p\in\Sigma$ has the form $t\tilde{a}_p$,
where $\tilde{a}_p$ is the graph of a co-character $a_p$ in the $W$-orbit 
$\orbit_p$
and $t$ is a section of the trivial $T$-bundle over $\Sigma$. 
The order of the pole of $\chi_i(t\tilde{a}_p)$ at $p$ is 
\[
-\Ord_0(\chi_i(t\tilde{a}_p)) \ \ = \ \
-\min\{\lambda(a_p) \ : \ \lambda \ 
\mbox{is a weight of} \ \rho_i\}. 
\]
This non-negative integer is independent of the choice of the co-character 
$a_p$.
Choose the co-character $a_p$, so that $-a_p$ is dominant. 
Then the minimum, on the right hand side of the equation above, 
is obtained when $\lambda$ is the highest weight $\lambda_i$. 
We get the equality
\[
-\Ord_0(\chi_i(t\tilde{a}_p)) \ \ = \ \ -\lambda_i(a_p). 
\]
We let $N$ be the vector bundle 
\begin{eqnarray}
\nonumber
N & := & \oplus_{i=1}^r \StructureSheaf{\Sigma}(S_i), \ \ \ \mbox{where}
\\
\label{eq-singularity-divisor-S-i}
S_i & := &  \sum_{p\in \Sigma}-\lambda_i(a_p)\cdot p. 
\end{eqnarray}
The morphism $\chi$ admits a relative analogue
\[
\chi \ : \ X(\orbit) \ \ \rightarrow \ \ N.
\]
Being $W$ invariant, $\chi$ factors through a morphism
\[
\bar{\chi} \ : \ (T/W)(\orbit) \ \rightarrow \ N.
\]

We have the equality 
\[
\sum_{i=1}^r\deg(S_i) \ \ = \ \ 
\frac{1}{2}\sum_{p\in \Sigma}\gamma(\orbit_p),
\]
where  $\gamma(\orbit_p)$ is the contribution of the point $p$
to the dimension of $M(G,c,\orbit)$ (see Equation 
(\ref{eq-gamma-is-dot-product-with-delta})). 
In particular, 
\begin{equation}
\label{eq-degree-of-normal-bundle}
\deg(N) \ \ = \ \ \frac{1}{2}\sum_{p\in S}\gamma(\orbit_p).
\end{equation}
Assume that none of the divisors $S_i$ vanishes. Since the $S_i$ 
are effective, the dimension of $H^0(N)$ is given by Riemann-Roch 
$h^0(N)=\deg(N)$. We conclude, that 
{\em 
the dimension of $H^0(N)$ is equal to half the dimension of $M(G,c,\orbit)$. 
}

\begin{lemma}
\label{lemma-degree-of-normal-bundle} 
Assume that none of the divisors $S_i$ vanishes. 
Then the space of global sections of 
$(T/W)(\orbit)$ is either empty, or it is a Zariski open subset
of an affine space, whose dimension 
is half the dimension of $M(G,c,\orbit)$.
\end{lemma}

\noindent
{\bf Proof:}
Follows from equation (\ref{eq-degree-of-normal-bundle}) and 
Lemma \ref{lemma-chi-bar-is-an-imersion-away-from-the-ramification-locus}.
\EndProof

\begin{lemma}
\label{lemma-chi-bar-is-an-imersion-away-from-the-ramification-locus}
If $G$ is simply connected, then the morphism  $\bar{\chi}$ 
extends to an isomorphism $(\overline{T/W})(\orbit)\cong N$. 
Consequently,
the complement of $(T/W)(\orbit)$ in $N$ has codimension $\geq 2$. 
\end{lemma}

{\bf Proof:} 
The question is local over the base curve $\Sigma$. We will prove 
the analogous local toric statement.
We let the base curve $\Sigma$ be $\bbc$, $X(\overline{F})$ 
the affine toric variety 
defined in section \ref{sec-invariant-polynomials}, and
$\pi: X(\overline{F}) \rightarrow \Sigma$ the natural projection. 
The local model of $X(\orbit)$ is the Zarisky open subset 
$X(F)\subset X(\overline{F})$.
The local analogue of $N$ is the total space of
\[
\oplus_{i=1}^r \StructureSheaf{\Sigma}(-\lambda_i(a)\cdot p_0),
\]
where $-a$ is a dominant co-character in $\orbit$ 
and $p_0\in \Sigma$ is the origin of $\bbc$. 
Clearly, $N$ is isomorphic to affine $r+1$ space. 
We have a natural morphism $\chi:X(F)\rightarrow N$. The morphism $\chi$ 
extends to $X(\overline{F})$, because the complement 
$X(\overline{F})\setminus X(F)$ has 
codimension $\geq 2$ in $X(\overline{F})$. We will prove that 
\[
\bar{\chi} \ : \ X(\overline{F})/W \ \ \rightarrow \ \ N 
\]
is an isomorphism. 

The constant function $1$ defines a section $e_i$ of the $i$-th summand
of $N$, because the divisor $S_i$ is effective. 
We let $z$ be the coordinate on $\Sigma$ vanishing at $p_0$. 
The affine coordinate ring $B$ of $N$ is generated by $\{z,y_1, \dots, y_r\}$, 
where $y_i: N \rightarrow \bbc$ corresponds to the holomorphic section 
of $N^*$ given by 
\[
y_i \ := \ z^{\lambda_i(-a)}\cdot e_i^*
\]
and $e_i^*$ is the meromorphic section of $N^*$. 

Let $A$ be the affine coordinate ring of $X(\overline{F})$. 
Denote by $A_{(0)}$ 
its localization, along the (reducible) fiber over $p_0$. 
Recall, that the subset of $X(\overline{F})$ over $\Sigma\setminus\{p_0\}$
is the trivial $T$-bundle $[\Sigma\setminus\{p_0\}]\times T$. 
Denote by $x_i\in A_{(0)}$ the $W$-invariant function, corresponding
to the trace of the $i$-th fundamental representation $\rho_i$. 
Then
\begin{equation}
\label{eq-A-localized-along-special-fiber}
A_{(0)}^W \ \ \cong \ \ \bbc[z,z^{-1}]\otimes \bbc[x_1, \dots, x_r].
\end{equation}
The morphism $\bar{\chi}$ corresponds to the algebra homomorphism 
$\bar{\chi}^*:B\rightarrow A$ given by 
\[
\bar{\chi}^*y_i=z^{\lambda_i(-a)}x_i \ \ \ \mbox{and} \ \ \ 
\bar{\chi}^*z=z.
\] 
It suffices to prove that $A^W$ is generated by 
$\{z, z^{\lambda_1(-a)}x_1, \dots, z^{\lambda_r(-a)}x_r\}$. 

Let $f\in A^W$. Using equation 
(\ref{eq-A-localized-along-special-fiber}),
we can write
\[
f \ := \ 
\sum_{n\in \bbz}\sum_{I\in Ch(T)_+}f_{n,I}\cdot z^n\cdot x^I,
\]
where the coefficients $f_{n,I}\in \bbc$ are constants, 
$Ch(T)_+$ is the integral cone spanned by the basis 
$\{\lambda_1, \dots, \lambda_r\}$ of fundamental weights, 
$x^I=\prod_{k=1}^r x_k^{i_k}$, 
and $(i_1, \dots, i_r)$ are the coordinates of $I$ in the basis 
of fundamental dominant weights.
The coefficients $f_{n,I}$ are zero, for all but finitely many 
pairs $(n,I)$. We need to prove the implication
\begin{equation}
\label{eq-implication}
f_{n,I}\neq 0 \ \ \ \Longrightarrow \ \ \ \left[n+I(a) \geq 0\right].
\end{equation}
Let $d_f$ be the minimal value of $n+I(a)$, for which $f_{n,I}\neq 0$.
Set 
\[
\bbi_f \ \ := \ \ \{I \in Ch(T)_+ \ : \ f_{d_f-I(a),I}\neq 0\}. 
\]

The graph 
\begin{eqnarray*}
\tilde{a} \ : \ \Sigma & \rightarrow & X(\overline{F}),
\\
z & \mapsto & (z,a(z)),
\end{eqnarray*}
of the co-character $a$, is a regular section of 
$\pi:X(\overline{F})\rightarrow \Sigma$. 
We will prove the implication (\ref{eq-implication}) by contradiction.
Assuming that $d_f<0$, we will prove that 
there exists $t\in T$, such that $f(t\tilde{a})$ has a pole at $p_0$. 
This will contradict the assumption that $f$ is a regular function on
$X(\overline{F})$. 

Every  weight $w$, of the fundamental representation $\rho_i$, satisfies
the inequality $w(a)\geq \lambda_i(a)$, since $-a$ is dominant and 
$\lambda_i$ is the highest weight. Set 
\begin{equation}
\label{eq-r-algebraically-independent-functions}
\chi_{i,a}(t) \ := \ \sum_{w(a)=\lambda_i(a)} w(t),
\end{equation}
where the sum is over the weights $w$ of the fundamental representation 
$\rho_i$, which satisfy $w(a)=\lambda_i(a)$. 
The coefficient of $z^{d_f}$ in $f(t\tilde{a})$ is
\begin{equation}
\label{coefficient-of-z-to-the-d-f}
\sum_{I\in \bbi_f}f_I\cdot \chi_a^I(t),
\end{equation}
where $f_I:=f_{d_f-I(a),I}$ and 
$\chi_a^I(t):= \prod_{j=1}^r \chi_{j,a}(t)^{i_j}$.

Assume that $d_f<0$. Lemma \ref{lemma-algebraic-independent}
implies, that there exists $t\in T$, for which the 
coefficient (\ref{coefficient-of-z-to-the-d-f}) of $z^{d_f}$ 
does not vanish. Hence, $f(t\tilde{a})$ has a pole at $p_0$. 
This completes the proof of Lemma 
\ref{lemma-chi-bar-is-an-imersion-away-from-the-ramification-locus}.
\EndProof

\begin{lemma}
\label{lemma-algebraic-independent}
The functions
$\{\chi_{1,a}, \dots, \chi_{r,a}\}$ over $T$, given in
(\ref{eq-r-algebraically-independent-functions}),
are algebraically independent over $\bbc$. 
\end{lemma}

{\bf Proof:}
Let $\bbi\subset Ch(T)_+$ be a finite subset 
of characters, in the cone generated by the fundamental weights, 
$c:\bbi \rightarrow \bbc^\times$ a non-vanishing function, and set 
$g:=\sum_{I\in \bbi} c_I \chi_a^I$. We will prove that
$g$ is identically zero on $T$, if and only if
$\bbi$ is the empty set. 

Assume $\bbi$ is not empty. 
Choose a {\em strictly dominant} co-character $-b$, 
for which 
\[
I(b) \ \neq \ I'(b), \ \ \ \mbox{for any two distinct elements} \ \ 
I, I'\in \bbi.
\]
Let $I\in \bbi$ be the element, for which $I(b)$ is minimal.
We will show, that the composition of $g$ with the morphism 
$b:\bbc^\times \rightarrow T$ is not the constant zero function. 
Consider the Laurent expansion of $\chi_{i,a}(b(z))$. 
Since $-b$ is strictly dominant, the order of the pole of the 
summand $\lambda_i(b(z))$ of (\ref{eq-r-algebraically-independent-functions}) 
is strictly higher than the order of the pole of the other summands
of (\ref{eq-r-algebraically-independent-functions}). 
Consequently, the coefficient of $z^{\lambda_i(b)}$, in the Laurent 
expansion of $\chi_{i,a}(b(z))$, is $1$. Moreover,  
$z^{\lambda_i(b)}$ is the lowest power of $z$ with a non-zero coefficient. 
Thus, the coefficient of $z^{I(b)}$ in $g(b(z))$ is $c_I$. 
We conclude, that $g$ is not the constant zero function.
\EndProof

\subsection{Spectral and cameral covers}
\label{sec-spectral-and-cameral-covers}

There are two natural constructions of $W$-Galois covers, of the 
base curve $\Sigma$, associated to a pair $(P,\varphi)$ in
$M_{\Sigma}(G,c,\orbit)$. We call them the {\em spectral} and 
{\em cameral covers}
(Equation (\ref{eq-spectral-curve}) and
Definition \ref{def-cameral-cover}). 
The image of the invariant polynomial map
(\ref{eq-char}) is most conveniently described in terms of the 
spectral covers. 
On the other hand, the fibers of the map (\ref{eq-char}) are 
most conveniently described in terms of the cameral covers
(section \ref{sec-fibers-of-invariant-polynomial-map}). 
The main result of this section states,  
that the spectral and cameral covers are isomorphic
(Proposition \ref{prop-spectral-equal-cameral}). 
This fact will be used in the proof of the complete integrability theorems 
\ref{thm-complete-integrability} and \ref{thm-general-complete-integrability}.
Note, that we do not assume the spectral or cameral covers to be 
smooth. 

We assume, throughout this section, that the derived subgroup $G'$ of $G$
is simply connected. The reason for this assumption
is explained in Remark \ref{rem-reason-for-simply-connected-derived-sugroup}. 
The analogue of Proposition \ref{prop-spectral-equal-cameral},
for a general semisimple group, is treated in
Lemma \ref{lemma-spectral-curves-factor-through-X-tau}. 

\begin{defi}
\label{def-regular-centralizer}
{\rm
A subalgebra $\gz$ of $\gg$ is a {\em regular centralizer},
if it is the centralizer of some element of $\gg$ and 
the dimension of $\gz$ is equal to the rank of $\gg$. 
}
\end{defi}

Let $r$ be the rank of $G$. 
Denote by $G_{reg}$ the Zariski open subset of $G$ of elements,
whose centralizer in $G$ is $r$-dimensional. 
The Lie algebra $\gz(x)$ of an element $x\in G_{reg}$ need not be 
a regular centralizer in the sense of definition
\ref{def-regular-centralizer}. 
Let $G_0\subset G_{reg}$ consists of elements $x$, 
whose centralizing subalgebra $\gz(x)$ is a regular centralizer 
in the sense of definition \ref{def-regular-centralizer} 
(i.e., it is also the centralizer of some regular element $y$ of $\gg$). 

\begin{lemma}
\label{lemma-codimension-of-complement-is-at-least-2}
$G_0$ is a Zariski open subset of $G$. 
It contains every regular element $x\in G_{reg}$, 
whose semisimple part has a centralizer $C(x_s)$ of dimension $\leq r+2$.
\end{lemma}

\noindent
{\bf Proof:}
If $x\in G_{reg}$, then $\gz(x)$ is a commutative subalgebra of rank $r$.
The set of regular centralizers is a Zariski open subset of the set of
$r$-dimensional commutative subalgebras of $\gg$. 
Hence $G_0$ is open. 
If $C(x_s)$ has dimension $r+2$, then it is an extension of 
$SL(2)$ or $PGL(2)$ by a torus. The statement reduces to the 
case of $SL(2)$ or $PGL(2)$, where it is easy to check that $G_0=G_{reg}$. 
\EndProof

\smallskip
The complement of the image of $G_0$ has codimension $\geq 2$ in $T/W$ 
(by Lemma \ref{lemma-codimension-of-complement-is-at-least-2}).

\begin{defi}
\label{moduli-of-simple-and-regular-pairs}
{\rm
Let $M_\Sigma(G,c,\orbit)_{reg}$ be the open subset of 
$M_\Sigma(G,c,\orbit)$ of simple pairs $(P,\varphi)$ satisfying:
\begin{enumerate}
\item
\label{def-item-a-section-of-P-G-0}
Away from its singularities, $\varphi$ is a section of $P(G_0)$. 
\item
\label{def-item-spectral-cover-in-unramified-over-singularities}
The spectral cover of $(P,\varphi)$ is unramified over the 
singularity divisor of $\orbit$. 
\item
\label{def-item-fiber-over-singular-points-is-a-Cartan-subalgebra}
The fiber of the centralizing subsheaf $\gz(\varphi)$ is a Cartan subalgebra  
over every point in the singularity divisor of $\orbit$. 
\end{enumerate}
}
\end{defi}

\noindent
Iwahori's Theorem implies, that every point $p\in \Sigma$ admits a 
local section of $P(G)$, which satisfies all the local conditions  
in the definition of $M_\Sigma(G,c,\orbit)_{reg}$.
(See Section \ref{sec-Iwahori} for the statement of Iwahori's Theorem). 
Condition \ref{def-item-spectral-cover-in-unramified-over-singularities}
implies, that the spectral curve (\ref{eq-spectral-curve})
is contained in the Zariski open subset 
$X(\orbit)$ of $\overline{X}(\orbit)$. 

A pair $(P,\varphi)$ in $M_\Sigma(G,c,\orbit)_{reg}$ 
is regular and generically semi-simple. 
Thus,  $\gz(\varphi)$  is a subbundle of commutative subalgebras, 
the generic fiber being a Cartan. Following \cite{donagi,donagi-gaitsgory}, 
we define the {\em cameral cover} to be the
curve of pairs $(p,\gb)$, where $p$ is a point of $\Sigma$ and 
$\gb$ is a Borel subalgebra of $\restricted{P\gg}{p}$ 
containing the fiber of $\gz(\varphi)$ at $p$. 
The scheme theoretic definition requires some notation. 
Let $\overline{G/N}$ be the space of regular centralizing Lie 
algebras. Denote by $\overline{G/T}$ the space of pairs $(\gc,\gb)$ 
of a regular centralizer and a Borel containing it. Let 
$P(\overline{G/N})$ and $P(\overline{G/T})$ be the relative analogues
modeled after the group scheme $P(G)$. Then $\gz(\varphi)$ 
determines a section $\Sigma\hookrightarrow P(\overline{G/N})$. 

\begin{defi}
\label{def-cameral-cover}
The cameral cover is the pullback of the cover 
$P(\overline{G/T})\rightarrow P(\overline{G/N})$, 
via the above section.
\end{defi}

\begin{prop}
\label{prop-spectral-equal-cameral}
The spectral and cameral covers are isomorphic. 
\end{prop}

\noindent
{\bf Proof:}
Away from the singularities of $\orbit$, the isomorphism follows from
Lemma \ref{lemma-cartesian-diagram}. 
The isomorphism between the cameral and spectral covers extends 
over the singularities of $\orbit$. This follows from Conditions 
\ref{def-item-spectral-cover-in-unramified-over-singularities} and 
\ref{def-item-fiber-over-singular-points-is-a-Cartan-subalgebra} 
of Definition \ref{moduli-of-simple-and-regular-pairs}. 
\EndProof

\smallskip
Consider the following diagram:
\begin{equation}
\label{eq-cartesian-diagram}
\begin{array}{ccccccccc}
\gt & \longleftarrow & \tilde{\gg}_{reg} & \longrightarrow & 
\overline{G/T} & \LongLeftArrowOf{\alpha} & \widetilde{G}_{0} & 
\longrightarrow & T
\\
\downarrow & & \downarrow & & \downarrow & & \downarrow & & 
\hspace{1ex} \ \downarrow \ \pi
\\
\gt/W & \longleftarrow & \gg_{reg} & \longrightarrow & \overline{G/N} & 
\LongLeftArrowOf{\delta} & G_{0} & \LongRightArrowOf{\gamma} & T/W.
\end{array}
\end{equation}
Above, $\gg_{reg}$ is the Zariski open subset of $\gg$ of regular elements
and $\tilde{\gg}_{reg}$ is the incidence variety in $\gg_{reg}\times G/B$. 
This incidence variety consists of pairs $(x,\gb)$ of a regular element $x$ 
and a Borel $\gb$ containing it. 
Let $\widetilde{G}_{reg}$ be the analogous incidence variety 
for the group $G$ and $\widetilde{G}_{0}$ the subset 
of $\widetilde{G}_{reg}$ over $G_{0}$. The top right arrow
takes a pair $(x,B)$, of a regular element $x\in G$ and a Borel $B$
containing it, to the image of $x$ under the canonical homomorphism
$B\rightarrow T$. 
The morphism $\delta$ sends an element $x$ in $G_{0}$ to 
its centralizing Lie algebra $\gz(x)$.
The morphism $\alpha$ sends a pair $(x,B)$, of an element $x\in G_0$ and a 
Borel $B$ containing it, to the pair $(\gz(x),\gb)$, consisting of
the centralizing Lie algebra of $x$ and the Lie algebra of $B$. 
Lemma \ref{lemma-borels-and-regular-centralizers} implies, that 
the morphism $\alpha$ is well defined. 

\begin{lemma}
\label{lemma-borels-and-regular-centralizers}
Let $x$ be an element of $G_0$.
Then a Borel subgroup $B$ contains $x$, if and only if
the centralizing Lie algebra $\gz(x)$ is contained in $\gb$.
\end{lemma}

\noindent
{\bf Proof:}
Denote by $\B_x$ the set of Borel subgroups of $G$ containing $x$.
Let $\B_{\gz(x)}$ be the set of Borel subgroups of $G$, 
whose Lie algebra contains $\gz(x)$. 
The Proposition in section 1.6 of
\cite{humphreys-conjugacy} proves that $\B_x \subset \B_{\gz(x)}$. 
Let $x=su$ be the Jordan decomposition of $x$, with $s$ semisimple
and $u$ unipotent. 
The identity component $Z_0(s)$, of the centralizer of $s$, is reductive and 
$x$, $u$, $s$ are contained in $Z_0(s)$ (\cite{humphreys-conjugacy} 1.12). 
The centralizer $\gz(x)$ is contained in the Lie algebra $\gz(s)$ of $Z_0(s)$. 
There is an element $\xi\in \gz(x)$, with 
$\gz(\xi)=\gz(x)$ (because $x$ is in $G_0$). 
Let $\xi=\sigma+\nu$ be the Jordan decomposition of $\xi$, 
with $\sigma$ semisimple and $\nu$ nilpotent. 
The identity component $Z_0(\sigma)$, of the centralizing subgroup of 
$\sigma$, is reductive and it contains $Z_0(s)$. 
Let $B'$ be a Borel subgroup of $Z_0(\sigma)$ containing $x$. 
Then its Lie algebra $\gb'$ contains $\gz(x)$ 
(by the Proposition in section 1.6 of \cite{humphreys-conjugacy} again). 
Hence, $\gb'$ must be the unique Borel subalgebra of 
$\gz(\sigma)$ containing $\xi$. 
A Borel in $\B_{\gz(x)}$ intersects $Z_0(\sigma)$ in a borel subgroup of 
the latter (\cite{humphreys-groups} Section 22.4). 
Suppose $B$ is in $\B_{\gz(x)}$. Then $\gb\cap \gz(\sigma)$ 
is the unique Borel subalgebra of $\gz(\sigma)$ containing $\xi$. Hence,
$B\cap Z_0(\sigma)=B'$ contains $x$. The containment 
$\B_x \supset \B_{\gz(x)}$ follows.
\EndProof

\begin{lemma}
\label{lemma-cartesian-diagram}
All the squares in  diagram (\ref{eq-cartesian-diagram}) are cartesian.
\end{lemma}

\noindent
{\bf Proof:}
The left square is well known to be cartesian. 
The middle left square is cartesian by proposition 10.3 in
\cite{donagi-gaitsgory}. 
We claim that $\widetilde{G}_{0}$ is smooth and connected. The varieties 
$\overline{G/N}$ and 
$\overline{G/T}$ are smooth and connected (\cite{donagi-gaitsgory}
propositions 1.3 and 1.5). 
Denote by $\C\rightarrow \overline{G/N}$ 
the group subscheme of $G\times \overline{G/N}$, consisting of pairs 
$(x,\ga)$, such that $\ga$ is contained in the centralizer $\gz(x)$ of $x$.
Lemma 11.2 in \cite{donagi-gaitsgory} states, that $\C$ is commutative and
smooth over $\overline{G/N}$ and irreducible as a variety. 
Denote by $\C_0$ the intersection of $\C$ with $G_0\times \overline{G/N}$.
Projection from $\C_0$ onto $G_0$ is a projective morphism, 
because $\C_0$ is the preimage of $G_0$ in the closure of $\C$ in
$G\times Gr(r,\gg)$. It follows that $\C_0\rightarrow G_0$ is 
an isomorphism (by Zariski's Main Theorem). 
$\widetilde{G}_{0}$ is smooth and connected,
being a Zariski open subscheme of the pullback of $\C$ to $\overline{G/T}$. 

We prove next that the fiber product $G_{0}\times_{T/W} T$ 
is smooth. $T/W$ is smooth because $G'$ is simply connected. 
Let $(x,\tau)$ be a point in the fiber product $G_{0}\times_{T/W} T$
over $\bar{\tau}$ in $T/W$. 
The fiber product is smooth at $(x,\tau)$, if and only if 
the image of the differential of 
$(\gamma,\pi): G_0\times T \rightarrow T/W\times T/W$
surjects onto the fiber of the normal bundle of the diagonal at 
$(\bar{\tau},\bar{\tau})$. This happens if and only if 
${\rm Im}(d_x\gamma)+{\rm Im}(d_\tau\pi)$ spans the tangent space of $T/W$
at $\bar{\tau}$. 
We show that $\gamma:G_0\rightarrow T/W$ is submersive.
Let $Z_0$ be the connected component of the center of $G$ and set 
$\overline{G}:=Z_0\times G'$ its product with the derived subgroup of $G$. 
Choose a maximal torus $T'$ in $G'$. Then
$\overline{T}:=Z_0\times T'$ is a maximal torus in $\overline{G}$. 
We get the short exact sequence
\[
0 \rightarrow K \rightarrow \overline{T} \rightarrow T \rightarrow 0,
\]
where $K$ is the anti-diagonal embedding of the intersection of $Z_0$ 
with the center of $G'$. $W\times K$ acts on $Z_0\times T'$,
and the factor $K$ acts on $Z_0$ and on $T'$ fixed point freely. 
$W$ acts trivially on $Z_0$. 
Thus, the stabilizer in $W\times K$ of a point $(a,b)$ in $Z_0\times T'$ 
is equal to the stabilizer of $b$ in $W$. The stabilizer of $b$ in $W$ 
is equal to the stabilizer of $ab$. Consequently, 
the morphism $\overline{T}/W\rightarrow T/W$ is \'{e}tale. 
$\overline{G}_{reg}\rightarrow \overline{T}/W$ is submersive, 
because it has a section. Hence, the composition 
$\overline{G}_{reg}\rightarrow T/W$ is submersive
and it factors through $G_{reg}$. Consequently, 
$G_0\rightarrow T/W$ is submersive and the fiber product is smooth.

By the universal property of the fiber product, we have a
morphism from $\widetilde{G}_{0}$ to $G_{0}\times_{T/W}T$. 
This morphism is projective because $\widetilde{G}_0\rightarrow G_0$ is. 
It is clearly of degree $1$, thus an isomorphism
by Zariski's Main Theorem. 
We conclude that the right square is cartesian. 

Next we prove that the middle left square is cartesian. 
We have seen that $G_0$ is a Zariski open subset of
the group scheme $\C$. Lemma 11.2 in \cite{donagi-gaitsgory} implies, that
$\C$ is smooth over $\overline{G/N}$.
Hence, the morphism 
$\delta:G_{0}\rightarrow \overline{G/N}$ is submersive. 
It follows that the fiber product 
$\overline{G/T}\times_{\overline{G/N}}G_0$ 
is smooth.
We have seen that $\widetilde{G}_0$ is smooth and connected. 
The morphism from $\widetilde{G}_0$ to the fiber product is 
one-to-one and onto, by Lemma \ref{lemma-borels-and-regular-centralizers}. 
The morphism from $\widetilde{G}_0$ to the fiber product is 
projective, because 
the morphism $\widetilde{G}_0\rightarrow G_0$ is projective.
$\widetilde{G}_{0}$ is isomorphic to the fiber product
of $G_{0}$ and $\overline{G/T}$ by Zariski's Main Theorem. 
\EndProof

\begin{rem}
\label{rem-reason-for-simply-connected-derived-sugroup}
{\rm
If the derived subgroup $G'$ is 
not simply connected, then 
the right hand square of diagram (\ref{eq-cartesian-diagram}) 
may fail to be cartesian over some 
regular and semisimple elements of $G$.
For example, let $G=PGL(n)$, $\xi\in \bbc^\times$ 
a primitive $k$-th root of unity, where $k=n$ if $n$ is odd, and 
$k=2n$, if $n$ is even. Let $\bar{x}\in SL(n)$ be 
${\rm diag}(1,\xi,\xi^2,\dots, \xi^{n-1})$, if $n$ is odd, 
${\rm diag}(\xi,\xi^3,\dots, \xi^{2n-1})$, if $n$ is even, 
and let $x\in PGL(n)$ be its image. Then $x$ is regular and semisimple.
Translates of $\bar{x}$, by elements of the center of $SL(n)$, belong
to the $W$-orbit of $\bar{x}$. Hence, the stabilizer of 
$x$ in $W$ is non-trivial. 
The diagram may fail to be cartesian over a locus of codimension $1$
(e.g., in the $PGL(2)$ case).
}
\end{rem}

\subsection{The fibers of the invariant polynomial map}
\label{sec-fibers-of-invariant-polynomial-map}

We prove in this section the complete integrability of 
$M_\Sigma(G,c,\orbit)_{reg}$
(Definition \ref{moduli-of-simple-and-regular-pairs}), 
for a simply connected group $G$. 
An {\em abstract Higgs bundle} is a pair $(P,\gz)$, consisting of a
principal $G$-bundle and a subbundle $\gz\subset P\gg$ of regular
centralizers (see Definition \ref{def-regular-centralizer}). 
Donagi and Gaitsgory developed an 
{\em abelianization technique} for studying abstract Higgs bundles
on a variety $X$, in terms of $W$-automorphic  $T$-bundles 
on a branched $W$-cover of $X$. Related abelianization techniques 
were used earlier in the context of the Hitchin system
\cite{hitchin-integrable-system,donagi,faltings,scognamillo}.
A generic pair $(P,\varphi)$, in our 
moduli spaces $M_\Sigma(G,c,\orbit)$, determines an abstract Higgs bundle.
One simply replaces $\varphi$ by its centralizing sheaf of
Lie subalgebras $\gz(\varphi)$. We use the results of
\cite{donagi,donagi-gaitsgory} in order to study the fibers of the
characteristic polynomial map (\ref{eq-char}).

A pair $(P,\varphi)$ in $M(G,c,\orbit)_{reg}$  is determined by 
the data consisting of 1) its spectral curve $C\subset X(\orbit)$
and 2) its abstract Higgs pair $(P,\gz(\varphi))$. 
Since $\varphi$ is generically regular and semisimple, 
the above statement follows from Proposition 
\ref{prop-spectral-equal-cameral} and the fact, that 
a regular and semisimple element $x$ in $G$ is determined by 
1) its image in $G//G$ and 2) the pair $T\subset B$ consisting of 
the centralizer of $x$ and a Borel containing it.
The data $(P,\gz(\varphi))$ determines a principal 
$T$-bundle over the cameral cover 
\cite{donagi-gaitsgory,faltings,scognamillo}. 
One first notes, that the bundle $P$ has a canonical reduction $P_B$ to a 
Borel subgroup on the cameral cover. The 
projection, from the Borel to the Cartan subgroup, 
associates to $P_B$ a $T$-bundle.  
If all co-roots of $G$ are primitive, then this procedure realizes the fiber
of (\ref{eq-char}) in $M_\Sigma(G,c,\orbit)_{reg}$ as a collection of torsors 
of the group 
\begin{equation}
\label{eq-generalized-prym}
\Hom_W(Ch(T),\Pic^0(C)), 
\end{equation}
of $W$-equivariant homomorphisms from the character lattice of $T$ to the
Jacobian of $C$. 
If $G$ has non-primitive co-roots, then this procedure is only an isogeny
onto a collection of such torsors.
The precise determination of 1) the kernel of the isogeny and 
2) which torsors occur, is a special case of the main result of 
\cite{donagi-gaitsgory}. We will not need the precise description of the 
torsors. 

We are now ready to prove the complete integrability Theorem: 

\begin{thm}
\label{thm-complete-integrability} 
Assume  $G$ is simply connected and 
none of the singularity divisors $S_i$ of $\orbit$ vanishes
(see (\ref{eq-singularity-divisor-S-i})). 
Every fiber of the characteristic polynomial map (\ref{eq-char}), 
over a section of $(T/W)(\orbit)$, 
intersects $M(G,c,\orbit)_{reg}$ in a Lagrangian 
subvariety of the latter. 
\end{thm}


\noindent
{\bf Proof of Theorem \ref{thm-complete-integrability}:}
The image of $M(G,c,\orbit)_{reg}$, under the characteristic polynomial map 
(\ref{eq-char}), is at most half the dimension of
$M(G,c,\orbit)$ (Lemma \ref{lemma-degree-of-normal-bundle}). 
Hence, the dimension of the fiber is at least half the dimension of
$M(G,c,\orbit)$. 
It remains to prove that the fiber is isotropic. 

Fix a pair $(P,\varphi)$ in $M(G,c,\orbit)_{reg}$.
The generic regularity of $\varphi$ implies, that $Z(\varphi)$ is an abelian  
group subscheme of $P(G)$. We have a natural morphism
\[
H^1(\Sigma,Z(\varphi))_0 \ \ \LongRightArrowOf{\iota} \ \ 
M(G,c,\orbit)
\] 
from the identity component of the abelian group $H^1(\Sigma,Z(\varphi))$
into $M(G,c,\orbit)$. A class $a$ of $H^1(\Sigma,Z(\varphi))_0$ is 
a right $Z(\varphi)$-torsor. It is sent to 
$(a\times_{Z(\varphi)}P,\varphi_a)$, where $Z(\varphi)$ acts on 
$a$ on the right and on $P$ on the left via the $P(G)$-action. 
$\varphi$ is a meromorphic section of $Z(\varphi)$. 
The adjoint group bundle of $a\times_{Z(\varphi)}P$
has a group subscheme $Z_a$ isomorphic to $Z(\varphi)$. We define  
$\varphi_a$ as the meromorhic section corresponding to $\varphi$,
via the isomorphizm between $Z_a$ and $Z(\varphi)$. 

The image of $\iota$ is contained in the fiber of 
(\ref{eq-char}) through $(P,\varphi)$. 
In general, the image of $\iota$ may be a subscheme of the fiber of lower 
dimension (e.g., an orbit in the boundary of a compactified Picard group). 
The global regularity of $Z(\varphi)$ implies, however, that 
the image of $\iota$ is a connected component of 
the locus in $M(G,c,\orbit)_{reg}$ consisting of pairs, whose spectral 
cover is that of $(P,\varphi)$ and whose cameral cover is isomorphic to 
that of $(P,\varphi)$ (see \cite{donagi-gaitsgory}). 
But the cameral cover is isomorphic to the spectral cover by
Proposition \ref{prop-spectral-equal-cameral}.
Hence, the image of $\iota$ is the component through $(P,\varphi)$, 
of the intersection of the fiber of 
the characteristic polynomial map (\ref{eq-char}),
with $M(G,c,\orbit)_{reg}$. 
The intersection is a Zariski open subset of the fiber. 
Hence, the Zariski tangent space at $(P,\varphi)$ 
to the image of $\iota$, is equal to
the Zariski tangent space to the fiber.

The tangent space to $M(G,c,\orbit)$ fits into the short exact sequence
\[
0\rightarrow  \ \ H^1(\ker(ad_\varphi))  \ \ \longrightarrow  \ \ 
T_{(P,\varphi)}M(G,c,\orbit) \ \  \longrightarrow  \ \ 
H^0(\coker(ad_\varphi)) \rightarrow 0.
\]
This is part of the long exact sequence of hypercohomologies, 
coming from the short exact sequence of the complexes
\[
0\rightarrow \ \ 
[\ker(ad_\varphi)\rightarrow 0] \ \ \rightarrow \ \ 
[P\gg\LongRightArrowOf{ad_\varphi} ad(P,\varphi)] \ \ \rightarrow \ \ 
[0\rightarrow \coker(ad_\varphi)] \ \ \rightarrow 0.
\]
The Zariski tangent space to the image of
$\iota$, is equal to the image of $H^1(\ker(ad_\varphi))$. 

Consider the short exact sequence 
\[
0\rightarrow  \ \ H^1(\ker(ad^*_\varphi))  \ \ \longrightarrow  \ \ 
T^*_{(P,\varphi)}M(G,c,\orbit) \ \  \longrightarrow  \ \ 
H^0(\coker(ad^*_\varphi)) \rightarrow 0.
\]
It is part of the long exact sequence of hyper-cohomologies, 
coming from the cotangent complex 
(\ref{eq-cotangent-complex-of-a-symplectic-leaf}). 
The restriction homomorphism, from $P\gg^*$ to $\ker(ad_\varphi)^*$, 
factors through 
the cokernel of the differential $ad_\varphi^*$, of the cotangent complex.
Serre's Duality identifies the composition
\[
H^1([ad(P,\varphi)^* \LongRightArrowOf{ad_\varphi^*} P\gg^*])
\ \ \ \longrightarrow \ \ \ H^0(\coker(ad^*_\varphi)) \ \ \ 
\longrightarrow \ \ \ 
H^0(\ker(ad_\varphi)^*),
\]
as the codifferential $d^*\iota$. 
The proof of Theorem 
\ref{thm-algebraic-2-form} shows, that $\ker(ad^*_\varphi)$
is isomorphic to $\ker(ad_\varphi)$. Moreover, the Poisson structure
maps the image of $H^1(\ker(ad^*_\varphi))$ onto the image of 
$H^1(\ker(ad_\varphi))$. Hence, the fiber is isotropic.
%
\EndProof

\subsection{Spectral curves and group isogenies} 
\label{sec-spectral-curves-and-group-isogenies}

Let $\iota : G \rightarrow \bar{G}$ be an isogeny of semisimple groups, 
$K$ its kernel, and $c$ a class in $\pi_1(\bar{G})$. 
We compare the spectral and cameral covers of $\Sigma$, 
associated to pairs in the moduli space $M(\bar{G},c,\orbit)$. 
The spectral curves embed in the variety 
$X(\orbit)$, constructed in section
\ref{sec-invariant-polynomials}. The spectral curve 
(\ref{eq-spectral-curve}) 
is the inverse image of a section of the $W$-quotient 
$(\bar{T}/W)(\orbit)$ of $X(\orbit)$. 
The main result of this section is 
Lemma \ref{lemma-spectral-curves-factor-through-X-tau}. It shows, 
roughly, that the morphism, from 
the cameral cover onto the spectral curve in $X(\orbit)$, factors through
an unramified Galois cover $X(\orbit,\tau)\rightarrow X(\orbit)$, 
with Galois group $K$. 
The cover $X(\orbit,\tau)$ is a partial compactification of the torus bundle, 
which is determined by the topological invariant $\tau$ in 
(\ref{eq-topological-invariant}) in section \ref{sec-topological-invariant}. 
The observation in Lemma \ref{lemma-spectral-curves-factor-through-X-tau}
is significant for several reasons: 
1) Proposition \ref{prop-spectral-equal-cameral}, 
relating the spectral and cameral covers, 
extends to the non-simply-connected case, when we replace 
the spectral curve in $X(\orbit)$ by the one in $X(\orbit,\tau)$. 
2) We will calculate the dimension, of the 
image of the invariant polynomial map (\ref{eq-char}), 
by studying the normal bundle of sections of $X(\orbit,\tau)/W$ 
(Lemma \ref{lemma-dimension-of-space-of-sections-of-X-orbit-tau-mod-W}).
3) The result will be applied in section \ref{sec-symplectic-surfaces}, 
in order to study the fibers of the invariant polynomial map (\ref{eq-char}), 
in terms of curves on a surface.

For a non-simply connected group, the relationship between the cameral 
cover and the spectral cover is more subtle. To fix ideas, let 
us consider a simple example. We work locally, letting $\Sigma = \bbc$,
and consider the 
co-character $(z^{1/2},z^{-1/2})$ of $PGL(2, \bbc)$. Assume that we 
have a morphism 
$\bar{\varphi}:\Sigma\setminus\{0\}\rightarrow PGL(2, \bbc)$, which is
conjugate over $PGL(2, \bbc[[z]])$ to 
$\pm {\rm diag}(z^{1/2},z^{-1/2})$. 
The spectral two-sheeted cover of $\Sigma$ lives in a toric variety 
$X(\orbit)$, which maps to  $\Sigma$ 
with fibers $\bbc^\times$ away from the origin, and 
a disjoint union $\bbc^\times\sqcup\bbc^\times$ at the origin. 
The spectral cover is unramified over the origin.

Let us relate the spectral cover to the spectral and so cameral curve for an
$Sl(2, \bbc)$-valued section.  The construction of the 
topological invariant (\ref{eq-topological-invariant})  
associates to $\bar{\varphi}$ a double cover $f: D\rightarrow \bbc$, branched 
at the origin. Choose one of the two lifts 
$\varphi:D\setminus\{0\}\rightarrow Sl(2, \bbc)$ of $\bar{\varphi}$. 
Denote by $f^{-1}(\orbit)$ the $W$-orbit of
co-character of $Sl(2, \bbc)$ corresponding to the lift $\varphi$. 
We get a variety $X(f^{-1}(\orbit))$, whose generic fiber over $D$
is the maximal torus of $Sl(2, \bbc)$. 
$X(f^{-1}(\orbit))$ is a $\bbz_{/2}\times \bbz_{/2}$ covering of
$X(\orbit)$; one factor corresponding to the lift $\psi$  of the deck 
transformation of the covering $f: D\rightarrow \bbc$, and the 
other to the deck transformation $\mu(\varphi)=-\varphi$ of 
the two-fold covering of $Sl(2, \bbc)$ over $PGL(2, \bbc)$, 
which gets reflected in the
maximal tori and so in the toric varieties. The 
spectral and so cameral cover $C_\varphi$ of $\varphi$ is invariant 
under the action of $\tau = \psi\circ\mu$. 

The centralizers and Borels for a $PGL(2, \bbc)$ section $\bar{\varphi}$ 
and those for its $SL(2, \bbc)$-valued 
lift $\varphi$  over $D$ correspond,  
and so the cameral cover  $C_{\bar{\varphi}}$ 
for  $\bar{\varphi}$ is the quotient of $C_\varphi$ by the action of $\tau$. 
The key point is that $C_{\bar{\varphi}}$ lives naturally in the quotient 
$X(\orbit,\tau)$ of $X(f^{-1}(\orbit))$ by $\tau$. 
There is a sequence of two-fold covers
$X(f^{-1}(\orbit))\rightarrow X(\orbit, \tau)\rightarrow X(\orbit)$ 
The last two spaces map to $\Sigma$ with fibers $\bbc^\times$ away
from the origin, and $\bbc^\times\sqcup\bbc^\times$ at 
the origin; the map between the two preserves fibers and 
is given fiberwise by $t\mapsto t^2$.

The map $X(\orbit, \tau)\rightarrow X(\orbit)$ takes the cameral cover 
$C_{\bar{\varphi}}$ to the spectral cover, in a generically bijective way. 
Over the point $z=-1\in\Sigma$ with eigenvalues 
$(z^{1/2},z^{-1/2}) = (i, -i)$, however, 
the two points in the fiber of the cameral cover map to a node of the 
spectral curve (as $i^2 = (-i)^2$). 
It is thus better to work with the space $ X(\orbit, \tau)$, 
when considering, for example, the deformation theory of the 
spectral and cameral covers. 

\medskip
We proceed to treat the general case in several steps.
In step I we construct the canonical branched cover (\ref{eq-D-orbit})
of $\Sigma$. We then lift a $\bar{G}$-pair $(P,\bar{\varphi})$ to a
$G$-pair $(P,\varphi)$ on the branched cover. 
In step II we study how the lift $(P,\varphi)$ varies 
with respect to two $K$-actions. 
In step III we construct the Galois $K$-cover of $X(\orbit)$
and prove, that the cameral covers embed in the $K$-cover. 

{\bf Step I:} 
Fix a component of $M_\Sigma(\bar{G},c,\orbit)$ and let
$\tau: \pi_1(\Sigma^0)\rightarrow \pi_1(\bar{G})$ be the 
topological invariant given in (\ref{eq-topological-invariant}), 
where $\Sigma^0$ is the complement of the support of $\orbit$.
The subgroup $K$ of $G$ is isomorphic to the quotient 
$\pi_1(\bar{G})/\pi_1(G)$, 
and we denote by 
\begin{equation}
\label{eq-bar-tau}
\bar{\tau}:\pi_1(\Sigma^0) \ \ \rightarrow \ \ K
\end{equation} 
the composition of $\tau$ with the quotient homomorphism. 
Let 
\begin{equation}
\label{eq-D-orbit}
f \ : \ D_\tau \ \ \ \longrightarrow \ \ \ \Sigma
\end{equation}
be the connected branched cover of $\Sigma$, 
compactifying the unramified Galois cover 
$D_\tau^0\rightarrow \Sigma^0$, which is 
determined by the following equality of subgroups of $\pi_1(\Sigma^0)$
\[
f_*\pi_1(D_\tau^0) \ \ \ = \ \ \ \tau^{-1}(\pi_1(G)).
\]
The Galois group of $D_\tau$ is the image of $\bar{\tau}$ in $K$.

Let $(P,\bar{\varphi})$ be a pair in $M_\Sigma(\bar{G},c,\orbit)$. 
Denote by $P(G)$ the bundle associated to $P$ via the conjugation
action of $\bar{G}$ on $G$. The topological invariant 
\[
\tau_{f^*\bar{\varphi}} \ : \ \pi_1(D_\tau^0) \ \ \ \rightarrow \ \ \
\pi_1(\bar{G}),
\]
of the pulled-back section $f^*(\bar{\varphi})$, is the composition 
$\tau\circ f_*$. 
By the definition of $f$, the image of $\tau_{f^*\bar{\varphi}}$ 
is contained in $\pi_1(G)$. Hence, the image of
$(f^*\bar{\varphi})_*:\pi_1(D_\tau^0) \rightarrow \pi_1(f^*P(\bar{G}))$
is contained in the image of $\pi_1(f^*P(G))$, under the 
covering map $f^*P(G)\rightarrow f^*P(\bar{G})$.
Consequently, there exists a lift 
\begin{equation}
\label{eq-tautologival-varphi}
\varphi : D_\tau \ \ \ \longrightarrow \ \ \ f^*(P)(G),
\end{equation}
of $f^*\bar{\varphi}$ to a section $\varphi$ of $f^*P(G)$. 
The lift $\varphi$ is unique up to the action of the Galois group $K$,
acting on $f^*P(G)$ via translations by global sections.

We will assume, throughout the rest of this section, the following condition. 
\begin{condition}
\label{cond-bar-tau-is-surjective}
The homomorphism $\bar{\tau}$ in (\ref{eq-bar-tau}) is surjective.
\end{condition}

\noindent
The general case can be reduced to the case satisfying condition
\ref{cond-bar-tau-is-surjective} as follows.
If $\bar{\tau}$ is not surjective, then it determines 
an intermediate group $\bar{G}'$, the quotient of $G$ by the image of
$\bar{\tau}$. The discussion above shows, that $\bar{\varphi}$
can be lifted to a meromorphic section $\varphi'$ of $P(\bar{G}')$ 
over $\Sigma$. 
The topological invariant $\tau'$ of $\varphi'$ is a surjective homomorphism
onto the kernel of the isogeny $G\rightarrow\bar{G}'$. 

Condition \ref{cond-bar-tau-is-surjective} implies, that $D_\tau^0$ is the
fiber product 
\[
\begin{array}{ccc}
D_\tau^0 & \LongRightArrowOf{\varphi} & P(G)
\\
\downarrow & & \downarrow
\\
\Sigma& \LongRightArrowOf{\bar{\varphi}} & P(\bar{G}).
\end{array}
\]
In particular, the lift $\varphi$ is $K$-equivariant.

{\bf Step II:}
%
Let us calculate the singularity data of $\varphi$. 
The composition of a co-character $\bar{a}:\bbc^\times \rightarrow \bar{T}$,
by the branched $d$-th power covering $\bbc^\times\rightarrow \bbc^\times$, 
is $d\bar{a}$. We define on $D_\tau$ the pulled back singularity data
\[
f^{-1}(\orbit) 
\]
as follows. The orbit $f^{-1}(\orbit_p)$ consists of the co-characters 
$d_p\cdot \bar{a}_{f(p)}$, where $\bar{a}_{f(p)}\in \orbit_{f(p)}$ 
and $d_p$ is the order of the subgroup 
of $K$, stabilizing points in the fiber of $f$ over $p$. 
This subgroup is generated by the image in $K$ of any co-character $\bar{a}_p$
in $\orbit_p$.
By construction, 
$f^{-1}(\orbit_p)$ is a $W$-orbit of co-characters of $\bar{T}$,
which come from integral co-characters of the maximal torus $T$ of $G$. 
It is easy to check, that 
$f^{-1}(\orbit)$ is precisely the singularity data of $\varphi$.

Denote by $X(f^{-1}\orbit)$ the partial compactification of 
$D_\tau^0\times T$ corresponding to the singularity data 
$f^{-1}(\orbit)$. 
The local toric model of $X(f^{-1}\orbit)$, over a neighborhood of a 
ramification point $p\in D_\tau$, is related to that of 
$X(\orbit)$, over a neighborhood of $f(p)$. Recall, that 
the toric model of $X(\orbit)$, locally around $f(p)\in \Sigma$,
was determined by the triple 
$(\bbr\oplus\gt_\bbr,Ch(\bbc^\times\times\bar{T})^*,\bar{\sigma})$,
where $\bar{\sigma}$ is a convex rational polyhedral cone in 
$\bbr\oplus\gt_\bbr$, generated by the 
graphs of the co-characters in the $W$-orbit $\orbit_{f(p)}$. 
The toric variety $X(\overline{\bar{F}})$, corresponding to the fan 
$\overline{\bar{F}}$ of faces of $\bar{\sigma}$, is singular. 
The local model $X(\bar{F})$ of $X(\orbit)$ is the complement of 
orbits in $X(\overline{\bar{F}})$ of codimension $\geq 2$. 

Denote by $\bbc^\times_{\bar{e}}$ a copy of $\bbc^\times$, with a generator 
$\bar{e}$ of its co-character lattice. Set $e:=d_p\cdot \bar{e}$. 
The cone $f^{-1}\bar{\sigma}$ in $Ch(\bbc^\times_e\times T)^*\otimes\bbq$
is the same as $\bar{\sigma}$, under the natural identification 
$Ch(\bbc^\times_e\times T)^*\otimes\bbq=
Ch(\bbc^\times_{\bar{e}}\times \bar{T})^*\otimes\bbq$. 
The lattice, of the toric model of $X(f^{-1}\orbit)$, is 
the co-character lattice $Ch(\bbc^\times_e\times T)^*$.
The local model of the morphism to $D_\tau$ is the projection 
onto the lattice $Ch(\bbc^\times_e)^*$ with the cone spanned by $e$. 

The $K$-action on $D_\tau$ lifts to an action on $X(f^{-1}\orbit)$. 
It suffices to extend the action of the stabilizer $K_p$ of each
ramification point $p\in D_\tau$. The local toric model 
identifies a punctured neighborhood of $p$ with a punctured neighborhood 
of $0$ in $\bbc$ and identifies $K_p$ with $\bbz\bar{e}/\bbz e$. 
This identification embeds $K_p$ in the first
factor of the torus $\bbc^\times_e\times T$. The toric model admits
an action of the torus $\bbc^\times_e\times T$, which restricts to 
$K_p\times T$. Hence, $X(f^{-1}\orbit)$ admits a $K\times T$ action,
which restricts to a $K\times K$-action.

The quotient of $X(f^{-1}\orbit)$ by the $K\times K$ action 
is $X(\orbit)$. 

The group $K$ admits two commuting actions on the set of meromorphic sections
of $f^*(P)(G)$. The first is the push-forward action, regarding an element 
$g\in K$ as an automorphism of $D_\tau$
\[
g_*(\varphi) \ := \ (g^{-1})^*(\varphi).
\]
$K$ embeds as a group of (central) global sections of $f^*(P)(G)$. 
Multiplication by these sections provides a second $K$-action. 
The $K$-equivariance of the morphism (\ref{eq-tautologival-varphi}) 
translates to an invariance of the section $\varphi$ with
respect to the {\em diagonal} action
\begin{equation}
\label{eq-invariance-condition}
g_*(g\cdot\varphi) \ \ = \ \  \varphi.
\end{equation}
Consequently, the spectral cover of $\varphi$ in 
$X(f^{-1}\orbit)$ is invariant under the corresponding 
diagonal action of $K$ on $X(f^{-1}\orbit)$.

{\bf Step III:}
Let $X(\orbit,\tau)$ be the quotient of $X(f^{-1}\orbit)$
by the {\em diagonal} $K$-action
\begin{equation}
\label{eq-X-orbit-tau}
X(\orbit,\tau) \ \ := \ \ X(f^{-1}\orbit)/\Delta(K).
\end{equation}
Denote by 
\begin{equation}
\label{eq-X-orbit-tau-0}
X(\orbit,\tau)_0
\end{equation} 
the complement in $X(\orbit,\tau)$
of the intersection of 1) the branch locus of 
$X(\orbit,\tau)\rightarrow X(\orbit,\tau)/W$, with 2) the fibers of 
$X(\orbit,\tau)$ over the singularity divisor of $\orbit$. 
Let $X(f^{-1}\orbit)_0$ be the inverse image of 
$X(\orbit,\tau)_0$ in $X(f^{-1}\orbit)$.
We show next, that $X(\orbit,\tau)$ is the variety, 
in which the cameral covers of pairs in $M_\Sigma(\bar{G},c,\orbit)$
embed. More canonically, we have a natural $K$-orbit of embeddings, 
with respect to the $K$-action on $X(\orbit,\tau)$.

\begin{lemma}
\label{lemma-spectral-curves-factor-through-X-tau}
\begin{enumerate}
\item
\label{lemma-item-diagram-becomes-cartesion}
If $G$ is simply connected, then both 
$X(f^{-1}\orbit)_0/W$ and $X(\orbit,\tau)_0/W$ are smooth and the diagram
\begin{equation}
\label{eq-K-equivariant-cartesian-diagram-of-f-and-pi}
\begin{array}{ccc}
X(f^{-1}\orbit)_0 & \longrightarrow & X(\orbit,\tau)_0
\\ 
\downarrow  & &  \hspace{1ex} \ \downarrow \ q
\\
X(f^{-1}\orbit)_0/W & \LongRightArrowOf{\tilde{f}} & X(\orbit,\tau)_0/W
\\
\tilde{\pi} \ \downarrow \ \hspace{1ex} & & \hspace{1ex} \ \downarrow \ \pi
\\
D_\tau & \LongRightArrowOf{f} & \Sigma
\end{array}
\end{equation}
is cartesian. 
\item
\label{lemma-item-spectral-curves-factor-through-X-tau}
Let  $(P,\bar{\varphi})$ be a pair
in $M(\bar{G},c,\orbit)$ with topological invariant $\tau$,
satisfying conditions 
\ref{def-item-a-section-of-P-G-0} and 
\ref{def-item-fiber-over-singular-points-is-a-Cartan-subalgebra}
of Definition \ref{moduli-of-simple-and-regular-pairs}. 
Let $C\subset X(\orbit,\tau)$ be the $K$-quotient 
of the spectral curve 
$C_\varphi\subset X(f^{-1}\orbit)$ 
of the pair $(f^*P,\varphi)$. 
Assume, that $G$ is simply connected and $C$ is contained in 
$X(\orbit,\tau)_0$. 
Then $C$ is isomorphic to the cameral cover of $\Sigma$,
associated to the pair $(P,\bar{\varphi})$. Furthermore, 
$C$ is the scheme theoretic inverse image, of the section 
$C_\varphi/[W\times K]$ of $\pi$. 
\item
\label{lemma-item-X-tau-is-unramified-cover}
The morphism $X(\orbit,\tau)\rightarrow X(\orbit)$ 
is unramified, provided condition 
\ref{cond-weak-primitivity-of-co-characters} below holds.
\end{enumerate}
\end{lemma}

\noindent
{\bf Proof:}
\ref{lemma-item-diagram-becomes-cartesion}) 
The equality (\ref{eq-X-orbit-tau}) yields, that the diagram is 
commutative and the lower square is $K$-equivariant. 
Let $T$ be the maximal torus of $G$.
$X(\orbit,\tau)/W$ is smooth over 
$\Sigma\setminus\mbox{supp}(\orbit)$, because $T/W$ is smooth. 
By definition, $X(\orbit,\tau)_0/W$ is smooth over $\mbox{supp}(\orbit)$. 
We conclude, that $\pi: X(\orbit,\tau)_0/W\rightarrow \Sigma$
is submersive. Hence, the fiber product 
$D_\tau\times_\Sigma X(\orbit,\tau)_0/W$ is smooth. The natural morphism from
$X(f^{-1}\orbit)_0/W$ to the fiber product is surjective and of degree $1$, 
since $\tilde{f}$ is surjective and the diagram
(\ref{eq-K-equivariant-cartesian-diagram-of-f-and-pi}) 
is $K$-equivariant. It follows from Zariski's Main Theorem, that 
$X(f^{-1}\orbit)_0/W$ is isomorphic to the fiber product. 

We consider next the upper square of the diagram. 
Over every point of $X(\orbit,\tau)_0/W$, 
one of the two morphisms $\tilde{f}$ or $q$ is submersive. 
It follows, that the fiber product of $X(f^{-1}\orbit)_0/W$ and 
$X(\orbit,\tau)_0$ over $X(\orbit,\tau)_0/W$ is smooth. 
Using Zariski's Main Theorem, we conclude, 
that the upper square in the diagram is cartesian.

\ref{lemma-item-spectral-curves-factor-through-X-tau})
Proposition \ref{prop-spectral-equal-cameral} implies, 
that $C_\varphi$ is both the spectral and cameral cover
of $(f^*P,\varphi)$. 
Let $C_{\bar{\varphi}}$ be the cameral cover of $(P,\bar{\varphi})$. 
The centralizer $\gz(\varphi)$ in $f^*P(\gg)$ is the pullback of
the centralizer $\gz(\bar{\varphi})$. 
Thus, $C_\varphi$, being the cameral cover of $(f^*P,\gz(\varphi))$,
is the fiber product of $D_\tau$ with $C_{\bar{\varphi}}$ over $\Sigma$.
It follows, that $C_\varphi/K$ is isomorphic to 
$C_{\bar{\varphi}}$.  
The statement, that $C$ is the inverse image of a 
section of $\pi$, follows from the cartesian nature of
diagram (\ref{eq-K-equivariant-cartesian-diagram-of-f-and-pi}).

\ref{lemma-item-X-tau-is-unramified-cover}) 
It suffices to show, that the diagonal subgroup of $K\times K$ contains 
the subgroup of $K\times K$ stabilizing $x$, 
for every point $x\in X(f^{-1}\orbit)$.
We use the local toric model extending the $K\times T$-action
on $X(f^{-1}\orbit)$ to a $\bbc^\times\times T$-action. 
Let $x$ be a point in $X(f^{-1}\orbit)$,
over a point $p\in D_\tau$. Denote by $K_p$ the stabilizer of $p$ in $K$. 
Let $x$ belong to the component of the
fiber of $X(f^{-1}\orbit)$ over $p$, corresponding to the 
co-character $(e,d_p\cdot \bar{a}_{f(p)})$. 
The stabilizer of $x$ in $\bbc^\times\times T$ is the image of the 
co-character $(e,d_p\cdot \bar{a}_{f(p)})$. 
The stabilizer of $x$, under the $K\times K$-action, is the intersection 
\begin{equation}
\label{eq-intersection-of-images-of-graph-of-cocharacter-with-Kp-times-K}
K_p\times K \ \ \ \cap \ \ \ (e,d_p\cdot \bar{a}_{f(p)})(\bbc^\times),
\end{equation}
of $K_p\times K$ with the image of $(e,d_p\cdot \bar{a}_{f(p)})$
in $\bbc^\times\times T$. 
The stabilizer $K_p$ is the subgroup of $K=Ch(\bar{T})^*/Ch(T)^*$ 
generated by $\bar{a}_{f(p)}$. The intersection 
$K\cap d_p\bar{a}_{f(p)}(\bbc^\times)$ is the image in $K$ 
of the rank $1$ subgroup $Ch(\bar{T})^*\cap {\rm span}_\bbq\{\bar{a}_{f(p)}\}$.
Condition
\ref{cond-weak-primitivity-of-co-characters} implies
the equality $K_p=K\cap d_p\bar{a}_{f(p)}(\bbc^\times)$. 
Condition 
\ref{cond-weak-primitivity-of-co-characters} implies also, 
that $d_p\cdot \bar{a}_{f(p)}$ is a primitive co-character of $T$. 
Consequently, $d_p\cdot \bar{a}_{f(p)}$ maps the group, of $d_p$-roots of unity
in $\bbc^\times$, isomorphically onto $K\cap d_p\bar{a}_{f(p)}(\bbc^\times)$.
Hence, the intersection 
(\ref{eq-intersection-of-images-of-graph-of-cocharacter-with-Kp-times-K}) 
is the diagonal embedding of $K_p$
in $K\times K$. 
\EndProof


\begin{condition}
\label{cond-weak-primitivity-of-co-characters}
The non-zero $W$-orbits $\orbit_p$ in $Ch(\bar{T})^*$, of the 
singularity data $\orbit$, consist of 
primitive (indivisible) co-characters. 
\end{condition}

\subsection{Complete integrability for semisimple groups}
\label{sec-dim-space-of-invariant-polynomials}

We prove in this section the complete integrability of 
$M(\bar{G},c,\orbit)$, for a semisimple group $\bar{G}$
(Theorem \ref{thm-general-complete-integrability}). 
Let $G$ be a simply-connected semisimple group, $\iota:G\rightarrow \bar{G}$
an isogeny, and $K$ its kernel. Fix a singularity data $\orbit$ for $\bar{G}$
and a topological type $\tau$ of pairs in $M(\bar{G},c,\orbit)$. 
We adopt condition \ref{cond-bar-tau-is-surjective}
(which does not lead to a loss of generality, due to the remark following
condition \ref{cond-bar-tau-is-surjective}). 

In Lemma \ref{lemma-dimension-of-space-of-sections-of-X-orbit-tau-mod-W} 
we calculate the dimension of the space of 
sections of $X(\orbit,\tau)/W$, where $X(\orbit,\tau)$ is 
given in (\ref{eq-X-orbit-tau}). As expected, the dimension of the space of 
sections is half the dimension of $M(\bar{G},c,\orbit)$. This 
generalizes Lemma \ref{lemma-degree-of-normal-bundle} 
for an arbitrary semisimple group $\bar{G}$. 

Let $f:D_\tau\rightarrow \Sigma$ be the $K$-Galois branched cover given in 
(\ref{eq-D-orbit}). Given a point $q\in \Sigma$, we let $a_q$
denote a co-character in the 
$W$-orbit $\orbit_q$, such that $-a_q$ is dominant.  
Given a point $p\in D_\tau$, we let $d_p$ be the multiplicity
of $p$ in the fiber of $f$. 
Set
\begin{eqnarray*}
\nonumber
N_{\orbit,\tau} & := & 
\oplus_{i=1}^r \StructureSheaf{\D_\tau}(S_i), \ \ \ \mbox{where}
\\
S_i & := &  \sum_{p\in D_\tau}[d_p\cdot \lambda_i(-a_{f(p)})]\cdot p. 
\end{eqnarray*}
The divisors $S_i$ are all $K$-invariant.
$S_i$ need {\em not}, however, be the pullback of a divisor on $\Sigma$,
because $\lambda_i(a_{f(p)})$ need not be integral.
Nevertheless, we will show below in 
(\ref{eq-normal-bundle-is-a-pullback-of-N-b}),
that the vector bundle $N_{\orbit,\tau}$ is the pullback of a vector bundle 
on $\Sigma$. 

\begin{claim}
\label{claim-normal-bundle-is-a-pull-back}
\begin{enumerate}
\item
\label{claim-item-degree-of-normal-bundle-to-b}
Let $b : \Sigma \hookrightarrow X(\orbit,\tau)_0/W$ 
be a section, where $X(\orbit,\tau)_0$ is given in (\ref{eq-X-orbit-tau-0}). 
Denote by $N_{b(\Sigma)}$ normal bundle to the section. Then 
\begin{equation}
\label{eq-degree-N-b-is-quotient-of-degree-N-orbit-tau}
\deg\left(N_{b(\Sigma)}\right) \ \ = \ \ 
\frac{1}{\deg(f)}\deg\left(N_{\orbit,\tau}\right) \ \ = \ \ 
\frac{1}{2}\dim M(\bar{G},c,\orbit). 
\end{equation}
\item
\label{claim-item-a-sufficient-vanishing-condition-for-H-1-N-b}
$H^1(\Sigma,N_{b(\Sigma)})$ vanishes, provided the 
following inequality holds:
\begin{equation}
\label{eq-H-1-of-N-is-less-than-genus-D-tau}
\dim H^1(D_\tau,N_{\orbit,\tau}) \ \ < \ \ g_{D_\tau},
\end{equation}
where $g_{D_\tau}$ is the genus of $D_\tau$.
\end{enumerate}
\end{claim}

Let us verify inequality (\ref{eq-H-1-of-N-is-less-than-genus-D-tau}) 
is some cases. Let $\mbox{supp}(\orbit)$ be the set theoretic support of 
the singularity divisor of $\orbit$. 
Assume that 1) $K$ has order $2$ and 
2) $\lambda_i(a_q)$ does not vanish, for any fundamental weight $\lambda_i$
and for any point $q\in \mbox{supp}(\orbit)$. 
Assumption 1) implies, that the genus of $D_\tau$ satisfies the inequality
\begin{equation}
\label{eq-upper-bound-for-genus-of-D-tau}
2g_{D_\tau}-2 \ \ \ \leq \ \ \ \deg(\mbox{supp}(\orbit)).
\end{equation}
Moreover, equality holds, if and only if every singularity point $q$ of 
$\orbit$ corresponds to a $W$-orbit of co-characters of $\bar{G}$, 
which do not come from co-characters of $G$. 
We assume, in addition: 
3) If equality holds in (\ref{eq-upper-bound-for-genus-of-D-tau}),
then for at least $r-\frac{\deg({\rm supp}(\orbit))}{2}$ 
fundamental weights $\lambda_i$, 
there exists a point $q\in \mbox{supp}(\orbit)$
($q$ may depend on $\lambda_i$), such that $-\lambda_i(a_q)>\frac{1}{2}$. 
Assumption 2) implies the inequality
\[
\deg(S_i) \ \ \geq \ \ \deg(\mbox{supp}(\orbit)).
\]
In particular, $\deg(S_i)\geq 2g_{D_\tau}-2$, for $1\leq i \leq r$. 
Assumption 3) implies, that the latter inequality is strict, for at least 
$r+1-g_{D_\tau}$ fundamental weights. 
Thus, $\dim H^1(D_\tau,N_{\orbit,\tau})< g_{D_\tau}$.

\medskip
\noindent
{\bf Proof of Claim \ref{claim-normal-bundle-is-a-pull-back}:}
\ref{claim-item-degree-of-normal-bundle-to-b})
Denote by 
$\tilde{b}:D_\tau \hookrightarrow X(f^{-1}\orbit)/W$
the section lifting $b$. 
The morphism $\tilde{f}$, 
given in (\ref{eq-K-equivariant-cartesian-diagram-of-f-and-pi}), 
induces an isomorphism
\begin{equation}
\label{eq-normal-bundle-is-a-pullback-of-N-b}
N_{\tilde{b}(D_\tau)} \ \ \ \cong \ \ \ \tilde{f}^*N_{b(\Sigma)} 
\end{equation}
(Use part \ref{lemma-item-diagram-becomes-cartesion} 
of lemma \ref{lemma-spectral-curves-factor-through-X-tau}).
Lemma \ref{lemma-chi-bar-is-an-imersion-away-from-the-ramification-locus}
identifies $(T/W)(f^{-1}\orbit)$ as a Zariski open subset of the total space 
of $N_{\orbit,\tau}$. 
The first equality in 
(\ref{eq-degree-N-b-is-quotient-of-degree-N-orbit-tau}) follows. 
The second equality in 
(\ref{eq-degree-N-b-is-quotient-of-degree-N-orbit-tau}) 
follows from the definition of $N_{\orbit,\tau}$ 
and equation (\ref{eq-gamma-is-dot-product-with-delta}). 

\ref{claim-item-a-sufficient-vanishing-condition-for-H-1-N-b})
Suppose $H^1(\Sigma,N_{b(\Sigma)})$ does not vanish. Then
$N^*_{b(\Sigma)}$ has a line sub-bundle $L$ of non-negative degree. 
Using Serre's Duality, the isomorphism 
(\ref{eq-normal-bundle-is-a-pullback-of-N-b})
and the projection formula, we get
\begin{eqnarray*}
h^1(D_\tau,N_{\tilde{b}(D_\tau)}) & = & 
h^0(D_\tau,N^*_{\tilde{b}(D_\tau)}\otimes\omega_{D_\tau}) \ = \ 
h^0(\Sigma,N^*_{b(\Sigma)}\otimes f_*\omega_{D_\tau}) 
\\
& \geq & 
h^0(\Sigma,L\otimes f_*\omega_{D_\tau})  \ \geq \ g_{D_\tau}.
\end{eqnarray*}
This contradicts the 
inequality (\ref{eq-H-1-of-N-is-less-than-genus-D-tau}). 
\EndProof

\begin{lemma}
\label{lemma-dimension-of-space-of-sections-of-X-orbit-tau-mod-W}
Assume that the  inequality 
(\ref{eq-H-1-of-N-is-less-than-genus-D-tau}) holds. 
Then the space of sections of $\pi:X(\orbit,\tau)_0/W\rightarrow \Sigma$
is either empty, or a smooth quasi-projective variety,
whose dimension is half the dimension of $M(\bar{G},c,\orbit)$.
\end{lemma}

\noindent
{\bf Proof:}
Part \ref{lemma-item-diagram-becomes-cartesion}
of Lemma \ref{lemma-spectral-curves-factor-through-X-tau} implies, 
that the space of global sections of 
$\pi:X(\orbit,\tau)_0/W\rightarrow \Sigma$ is a Zariski open subset of the 
$K$-invariant locus of the space of global sections of $\tilde{\pi}$.
The space of global sections of
$\tilde{\pi}$, in turn, is a Zariski open subset of
$H^0(D_\tau,N_{\orbit,\tau})$
(Lemma \ref{lemma-chi-bar-is-an-imersion-away-from-the-ramification-locus}). 
In particular, if non-empty, 
the space of sections of $\pi:X(\orbit,\tau)_0/W\rightarrow \Sigma$
is a smooth quasi-projective variety. 
The formula, for the dimension of 
the space of sections of $\pi$,  follows from 
Claim \ref{claim-normal-bundle-is-a-pull-back}.
\EndProof

We can now generalize the complete integrability theorem
\ref{thm-complete-integrability}. 
A pair $(P,\bar{\varphi})$ in $M_\Sigma(\bar{G},c,\orbit)$,
of topological type $\tau$, determines a $W$-invariant curve $C$ in
$X(\orbit,\tau)$, as in Lemma
\ref{lemma-spectral-curves-factor-through-X-tau}. 

\begin{defi}
Let $M_\Sigma(\bar{G},c,\orbit,\tau)_{reg}$ be the open subset of
$M_\Sigma(\bar{G},c,\orbit)$, consisting of 
simple pairs $(P,\bar{\varphi})$, with topological type $\tau$, 
satisfying conditions 
\ref{def-item-a-section-of-P-G-0}, 
\ref{def-item-fiber-over-singular-points-is-a-Cartan-subalgebra}, 
and the following analogue of condition 
\ref{def-item-spectral-cover-in-unramified-over-singularities} 
of Definition \ref{moduli-of-simple-and-regular-pairs}. 
We replace condition 
\ref{def-item-spectral-cover-in-unramified-over-singularities} of 
Definition \ref{moduli-of-simple-and-regular-pairs} by the condition, that 
the curve $C$ is contained in $X(\orbit,\tau)_0$. 
\end{defi}

\begin{thm}
\label{thm-general-complete-integrability}
Assume that the inequality (\ref{eq-H-1-of-N-is-less-than-genus-D-tau}) holds. 
Every fiber of the invariant polynomial map 
(\ref{eq-char}), over a section of $(T/W)(\orbit)$,
intersects $M_\Sigma(\bar{G},c,\orbit,\tau)_{reg}$ in a 
Lagrangian subvariety of the latter.
\end{thm}

\noindent
{\bf Proof:}
The section $\bar{b}$ of $(T/W)(\orbit)$, coming from
a pair $(P,\bar{\varphi})$ of topological type $\tau$, lifts to a $K$ orbit 
of sections of $X(\orbit,\tau)/W$. 
Fix a lift $b:\Sigma\hookrightarrow X(\orbit,\tau)_0/W$ of $\bar{b}$. 
Part \ref{lemma-item-spectral-curves-factor-through-X-tau} of 
lemma \ref{lemma-spectral-curves-factor-through-X-tau} implies, that the
cameral cover of $(P,\bar{\varphi})$, is isomorphic to 
the inverse image $C$ of $b$ in $X(\orbit,\tau)_0$. 
Hence, the cameral cover, of pairs $(P,\bar{\varphi})$ in the fiber over 
the section $\bar{b}$, depends only on $\bar{b}$. 
The rest of the proof is identical to that of Theorem
\ref{thm-complete-integrability}, except that 
lemma \ref{lemma-degree-of-normal-bundle} is replaced by lemma 
\ref{lemma-dimension-of-space-of-sections-of-X-orbit-tau-mod-W}. 
\EndProof

\begin{rem}
\label{rem-comlete-integrability-for-reductive-groups}
{\rm
One can easily extend the proof of the complete integrability 
theorem \ref{thm-complete-integrability}, for a simply connected $G$, to the 
case where $G$ is a product of a torus with a simply connected group. 
Similarly, the arguments in the proof of theorem 
\ref{thm-general-complete-integrability}, 
reducing the proof from the semisimple case to the simply connected case, 
can be extended to a reduction from the general reductive group case
to the case of a product of a torus with a simply connected group. 
}
\end{rem}

\subsection{Symplectic surfaces and group isogenies} 
\label{sec-symplectic-surfaces}

The fibers of the invariant polynomial map (\ref{eq-char}) 
can be described in terms of spectral curves associated to 
a choice of a representation $\rho$. 
The resulting spectral curves are branched covers $C_\rho\rightarrow \Sigma$
of degree equal to the dimension of the representation. 
The generic such curve is birational to the quotient, of the $W$-Galois
cameral cover, by the subgroup of $W$ stabilizing the dominant weight of the
representation. 
The spectral curves naturally embed in a quasiprojective symplectic surface 
(\ref{eq-S-lambda-bar-orbit})
(a birational model of $\Sigma\times \bbc^\times$). 
When the group $G$ is $SL(n)$, 
$SO(n)$, or $Sp(2n)$, a natural choice is the standard representation. In the
$SO(n)$ and $Sp(2n)$ cases, the spectral curves are invariant with respect to 
the anti-symplectic involution of $\Sigma\times \bbc^\times$, 
acting on the $\bbc^\times$ factor by the inversion $t\mapsto \frac{1}{t}$. 
The fiber of (\ref{eq-char}) is then the Prym variety of this double cover. 
The analysis is analogous to Hitchin's discussion in the 
case of (Lie algebra-valued) Higgs pairs \cite{hitchin-integrable-system}.

Let $\iota : G \rightarrow \bar{G}$ be an isogeny of semisimple groups, 
$K$ its kernel, 
and $c$ a class in $\pi_1(\bar{G})$. 
A new feature arises, when we study the spectral curves of pairs in the 
moduli space $M(\bar{G},c,\orbit)$, in terms of a faithful 
representation of $G$. 
(the singularity data $\orbit$ need not lift to one for $G$
and the representation of $G$ does not factor through $\bar{G}$). 
The main result of this section, Lemma \ref{lemma-symplectic-covering}, 
clarifies the subtle role of another symplectic surface 
(\ref{eq-S-lambda-orbit}), a Galois covering of the symplectic surface
(\ref{eq-S-lambda-bar-orbit}) with Galois group $K$. 
The map, from each spectral curve to the 
symplectic surface (\ref{eq-S-lambda-bar-orbit}), 
factors through the symplectic covering (\ref{eq-S-lambda-orbit}). 
The symplectic covering is determined by the topological 
invariant (\ref{eq-topological-invariant}) of the 
component of the moduli $M(\bar{G},c,\orbit)$.
When $G$ is $SO(n)$ or $Sp(2n)$ and $\bar{G}$ its adjoint group, 
the surface  (\ref{eq-S-lambda-orbit}) admits an anti-symplectic involution
as well, but its fixed locus need no longer be the union of two sections. 
Rather, the fixed locus is the branched double cover 
of $\Sigma$ defined earlier in (\ref{eq-D-orbit}). 
We construct the two surfaces in two steps. 
In step I we construct the symplectic surface
associated to the choice of a representation. 
In step II we construct the symplectic covering 
(\ref{eq-S-lambda-orbit}) of the symplectic surface.
We show, that a quotient of the spectral curve of $\varphi$ 
maps into this surface. 
The general method developed in this section will be applied in
Section \ref{sec-type-D} in particular examples.

{\bf Step I:} 
Let $\lambda:T\rightarrow \bbc^\times$ be a dominant weight of a faithful
representation of $G$. The representation restricts as an injective
homomorphism from $K$ to the center of $GL(N)$. Finite subgroups of the
latter are cyclic. In particular, the intersection $K\cap \ker\lambda$ 
is trivial and $K$ must be cyclic. 
The kernel $K$ of the isogeny $\iota$ fits in the exact sequence
\[
0 \rightarrow K \rightarrow T \rightarrow \bar{T} \rightarrow 0
\]
of the corresponding maximal tori. 
Let $\{\bar{e},\bar{\alpha}\}$ be a basis of the co-character lattice of 
$\bbc^\times \times [T/(K+\ker\lambda)]$. 
The set 
\begin{equation}
\label{eq-fan-of-symplectic-surface}
\lambda(\orbit_p) \ \ := \ \ 
\{(\bar{e},\lambda(\bar{a}_p)) \ : \ \ \ \ \bar{a}_p \in \orbit_p\}.
\end{equation}
is the image via 
\[
id\times \bar{\lambda} \ : \  \bbc^\times\times \bar{T} \ \ 
\longrightarrow \ \ \bbc^\times\times [T/(K+\ker\lambda)],
\]
of the $W$-orbit of graphs of co-characters in $\orbit_p$. 
The set $\lambda(\orbit_p)$ determines a fan
$id\times \bar{\lambda}(\bar{\sigma})$, which is supported on
the $id\times \bar{\lambda}$ image of $\bar{\sigma}$ in 
$\Span_\bbr\{\bar{e},\bar{\alpha}\}$. 
The toric surface $X(id\times \bar{\lambda}(\bar{\sigma}))$ is singular.
The local model of 
\begin{equation}
\label{eq-S-lambda-bar-orbit}
S_\lambda(\orbit)
\end{equation}
is the complement in $X(id\times \bar{\lambda}(\bar{\sigma}))$
of the fixed points for the action of 
$\bbc^\times\times [T/(K+\ker\lambda)]$. 
We have a natural morphism
\begin{equation}
\label{eq-bar-lambda}
\bar{\lambda} \ : \ 
X(\orbit) \ \ \ \longrightarrow \ \ \ S_\lambda(\orbit), 
\end{equation}
extending the morphism 
$id\times\bar{\lambda}:\Sigma^0\times \bar{T}\rightarrow 
\Sigma^0\times \bbc^\times$, where
$\Sigma^0$ is the complement of the support of $\orbit$. 
The $2$-form $dz\wedge \frac{dt}{t}$ on $\Sigma\times \bbc^\times$
determines a (non-degenerate) holomorphic symplectic structure on the surface
$S_\lambda(\orbit)$ as well. The $2$-form is holomorphic, because
the morphism $\bar{\lambda}$ in (\ref{eq-bar-lambda}) is surjective 
and it pulls it back to the contraction $i(\Omega,\bar{\lambda})$, where
$\Omega$ is the $\gt$-valued $2$-form of Lemma \ref{lemma-t-valued-2-form}.
The $2$-form is non-degenerate, because its pullback 
$i(\Omega,\bar{\lambda})$ is nowhere vanishing (Lemma
\ref{lemma-t-valued-2-form}).

{\bf Step II:}
We extend in this step the morphism (\ref{eq-bar-lambda})
to a commutative $K\times K$-equivariant diagram
\[
\begin{array}{ccccc}
X(f^{-1}\orbit) & \rightarrow & X(\orbit,\tau) & \rightarrow & X(\orbit)
\\
\downarrow & & \downarrow & & \hspace{1ex} \downarrow \bar{\lambda} 
\\
S_\lambda(f^{-1}\orbit) & \rightarrow & S_\lambda(\orbit,\tau) & \rightarrow & 
S_\lambda(\orbit),
\end{array}
\]
where $X(f^{-1}\orbit)$ and $X(\orbit,\tau)$ are the varieties defined in 
section \ref{sec-spectral-curves-and-group-isogenies}. The surfaces 
$S_\lambda(f^{-1}\orbit)$ and $S_\lambda(\orbit,\tau)$ are defined below.

We have a $K\times K$ action on 
$D_\tau\times \bbc^\times$. The action on $\bbc^\times$ is the
restriction of $\lambda$ to $K$. 
Over $D_\tau^0$, we get a $K\times T$-equivariant morphism 
\[
X(f^{-1}\orbit\restricted{)}{D_\tau^0} \ \ \ = \ \ \ 
D_\tau^0\times T \ \ \ \longrightarrow \ \ \ D_\tau^0\times \bbc^\times.
\]
It extends to a regular 
$K\times T$-equivariant morphism into a surface
$S_\lambda(f^{-1}\orbit)$
\begin{equation}
\label{eq-S-lambda-f-1-orbit}
\lambda \ : \ X(f^{-1}\orbit) \ \ \ \longrightarrow \ \ \ 
S_\lambda(f^{-1}\orbit).
\end{equation}
The latter surface is a partial compactification
of $D_\tau^0\times \bbc^\times$. The local toric model of 
$S_\lambda(f^{-1}\orbit)$, over a point $p$ in $D_\tau$, is 
determined by the fan in $Ch(\bbc^\times_e\times [T/\ker\lambda])^*$,
which is the image of the cone in $Ch(\bbc^\times_e\times T)^*$ 
determined by $f^{-1}\orbit_p$.
$S_\lambda(\orbit)$ is the quotient of
$S_\lambda(f^{-1}\orbit)$ by the $K\times K$ action. 
Let $\{\bar{e},\bar{\alpha}\}$ be a basis of the co-character lattice of 
$\bbc^\times_{\bar{e}} \times [T/K+\ker\lambda]$. 
Recall, that $e:=d_p\cdot \bar{e}$, where $d_p$ is the order of the 
stabilizing subgroup  of the 
ramification point of $D_\tau$ over $p\in \Sigma$. 
The co-character lattice of $T/\ker\lambda$ is generated by 
$\alpha:=\Abs{K}\cdot \bar{\alpha}$. The toric models of 
$S_\lambda(\orbit)$ and $S_\lambda(f^{-1}\orbit)$ are
related as follows. Their fans in $\Span_\bbr\{\bar{e},\bar{\alpha}\}$
are identical; they are determined by the set of vectors
(\ref{eq-fan-of-symplectic-surface}).
The lattice of $S_\lambda(\orbit)$ is 
$\Span_\bbz\{\bar{e},\bar{\alpha}\}$.
The lattice of $S_\lambda(f^{-1}\orbit)$ is the co-character sublattice
spanned by $\{e,\alpha\}$. 

Let 
\begin{equation}
\label{eq-S-lambda-orbit}
S_\lambda(\orbit,\tau)
\end{equation}
be the quotient of $S_\lambda(f^{-1}\orbit)$ by the 
{\em diagonal} action of $K$. 
$S_\lambda(\orbit,\tau)$ is a partial compactification of the quotient of 
$D_\tau^0\times \bbc^\times$ by the 
diagonal action of $K$. The latter quotient is the total space of a  
$\bbc^\times$ local system over $\Sigma^0$. The local system is 
determined by the composition of 1) the topological invariant
(\ref{eq-topological-invariant}) with 2) the restriction of the character 
$\lambda$ to the center of the universal cover of $\bar{G}$.
The fiber of $S_\lambda(\orbit,\tau)$, over a point in $\Sigma_0$, 
is $\bbc^\times$, 
but fibers over the support of $\orbit$ have multiplicities. 
$S_\lambda(\orbit,\tau)$ admits, locally over a point in $\Sigma$, a toric model 
as well. 
Again, the fan in $\Span_\bbr\{\bar{e},\bar{\alpha}\}$ is equal to the one 
in the definition of $S_\lambda(\orbit)$ in 
(\ref{eq-S-lambda-bar-orbit}).
The sublattice is generated by $\{\bar{e}+\bar{\alpha},e,\alpha\}$. 
The local model, of the morphism to $\Sigma$, is the projection
(modulo $\alpha$) to the lattice $\Span_\bbz\{\bar{e}\}$ with 
the cone spanned by $\bar{e}$.

\begin{lemma}
\label{lemma-symplectic-covering} 
\begin{enumerate}
\item
\label{lemma-item-factorization-of-bar-lambda}
Let $(P,\bar{\varphi})$ be a pair in $M(\bar{G},c,\orbit)$. Let
$\bar{\lambda}:C\rightarrow S_\lambda(\orbit)$ be the 
restriction of the morphism (\ref{eq-bar-lambda}) to the spectral cover $C$
of $(P,\bar{\varphi})$. Over the smooth locus of $C$, 
the morphism $\bar{\lambda}$ factors through a morphism
$C/W_\lambda\rightarrow S_\lambda(\orbit,\tau)$
\[
C \ \rightarrow \ C/W_\lambda \ \rightarrow \ S_\lambda(\orbit,\tau) \ \rightarrow \
S_\lambda(\orbit).
\]
\item
\label{lemma-item-symplectic-covering}
Assume that condition \ref{cond-weak-primitivity-of-co-characters} 
holds. The quotient morphism 
$S_\lambda(\orbit,\tau)\rightarrow S_\lambda(\orbit)$, via the $K$-action, 
is an unramified covering. 
Consequently, $S_\lambda(\orbit,\tau)$ is symplectic as well. 
\item
\label{lemma-item-anti-symplectic-involution}
Assume $G$ is $SO(n)$, or $Sp(2n)$, $\lambda$ is the 
standard representation, and $\bar{G}$ is the adjoint group.
The stabilizer $W_\lambda$, of $\lambda$, has index $2$ in the subgroup 
$\widetilde{W}_\lambda$ of $W$ stabilizing the set 
$\{\lambda,-\lambda\}$. The order $2$ element $\eta$ of 
$\widetilde{W}_\lambda/W_\lambda$ acts on 
$S_\lambda(\orbit,\tau)$ as an anti-symplectic involution. The fixed locus of
$\eta$ is a curve, which is isomorphic to a dense open subset 
(\ref{eq-D-tau-prime}) of $D_\tau$. 
\end{enumerate}
\end{lemma}

\noindent
{\bf Proof:}
\ref{lemma-item-factorization-of-bar-lambda}) is clear.
\ref{lemma-item-symplectic-covering})
The proof is similar to that of part 
\ref{lemma-item-X-tau-is-unramified-cover} 
of Lemma \ref{lemma-spectral-curves-factor-through-X-tau}.

\ref{lemma-item-anti-symplectic-involution})
Let $D_\tau'$ be the Zariski open subset of $D_\tau$ defined by 
\begin{equation}
\label{eq-D-tau-prime}
D_\tau' \ := \ 
\{p\in D_\tau: \lambda(d_p\bar{a}_p)=0, \ \ \mbox{for some}
\ \ \bar{a}_p\in \orbit_{f(p)}\}.
\end{equation}
Clearly, $D_\tau'$ contains $D_\tau^0$. We claim, that $D_\tau'$ 
is disjoint from the ramification locus of $f:D_\tau\rightarrow \Sigma$. 
If $p$ is a ramification point, then $K$ intersects non-trivially 
the image in $T$ of the co-character $d_p\bar{a}_p$. 
Since $\ker\lambda$ intersects
$K$ trivially, then $d_p\bar{a}_p$ can not be in the kernel of
$\lambda$. 

Let $p$ be a point in the complement of $D_\tau'$. 
The involution $\eta$ takes $\lambda(d_p\bar{a}_p)$ to 
$-\lambda(d_p\bar{a}_p)$. Hence, $\eta$ takes every component, of the fiber of
$S_\lambda(\orbit,\tau)$ over $p$, to a different component. 

Let $x$ be a point in $S_\lambda(\orbit,\tau)$ over 
a point $p$ in $D_\tau'$. The component, of the fiber of 
$S_\lambda(\orbit,\tau)$ containing $x$, corresponds to a co-character $\bar{a}_p$,
which comes from $T$. If $\lambda(\bar{a}_p)\neq 0$, then $\eta(x)\neq x$. 
If $\lambda(\bar{a}_p)=0$, then the component of $x$ is 
naturally identified with $T/\ker\lambda$. 
Let $\kappa$ be the diagonal involution of 
$S_\lambda(f^{-1}\orbit)$ in $K\times K$. 
Then $\kappa(p,t)=(\kappa(p),-t)$. 
Represent $x$ as a $K$-orbit $\{(p,t), (\kappa(p),-t)\}$. 
The involution $\eta$ in $\widetilde{W}_\lambda/W_\lambda$
takes this $K$-orbit to $\{(p,\frac{1}{t}), (\kappa(p),-\frac{1}{t})\}$. 
The two orbits are equal precisely when $t^2=1$. 
We conclude, that the fixed locus of $\eta$ is isomorphic to the 
$\kappa$-quotient of $D_\tau'\times \{1,-1\}$, which is 
isomorphic to $D_\tau'$. 
\EndProof


\medskip
Assume now, that $G$ is $SL(n)$, $SO(n)$, or $Sp(2n)$, $\lambda$ is the 
standard representation, and $\bar{G}$ is the adjoint group. 
When $G=SL(n)$ and $\bar{G}=PGL(n)$, then the fiber of (\ref{eq-char})
is isogenous to the prym variety of $\pi:C/W_\lambda\rightarrow\Sigma$ 
(the kernel of the norm map). 
When $G=SO(n)$ or $Sp(2n)$, then 
$C/W_\lambda$ is $\widetilde{W}_\lambda/W_\lambda$-invariant. 
In the $SO(2n+1)$ and $Sp(2n)$ cases, 
the generic $C/W_\lambda$ meets the fixed locus of 
$\widetilde{W}_\lambda/W_\lambda$ transversally along a ramification
point of $C/W_\lambda$ over its $\widetilde{W}_\lambda/W_\lambda$-quotient. 
The fiber of (\ref{eq-char}) is isogenous to the Prym variety of the
involution of $C/W_\lambda$ (via the element $\eta$ of 
$\widetilde{W}_\lambda/W_\lambda$). 
In the $SO(2n)$ case, the generic $C/W_\lambda$ has nodes along 
the fixed locus of $\eta$. 
$\widetilde{W}_\lambda/W_\lambda$ acts freely on its normalization. 
The fiber of (\ref{eq-char}) is isogenous to the Prym variety of the 
unramified double cover. See section \ref{sec-max-isotropic-grassmannian} 
for an application of lemma \ref{lemma-symplectic-covering}
in the $SO(2n)$ case.

\section{Infinitesimal study of $M(G,c)$} 
\label{sec-infinitesimal-study-of-M-infinity}

We study the Poisson geometry of 
the infinite dimensional moduli space 
$M(G,c)$ in Theorem \ref{thm-infinite-dimensional-Poisson-moduli}.
$M(G,c)$ is the inductive limit (union) of the 
finite dimensional moduli spaces $M(G,c,\rho,D)$, as $D$ varies over all
effective divisors on $\Sigma$. 
The discussion is slightly complicated by the fact, that the moduli spaces
$M(G,c,\rho,D)$ are singular in general. 

\subsection{An anti-symmetric tensor}
\label{sec-an-anti-symmetric-tensor-D}
Next, we identify the tangent space to the moduli space $M(G,c,\rho,D)$ 
in Theorem \ref{thm-infinite-dimensional-Poisson-moduli}.
Given a pair $(P,\varphi)$, we define the vector bundle 
$$ad(P,\varphi,\rho)$$ as follows:
Away from the singularities of $\varphi$, the bundle $ad(P,\varphi,\rho)$ 
is the pullback via 
$\varphi : \Sigma \hookrightarrow P(G)$ 
of the vertical tangent bundle of $P(G)$. 
Along the divisor of singularities of $\varphi$, we define 
$ad(P,\varphi,\rho)$ in terms of the faithful representation: 
It is the unique extension to a subbundle of $End(E)$,
where $E$ is the vector bundle associated to $P$ via $\rho$. 
(This definition breaks the symmetry between the poles of $\varphi$ and 
$\varphi^{-1}$ as meromorphic sections of $End(E)$.) 

\begin{example} 
\label{example-ad-P-phi-rho}
{\rm
If $G=GL(n)$ and $\rho$ is the standard representation, then
$ad(P,\varphi,\rho)=End(E)$. 
If $G=SL(n)$, $E$ a vector bundle with trivial determinant,
and $\varphi$  a meromorphic section of $End(E)$ with determinant $1$,
then $ad(E,\varphi,\rho)$ is the subsheaf of $End(E)$ of sections 
satisfying
\[
\{f \in End(E) \ : \ {\rm tr}(\varphi^{-1}f)=0\}. 
\]
If $\varphi^{-1}$ is a nowhere vanishing holomorphic section of $End(E)(D')$, 
then it defines a line subbundle $L$ of $End(E)$ isomorphic to 
$\StructureSheaf{\Sigma}(-D')$. 
$ad(E,\varphi,\rho)$ is the subbundle $L^\perp$ orthogonal with respect to the 
trace pairing. It is isomorphic to the dual of the quotient
$End(E)/L$. Thus, ${\rm deg}(ad(E,\varphi,\rho))={\rm deg}(L)=-{\rm deg}(D')$. 

Let us describe $ad(P,\varphi,\rho)$ in the case where 
$G$ is one of the classical groups $Sp(2n)$ or $SO(n)$, 
$P$ a principal $G$-bundle, and 
$E$ the associated vector bundle via the standard 
representation. $E$ is endowed with a non-degenerate 
bilinear form $J : E \rightarrow E^*$
(anti-symmetric or symmetric accordingly). Then 
$P(G)$ is the subsheaf of $End(E)$ of 
invertible sections satisfying 
$
J \ \ = \ \ \varphi^*\circ J\circ \varphi.
$
Given a meromorphic section $\varphi$ of $P(G)$, the subbundle 
$ad(P,\varphi,\rho)$ of $End(E)$ corresponds to the subsheaf of sections $\phi$ 
satisfying 
$$
\phi^*\circ J\circ \varphi + \varphi^*\circ J\circ \phi \ \ = \ \ 0.
$$
\EndProof
}
\end{example}

If $\bbi$ is the identity section of $P(G)$, then 
$ad(P,\bbi,\rho)$ is the adjoint Lie algebra bundle $P\gg$. 
We have two natural homomorphisms
$$
L_\varphi, R_\varphi \ \ : \ \ P\gg \ \ \longrightarrow \ \ 
ad(P,\varphi,\rho)\otimes
{\cal O}_{\Sigma}(D)
$$
corresponding to left and right multiplication by the meromorphic 
section $\varphi$ of the group bundle $P(G)$. 
Denote their difference by $ad(\varphi)$
$$
ad(\varphi) \ := \ L_\varphi-R_\varphi \ \ : \ \ 
P\gg \ \ \longrightarrow ad(P,\varphi,\rho)\otimes
{\cal O}_{\Sigma}(D).
$$
In terms of the associated vector bundle $E$
in the above example, $L_\varphi$ and $R_\varphi$ are induced by left and right
multiplication in $End(E)$, which map the subbundle $P\gg$ of $End(E)$ 
into the subbundle $ad(P,\varphi,\rho)(D)$ of $End(E)(D)$. 

Standard deformation theory shows, that 
the tangent space to $M(G,c,\rho,D)$ at $(P,\varphi)$ is naturally identified with 
the first hypercohomology of the complex (in degrees $0$ and $1$)
\begin{equation}
\label{eq-tangent-complex-of-M-D}
P\gg 
\ \ {\buildrel{ad_\varphi}\over {\longrightarrow}} \ \ 
ad(P,\varphi,\rho)\otimes
{\cal O}_{\Sigma}(D).
\end{equation}
The proof is a repetition of the arguments in \cite{Bo1,biswas-ramanan,M}.
The cotangent space $T^*_{(P,\varphi)}M(G,c,\rho,D)$ is given by the first 
hypercohomology of the complex (in degrees $0$ and $1$)
\begin{equation}
\label{eq-cotangent-complex-M-D}
ad(P,\varphi,\rho)^*(-D) 
\ \ {\buildrel{ad_\varphi^*}\over {\longrightarrow}} \ \ 
P\gg^*
\end{equation}
(here we used a trivialization of the canonical line bundle of $\Sigma$). 
Once again, we have a natural homomorphism $\Psi$ from the cotangent 
complex to the tangent complex. In degree $1$ the homomorphism $\Psi_1$
is the composition of the isomorphism $P\gg^*\cong P\gg$ with 
the homomorphism $L_\varphi$ from $P\gg$ to $ad(P,\varphi,\rho)(D)$. 
It is easier to define the dual of the homomorphism $\Psi_0$ in degree
zero. The negative $-\Psi_0^*$ of the dual is the composition of the 
isomorphism $P\gg^*\cong P\gg$ with the {\it right} multiplication 
homomorphism $R_\varphi$ from $P\gg$ to $ad(P,\varphi,\rho)(D)$. 

We need to check that $\Psi$ is a homomorphism of complexes.
It suffices to check it away from the singular divisor of $\varphi$.
Over this open subset of $\Sigma$, the invariant bilinear form on $\gg$ 
induces an isomorphism
\begin{equation}
\label{eq-the-dual-of-ad-phi-rho-is-the-one-for-phi-inverse}
ad(P,\varphi,\rho)^*\ \ \cong \ \  ad(P,\varphi^{-1},\rho).
\end{equation}
The composition, of $\Psi_0$ with the generic isomorphism
$ad(P,\varphi,\rho)^*(-D) \cong  ad(P,\varphi^{-1},\rho)(-D)$, is equal to 
the negative of left 
multiplication by $\varphi$. Similarly, $ad_\varphi^*$
becomes $-ad_\varphi$. The commutativity
$ad_g\circ \Psi_0=\Psi_1\circ ad_g^*$ follows.

We get an induced homomorphism from the Zariski cotangent space to
the tangent space on the level of first hypercohomology 
\begin{equation}
\label{eq-infinitesimal-Psi-M-D}
\Psi \ : \ T^*_{(P,\varphi)}M(G,c,\rho,D) \ \ \longrightarrow \ \ 
T_{(P,\varphi)}M(G,c,\rho,D). 
\end{equation}
If we interchange the roles of
left and right multiplication in the above construction, we get  
another homomorphism of complexes. 
The two homomorphism of complexes are homotopic (see section
\ref{sec-symplectic-structure}). 

\begin{rem}
{\rm
The moduli space $M(G,c,\rho,D)$ is singular in general.
(It is smooth, for $G=GL(n)$, being an open subset of the 
moduli space of simple Higgs bundles, which is smooth \cite{nitsure}).
Let us determine the dimension of the Zariski tangent space 
to $M(G,c,\rho,D)$ at $(P,\varphi)$. 
Denote the complex (\ref{eq-tangent-complex-of-M-D}) by $B_\bullet$.
Assume that $D$ is sufficiently large, so that 
$H^1(ad(P,\varphi,\rho)(D))$ vanishes. 
Then $\bbh^2(B_\bullet)$ vanishes as well. 
Denote the center of $\gg$ by $\gz$. If $(P,\varphi)$ is simple, then 
$\bbh^0(B_\bullet)=\gz$. 
The dimension of $T_{(P,\varphi)}M(G,c,\rho,D)$ is given by 
$\bbh^1(B_\bullet)$ 
\[
\begin{array}{ccl}
\dim T_{(P,\varphi)}M(G,c,\rho,D) & = & \dim(\gz)-\chi(B_\bullet) 
\\
& = & \dim(\gz)-\chi(P\gg)+\chi(ad(P,\varphi,\rho)(D)) 
\\
& = & \dim(\gz)+\chi(ad(P,\varphi,\rho))+\deg(D)\cdot\dim(\gg). 
\end{array}
\]
The Euler characteristic $\chi(ad(P,\varphi,\rho))$ 
depends on $\rho$ and on the singularity data $\orbit$ of 
$\varphi$. 
We see that the dimension, of the Zariski tangent space to 
$M(G,c,\rho,D)$, does depend on the singularity data $\orbit$.
Hence, we expect the moduli space $M(G,c,\rho,D)$ to be singular. 
If, for example, $G=SL(n)$ and $\rho$ is the standard
representation, then each moduli $M(SL(n),0,\orbit)$
can be contained in the smooth locus of 
$M(SL(n),0,\rho,D)$, for at most one choice of $D$ 
(The Euler characteristic 
is calculated in example \ref{example-ad-P-phi-rho}.)
}
\end{rem}

Given two divisors $D_1\subset D_2$ on $\Sigma$, we have a natural 
embedding 
\[
M(G,c,\rho,D_1) \ \ \ \longrightarrow \ \ \ M(G,c,\rho,D_2).
\]
We get a direct system of moduli spaces $M(G,c,\rho,D)$, 
with respect to the partial ordering of effective divisors on $\Sigma$. 
We denote the inductive limit by $M(G,c)$. Similarly, given a simple
pair $(P,\varphi)$ in $M(G,c)$, we get a direct system of
Zarisky tangent spaces
\[
\bbh^1(\Sigma,[P\gg \LongRightArrowOf{ad_\varphi} ad(P,\varphi,\rho)(D_1)])
\ \ \ \hookrightarrow \ \ \ 
\bbh^1(\Sigma,[P\gg \LongRightArrowOf{ad_\varphi} ad(P,\varphi,\rho)(D_2)])
\]
with injective homomorphisms. 
Above, the two divisors $D_1\subset D_2$ are assumed to contain the divisor
$D(\orbit,\rho)$ defined in (\ref{eq-polar-divisor}). 
The cotangent spaces fit into an inverse system with
surjective homomorphisms. 
Denote by $T_{(P,\varphi)}M(G,c)$ and $T^*_{(P,\varphi)}M(G,c)$
the corresponding direct and inverse limits. 
We get a natural homomorphism
\begin{equation}
\label{eq-Psi-infinity}
\Psi_\infty \ : \  T^*_{(P,\varphi)}M(G,c) \ \ \ \longrightarrow \ \ \ 
T_{(P,\varphi)}M(G,c),
\end{equation}
which has finite rank.
If $(P,\varphi)$ belongs to $M(G,c,\rho,D)$, then $\Psi_\infty$
factors through (\ref{eq-infinitesimal-Psi-M-D}) via the natural 
projection $T^*_{(P,\varphi)}M(G,c)\rightarrow T^*_{(P,\varphi)}M(G,c,\rho,D)$
and the natural inclusion 
$T_{(P,\varphi)}M(G,c,\rho,D)\hookrightarrow T_{(P,\varphi)}M(G,c)$. 
We will see below that the homomorphism factors through
the cotangent and tangent spaces to $M(G,c,\orbit)$ (see Lemma 
\ref{lemma-symplectic-leaves-are-loci-with-fixed-singularity-type}).

The construction of the sheaf $ad(P,\varphi,\rho)$ does not behave well 
in families (it leads to a sheaf, which is not flat).
Nevertheless, 
the construction (of the pull-back) of the tangent and cotangent bundles 
of $M(G,c)$ can be carried out in the relative setting
of finite dimensional families of pairs. 
We indicate briefly this relative construction. 
First, we provide an alternative construction of $T_{(P,\varphi)}M(G,c)$. 
Assume 
that $(P,\varphi)$ is in $M(G,c,\rho,D_0)$. Choose effective divisors
$D_1$ and $D_2$, such that $Ad_{\varphi^{-1}}$ is a regular section
of $\End(P\gg)(D_1)$ and 
\begin{eqnarray*}
L_{\varphi^{-1}} \ : \ ad(P,\varphi,\rho)(D_0) & \longrightarrow & P\gg(D_1),
\\
L_{\varphi} \ \hspace{3ex} : \hspace{5ex} \ P\gg(D_1) & \longrightarrow & 
ad(P,\varphi,\rho)(D_2)
\end{eqnarray*}
are regular homomorphisms. Let
$B_\bullet(D)$ be the complex (\ref{eq-tangent-complex-of-M-D}). 
Let $A_\bullet(D)$ be the complex
\begin{equation}
\label{eq-complex-approaximating-T-M-infinity}
P\gg \ \LongRightArrowOf{\bbi-Ad_{\varphi^{-1}}} \ P\gg(D).
\end{equation}
We have the complex homomomorphisms
\begin{eqnarray*}
(\bbi,L_{\varphi^{-1}}) \ : \ B_\bullet(D_0) & \longrightarrow & A_\bullet(D_1)
\\
(\bbi,L_\varphi) \ : \ A_\bullet(D_1) & \longrightarrow & B_\bullet(D_2). 
\end{eqnarray*}
It is easy to see that both complex homomorphisms induce injective 
homomorphisms on the level of first hyper-cohomologies
\[
\bbh^1(B_\bullet(D_0)) \ \subset \ \bbh^1(A_\bullet(D_1))  \ \subset \ 
\bbh^1(B_\bullet(D_2)). 
\]
Consider the directed system of first hyper-cohomologies of 
(\ref{eq-complex-approaximating-T-M-infinity}),
where $D$ varies over all effective divisors, which contain $D_1$. 
The above comparison implies, that there is a natural isomorphism 
between the direct limits of $\bbh^1(A_\bullet(D))$ and 
$\bbh^1(B_\bullet(D))$. 

The complex (\ref{eq-complex-approaximating-T-M-infinity}) 
behaves well in families and the direct limit
provides the relative construction
of the pull-back of $TM(G,c)$ to finite dimensional families of pairs
parametrized by $M(G,c,\rho,D)$. Similarly, we get the pull-back of $T^*M(G,c)$
and the relative version of the homomorphism (\ref{eq-Psi-infinity}).
Over the smooth locus of 
$M(G,c,\rho,D)$, it gives rise to a global homomorphism from the 
cotangent to the tangent bundles:
\begin{equation}
\label{eq-Psi-M-D}
\Psi \ : \ T^*M(G,c,\rho,D) \ \ \longrightarrow \ \ TM(G,c,\rho,D).
\end{equation}

\begin{thm}
\label{thm-Lambda-is-a-Poisson-structure}
$\Psi$ defines a Poisson structure on $M(G,c)$ and on the smooth
locus of $M(G,c,\rho,D)$. 
\end{thm}

\noindent
{\bf Proof:}
$\Psi$ is anti-symmetric because 
$\Delta_0$ is homotopic to $-(\Delta_1^*)$. This is proven 
as in Theorem \ref{thm-algebraic-2-form}.
The Jacobi identity follows from the involutivity of the foliation induced by
the image of $\Psi$ (Lemma 
\ref{lemma-symplectic-leaves-are-loci-with-fixed-singularity-type})
and the Jacobi identity for the leaves $M_\Sigma(G,c,\orbit)$
(Theorem \ref{thm-algebraic-2-form}). 
\EndProof

\subsection{The symplectic leaves foliation}
\label{sec-symplectic-leaves-folliation}

In this section we prove that the symplectic 
leaves of $M(G,c,\rho,D)$ 
are precisely the moduli spaces $M(G,c,\orbit)$,  whose polar divisor
$D(\orbit,\rho)$ is contained in $D$. See (\ref{eq-polar-divisor})
for the definition of $D(\orbit,\rho)$. 

\begin{lemma}
\label{lemma-symplectic-leaves-are-loci-with-fixed-singularity-type}
$M(G,c,\orbit)$ 
is a finite union of symplectic leaves of $M(G,c,\rho,D)$. 
\end{lemma}

Proof: 
It suffices to show that the Zariski tangent space to $M(G,c,\orbit)$ 
is equal 
to the image of the homomorphism $\Psi$ induced by the Poisson structure. 
We have the following short exact sequence of (column) complexes:

\[
\begin{array}{rcccl}
\left[\!\!
\begin{array}{c}
ad(P,\varphi,\rho)^*(\!-\!D) \\ \downarrow \ ad_\varphi^* \\
P\gg^* 
\end{array}
\!\!\right]
& 
\begin{array}{c}
{\buildrel{R_\varphi^*}\over {\longrightarrow}} \\ \\
{\buildrel{L_\varphi}\over {\longrightarrow}} 
\end{array}
& 
\left[\!\!
\begin{array}{c}
P\gg \\ \downarrow \ ad_\varphi \\ ad(P,\varphi,\rho)(D) 
\end{array}
\!\!\right]
& 
\begin{array}{c}
{\buildrel{q_0}\over {\longrightarrow}} \\ \\
{\buildrel{q_1}\over {\longrightarrow}} 
\end{array}
& 
\left[\!\!
\begin{array}{c}
P\gg/Im(R_\varphi^*) \\ \downarrow \\ 
\frac{ad(P,\varphi,\rho)(D)}{Im(L_\varphi)}
\end{array} 
\!\!\right]
\end{array}
\]


\noindent
which we denote, for short, by 
$$
A_\bullet \ \ \longrightarrow \ \ B_\bullet \ \ \longrightarrow \ \ 
C_\bullet.
$$
We get the long exact sequence of cohomologies:
\begin{equation}
\label{4.50}
\begin{array}{rcccccl}
& 0 & \longrightarrow & \bbh^0(B_\bullet) & \longrightarrow & 
\bbh^0(C_\bullet) & \rightarrow
\\
\rightarrow & T^*_{(P,\varphi)}M(G,c,\rho,D) & 
{\buildrel{\Psi}\over {\longrightarrow}} & T_{(P,\varphi)}M(G,c,\rho,D) & 
\longrightarrow & \bbh^1(C_\bullet) & \rightarrow 
\\
\rightarrow & \bbh^2(A_\bullet) & \longrightarrow & \bbh^2(B_\bullet) &
\longrightarrow & 0.
\end{array}
\end{equation}
The group $\bbh^1(C_\bullet)$ is the quotient of 
$ad(P,\varphi,\rho)(D)$ by the image of the homomorphism
$$
L_\varphi+ad_\varphi \ \ : \ \ P\gg\oplus P\gg \ \ 
\longrightarrow \ \ ad(P,\varphi,\rho)(D).
$$
The image is equal to that of the homomorphism 
$$
L_\varphi+R_\varphi \ \ : \ \ P\gg\oplus P\gg \ \ 
\longrightarrow \ \ ad(P,\varphi,\rho)(D).
$$
The image of the latter homomorphism is precisely $ad(P,\varphi)$. 
We conclude that the kernel of
$\bbh^1(B_\bullet) \rightarrow \bbh^1(C_\bullet)$ is precisely the
Zariski tangent space to $M(G,c,\orbit)$. 
The Lemma follows from the exactness of the sequence (\ref{4.50}). 
\EndProof

\section{Poisson Hecke correspondences}
\label{sec-Hecke-correspondences}

Given a Poisson group, Drinfeld constructed compatible Poisson
structures on its homogeneous spaces \cite{drinfeld-poisson-groups}.
We introduce next a Poisson structure on the moduli spaces 
$\Hecke(G,c_1,c_2,\orbit)$, which 
are analogous to $G[[t]]$-orbits in the loop Grassmannian $G((t))/G[[t]]$. 
The moduli space $\Hecke(G,c_1,c_2,\orbit)$ 
parametrizes isomorphism classes of simple  
triples $(P_1,P_2,\varphi)$, where $P_i$ is a principal $G$
bundle of topological type $c_i$ and $\varphi:P_1\rightarrow P_2$
is a meromorphic, $G$-equivariant, isomorphism (a meromorphic section of
$\Isom(P_1,P_2)$). The singularities of $\varphi$ are assumed to 
be in $\orbit=\sum_{p\in \Sigma}\orbit_p$. 
We will see that the locus $\Hecke(G,P_1,P_2,\orbit)$, 
with a fixed pair of bundles $(P_1,P_2)$, is a symplectic leaf.
Note that the orbit of the triple $(P_1,P_2,\varphi)$, under 
$\Aut(P_1)\times \Aut(P_2)$, is contained in a single isomorphism class.
The triple is {\em simple}, 
if the stabilizer of $\varphi$ in $\Aut(P_1)\times \Aut(P_2)$ is 
the diagonal embedding of the center of $G$. 
It is {\em infinitesimally simple} if the above statement holds on the level
of adjoint Lie algebra bundles.

\begin{rem}
\label{rem-non-emptiness-condition-for-Hecke}
{\rm
\begin{enumerate}
\item
The moduli spaces $\Hecke(G,c_1,c_2,\orbit)$ are related to the Hecke 
correspondences in the Geometric Langlands program \cite{beilinson-drinfeld}.
\item
Given a triple $(P_1,P_2,\varphi)$, 
the topological type $c_2$ of $P_2$ is determined by the type $c_1$ 
of $P_1$ and the singularity data 
$\orbit$ of $\varphi$. 
This leads to a necessary compatibility condition for 
$(c_1,c_2,\orbit)$ in order for $\Hecke(G,c_1,c_2,\orbit)$ to be non-empty. 
For each singularity point $p\in \Sigma$, choose a co-character $a_p$ 
in the $W$-orbit determined by $\orbit_p$. Then the sum 
$\sum_{p\in \Sigma}a_p$ projects, via
(\ref{eq-homomorphism-from-co-characters-to-pi-1-G}), 
to the difference $c_2-c_1$ in $\pi_1(G)$. 

Proof of the condition: 
Let $\U$ be the open covering of $\Sigma$, consisting of the complement 
$U_0$ of the singular points of $\varphi$, and small disks $U_p$,
one for each singular point. Choose an isomorphism 
$\eta_p:(P_1\restricted{)}{U_{p}}\rightarrow (P_2\restricted{)}{U_{p}}$
between the restrictions of the bundles to $U_p$.
The cohomology class of the $1$-cocycle $g_{p,0}:= \eta_p^{-1}\circ\varphi$ in
$Z^1(\U,P_1(G))$, represents the isomorphism class of the bundle $P_2$. 
The connecting homomorphism, of the long exact cohomology sequence
(\ref{eq-long-exact-seq-of-analytic-sheaf-coho}), with $P=P_1$, takes 
this $1$-cocycle to the class $c_2-c_1$ in $H^2(\U,\pi_1(G))$.
On the other hand, under the identification of $H^2(\Sigma,\pi_1(G))$
with $\pi_1(G)$, the connecting homomorphism takes the cocycle $(g_{p,0})$ 
to the image of $\orbit$ under 
(\ref{eq-homomorphism-from-co-characters-to-pi-1-G}). 
\EndProof
\end{enumerate}
}
\end{rem}

\subsection{Comparison between $\Hecke(G,c,c,\orbit)$ and $M(G,c,\orbit)$}
\label{sec-dimensions-of-Hecke-and-M}

We first relate the moduli space $\Hecke(G,c_1,c_2,\orbit)$ to products 
of orbits in the loop Grassmannian. Then we relate the moduli spaces 
$\Hecke(G,c,c,\orbit)$ and $M(G,c,\orbit)$. Combining the two relations, 
we will get a conceptual explanation of the formula for the dimension 
of $M(G,c,\orbit)$ in Corollary \ref{cor-dimension-of-a-symplectic-leaf}.

Consider the locus $\Hecke(G,c_1,P_2,\orbit)$ in 
$\Hecke(G,c_1,c_2,\orbit)$, with a  fixed bundle $P_2$. 

\begin{claim} 
\label{claim-Hecke-with-fixed-P-2}
$\Hecke(G,c_1,P_2,\orbit)$ 
is isomorphic to a Zariski open subset in the $\Aut(P_2)$-quotient of the
product of the $G[[t]]$-orbits in
$G((t))/G[[t]]$, parametrized by the data $\orbit$. 
\end{claim}

\noindent
{\bf Proof:} The bundle $P_2$
determines a bundle $B$ of loop grassmannians over $\Sigma$.
Using the notation of section
\ref{sec-main-results}, the fiber of $B$ over $p\in\Sigma$ is 
the quotient 
$P_2(G(\CompletedFunctionField{(p)}))/
P_2(G(\CompletedStructureSheaf{(p)}))$. 
A triple $(P_1,P_2,\varphi)$ determines a 
section $\sigma$ of $B$ as follows. Choose
a local isomorphism $\eta:P_2\rightarrow P_1$.
Then $\varphi\circ\eta$ is a local meromorphic section of $P_2(G)$.
A different $\eta$ yields the same section $\sigma$ of $B$.
Note that $B$ has a natural section, corresponding to the image
of $P_2(G)$, and $\sigma$ coincides with this section away from 
the singularities of $\varphi$. 
Suppose $(P_1',P_2,\varphi')$ is another triple, whose section $\sigma'$
is equal to $\sigma$. Then there is a global isomorphism 
$f:P_1\rightarrow P_1'$ such that $\varphi'\circ f=\varphi$
(the equality $\sigma=\sigma'$ implies that $f$ exists locally, but its 
uniqueness implies that it is global). 

Conversely, given a section $\sigma$, we can lift it to a triple
$(P_1,P_2,\varphi)$ as follows. For every point $p$, with 
$\orbit_p\neq G[[t]]$, choose an analytic (or formal) neighborhood $U_p$ of 
$p$. For a sufficiently small $U_p$,  
we can lift $\sigma$ to a meromorphic section of $P(G)$, which is holomorphic
on the punctured neighborhood $U_p\setminus\{p\}$. 
Over the complement in $\Sigma$, of the singularity divisor $S$ of $\orbit$,
we choose the section $1$ of $P(G)$. The open subsets $U_p$ and $U_q$
are assumed to be disjoint, for any two distinct points $p$ and $q$ in $S$. 
We get a \v{C}eck cocycle
representing a class in $H^1(\Sigma,P(G))$, i.e., an isomorphism
class of a principal $G$-bundle $P_1$. It comes with
a natural identification $\varphi$ with $P_2$ over $\Sigma\setminus S$. 
By construction, the triple $(P_1,P_2,\varphi)$ has singularities in 
$\orbit$. 
For $\sigma$ in an open subset, the triple would be simple. 
The claim follows.
\EndProof

\medskip
For a semi-stable and regular bundle $P_2$, the 
dimension of $\Aut(P_2)$ is equal to the dimension of the moduli
of deformations of $P_2$. If $G$ is semi-simple, 
it follows from the above claim, that the expected dimension of 
$\Hecke(G,c_1,c_2,\orbit)$ is equal to the sum of the dimensions of the
$G[[t]]$-orbits in $G((t))/G[[t]]$ parametrized by the data $\orbit$. 
The dimension of each such orbit is given by formula 
(\ref{eq-gamma-is-dot-product-with-delta}) (see \cite{lusztig}). 
If $G$ is reductive with a positive dimensional center, then 
the action of $\Aut(P_2)$, on each orbit in the loop grassmannian, 
factors through the quotient $\Aut(P_2)/Z$ by the center of $G$. 
Hence, the expected dimension of $\Hecke(G,c_1,c_2,\orbit)$ is larger
(see Theorem \ref{thm-Poisson-structure-on-Hecke}).

When $G=PGL(n)$ and $c_1=c_2=c$ is a generator of $\pi_1(G)$, then 
$M(G,c,\orbit)$ admits a rational symplectic morphism 
onto an open symplectic subset of $\Hecke(G,c,c,\orbit)$, 
because if $P_1$ and $P_2$ are 
stable, then they are isomorphic (compare with Lemma
\ref{lemma-symplectic-leaf-is-birational-to-product-of-flag-varieties}). 
For example, if the data $\orbit$ is non-singular at every point
(and we allow infinitesimally simple, but non-simple objects), then 
$\Hecke(G,c,c,\orbit)$ is a single point, while $M(G,c,\orbit)$ is 
the commutative group of global automorphisms of the unique stable 
$PGL(n)$ bundle of topological type $c$ 
(see Example \ref{example-PGL-n-moduli-is-disconnected}).

In general, we have two rational morphisms
\[
\Hecke(G,c,c,\orbit) \LongLeftArrowOf{\pi_1} 
\ M(G,c,\orbit) \ \LongRightArrowOf{\pi_2} 
U(G,c),
\]
where $U(G,c)$ is the moduli space of semi-stable $G$-bundles of topological 
type $c$, endowed with the trivial Poisson structure. 
The image of $\pi_1$
consists of triples $(P_1,P_2,\varphi)$, with isomorphic $P_1$ and $P_2$. 
The fiber of $\pi_1$ through a simple $(P,\varphi)$ is isomorphic to the 
quotient of 
$\Aut(P)\cdot\varphi\cdot\Aut(P)$ by the conjugation action of $\Aut(P)$. 
Note, that the dimensions of $H^0(P\gg)$ and $H^1(P\gg)$ are equal. Hence, 
if $G$ is semi-simple, then 
the expected dimensions of $M(G,c,\orbit)$ and $\Hecke(G,c,c,\orbit)$
are equal. {\em This explains the equality between the dimension of 
$M(G,c,\orbit)$ (Corollary \ref{cor-dimension-of-a-symplectic-leaf}), 
for semi-simple $G$, and the sum of dimensions of orbits in the loop 
Grassmannian}. If $G$ is reductive with a positive dimensional center, 
we will see that the deformations of the pair $(P_1,P_2)$, arising
from triples in $\Hecke(G,c,c,\orbit)$, are restricted. 
This restriction is the algebraic analogue, 
on the level of rational equivalence, of the topological condition
in Remark \ref{rem-non-emptiness-condition-for-Hecke}. 
This restriction leads to a difference
between the dimensions of the two moduli spaces.

\subsection{A Poisson structure on $\Hecke(G,c_1,c_2,\orbit)$}

The infinitesimal study of $\Hecke(G,c_1,c_2,\orbit)$ is similar 
to that of $M(G,c,\orbit)$. 
We will see that it is in fact simpler, since the 
infinitesimal deformations of the triples can be given by the first 
cohomology of a vector bundle (rather than a complex). 
Let $\Delta_\varphi\subset P_1\gg\oplus P_2\gg$ be the
subbundle corresponding to the subsheaf
\[
\{(a,b)\ : \ Ad_{\varphi}(a)+b=0
\}.
\]
Define $ad(P_1,P_2,\varphi)$ to be the quotient in the short exact sequence 
\begin{equation}
\label{eq-short-exact-seq-defining-ad-P1-P2-phi}
0\rightarrow \Delta_\varphi \RightArrowOf{e} 
P_1\gg\oplus P_2\gg \RightArrowOf{q}
ad(P_1,P_2,\varphi) \rightarrow 0.
\end{equation}
Set $d:=q$ 
\begin{equation}
\label{eq-tangent-complex-of-Hecke}
P_1\gg\oplus P_2\gg \ \ \LongRightArrowOf{d} \ \ 
ad(P_1,P_2,\varphi). 
\end{equation}
Lemmas \ref{lemma-sections-of-ad-P-varphi-are-infinitesimal-def}
and \ref{lemma-ad-P-phi-is-dual-to-Delta-phi} hold in the more
general setting, with $ad(P_1,P_2,\varphi)$ replacing $ad(P,\varphi)$
and $d$ replacing $ad_\varphi$.
The first order infinitesimal deformations of a triple are
parametrized by the  first hyper-cohomology of the complex 
(\ref{eq-tangent-complex-of-Hecke}).
The Zariski cotangent space, at a simple triple, is the first hyper-cohomology
of the dual complex
\begin{equation}
\label{eq-cotangent-complex-of-Hecke}
ad(P_1,P_2,\varphi)^* \ \ \LongRightArrowOf{d^*} \ \ 
[P_1\gg\oplus P_2\gg]^*.
\end{equation}

Polishchuk defined  a Poisson structure on
$\Hecke(G,c_1,c_2,\orbit)$ in the $GL(n)$ case \cite{Po}.
He then generalized his construction for moduli spaces of pairs $(P,s)$, 
where $P$ is a principal bundle, with a reductive structure group, and 
$s$ is a section of a vector bundle associated to $P$ via a
suitable representation (Theorem 6.1 in \cite{Po}). We get a moduli of triples,
once we choose the structure group of $P$ to be $G\times G$. 
However, Polishchuk's condition on the representation of $G\times G$
requires to replace a simple group $G$ by some central extension.
Here we use a slightly different generalization of Polishchuk's construction
to exhibit a natural Poisson structure on $\Hecke(G,c_1,c_2,\orbit)$, 
for arbitrary reductive $G$. 
Consider the homomorphism $\psi=(\psi_0,\psi_1)$ 
from (\ref{eq-cotangent-complex-of-Hecke}) to 
(\ref{eq-tangent-complex-of-Hecke}), which is zero in degree $1$, 
and in degree $0$ is given by the composition 
\[
ad(P_1,P_2,\varphi)^* \LongRightArrowOf{d^*} [P_1\gg\oplus P_2\gg]^*
\LongRightArrowOf{(\kappa,-\kappa)^{-1}} 
P_1\gg\oplus P_2\gg.
\]
Lemma \ref{lemma-ad-P-phi-is-dual-to-Delta-phi} implies the equality
\[
d\circ (\kappa,-\kappa)^{-1}\circ d^* \ \ = \ \ 0.
\]
Consequently, 
$\psi$ is indeed a complex homomorphism. 
The  degree $1$ sheaf in (\ref{eq-cotangent-complex-of-Hecke}) 
is isomorphic to the degree $0$ sheaf in (\ref{eq-tangent-complex-of-Hecke}).
A choice of an isomorphism 
$h: [P_1\gg\oplus P_2\gg]^* \rightarrow P_1\gg\oplus P_2\gg$,
gives rise to a 
homotopic complex homomorphism
$\psi+(h\circ d) + (d\circ h)=
(\psi_0+h\circ d^*,d\circ h)$. Hence, we could have chosen 
$\psi_0$ to vanish instead. Denote by
\begin{equation}
\label{eq-Poisson-str-on-Hecke}
\Psi \ : \ T^*_{(P_1,P_2,\varphi)}\Hecke \ \ \rightarrow \ \
T_{(P_1,P_2,\varphi)}\Hecke
\end{equation}
the homomorphism induced by $\psi$ on the level of first hyper-cohomologies. 
$\Psi$ is anti-self dual, because $\psi^*$ is homotopic to 
$-\psi$ via the homotopy $h=(\kappa,-\kappa)^{-1}$ ($\kappa$ as in
Lemma \ref{lemma-ad-P-phi-is-dual-to-Delta-phi}).
Indeed, $\psi_0=h\circ(d^*)$ and, since $h$ is self-dual, 
$\psi_0^*=d\circ h$. Thus, 
$\psi+\psi^*=h\circ d^*+d\circ h$ and
$\psi^*$ is homotopic to $-\psi$. 

Note, that $d$ is surjective, and its kernel is $(id,-id)(\Delta_\varphi)$. 
Consequently, (\ref{eq-tangent-complex-of-Hecke})
is quasi-isomorphic to $\Delta_\varphi$ and
(\ref{eq-cotangent-complex-of-Hecke}) is quasi-isomorphic to the complex 
$ad(P_1,P_2,\varphi)[-1]$ (with a single sheaf
$ad(P_1,P_2,\varphi)$ in degree $1$). 
Identify the first hyper-cohomology of 
(\ref{eq-tangent-complex-of-Hecke}) with $H^1(\Delta_\varphi)$ and
the first hyper-cohomology of (\ref{eq-cotangent-complex-of-Hecke}) 
with $H^0(ad(P_1,P_2,\varphi))$. Then the homomorphism $\Psi$ 
is the connecting homomorphism of the long exact sequence
\[
\begin{array}{ccccccl}
0 \ \rightarrow \ H^0(\Delta_\varphi) 
&\rightarrow  &H^0(P_1\gg\oplus P_2\gg) &\rightarrow &
T^*_{(P_1,P_2,\varphi)} \Hecke  &\LongRightArrowOf{\Psi}
\\
T_{(P_1,P_2,\varphi)}\Hecke  & \LongRightArrowOf{df}  &
H^1(P_1\gg\oplus P_2\gg) & \rightarrow  & 
H^1(ad(P_1,P_2,\varphi)) &\rightarrow  &0,
\end{array}
\] 
arising from 
the short exact sequence (\ref{eq-short-exact-seq-defining-ad-P1-P2-phi}). 
If the triple $(P_1,P_2,\varphi)$ is infinitesimally simple,
then $H^0(\Delta_\varphi)$ is isomorphic to $\gz$ and  
$H^1(ad(P_1,P_2,\varphi))$ to $\gz^*$ (by Serre's Duality). 
The homomorphism $df$ above sends an infinitesimal deformation of the triple
$(P_1,P_2,\varphi)$ to the infinitesimal deformation of the pair $(P_1,P_2)$.
It follow that the symplectic leaves foliation of $\Psi$ is induced by fixing 
the isomorphism class of the pair $(P_1,P_2)$.
Moreover, the deformations of the pair $(P_1,P_2)$, arising from
triples in $\Hecke(G,c_1,c_2,\orbit)$, are constrained to a subspace, whose
codimension is equal to the dimension of $\gz^*$. 

\begin{thm}
\label{thm-Poisson-structure-on-Hecke}
\begin{enumerate}
\item
\label{thm-item-dimension-of-Hecke}
The moduli space $\Hecke_\Sigma(G,c_1,c_2,\orbit)$ is either empty, or smooth  
of dimension 
\[
\dim(\gz)+\sum_{p\in\Sigma}\gamma(\orbit_p),
\]
where $\gamma(\orbit_p)$ is given in 
(\ref{eq-gamma-is-dot-product-with-delta}). 
\item
\label{thm-item-Hecke-is-Poisson}
$\Hecke_\Sigma(G,c_1,c_2,\orbit)$ admits a natural Poisson structure $\Psi$. 
Its symplectic leaves are
the connected components of the loci $\Hecke_\Sigma(G,P_1,P_2,\orbit)$, 
in which the isomorphism class of the pair $(P_1,P_2)$ is fixed. 
\end{enumerate}
\end{thm}

\noindent
{\bf Proof:} 
\ref{thm-item-dimension-of-Hecke}) 
Let $(P_1,P_2,\varphi)$ be a simple triple. 
Simplicity of the triple implies, that $H^0$ of the tangent complex
(\ref{eq-tangent-complex-of-Hecke}) is isomorphic to $\gz$. 
$H^2$ vanishes, because the complex is quasi-isomorphic to 
$\Delta_\varphi$. The Euler characteristic of $\Delta_\varphi$
is calculated, up to a change of sign, in Lemma
\ref{lemma-degree-of-ad-P-phi}. The dimension of the Zariski tangent space
is given by $h^1(\Delta_\varphi)$, which is equal to 
$h^0(\Delta_\varphi)-\chi(\Delta_\varphi)$. 
The proof of smoothness of $\Hecke_\Sigma(G,c_1,c_2,\orbit)$ is 
analogous to the proof of Theorem \ref{thm-smoothness}.

\ref{thm-item-Hecke-is-Poisson})
It remains to prove the Jacobi identity for the Poisson bracket. 
The Jacobi identity can be expressed, in terms of
the homomorphism (\ref{eq-Poisson-str-on-Hecke}), by the vanishing
\begin{equation}
\label{eq-jacobi-identity}
\Psi(\omega_1)\cdot \langle\Psi(\omega_2),\omega_3\rangle -
\langle[\Psi(\omega_1),\Psi(\omega_2)],\omega_3\rangle + 
cp(1,2,3) \ \ = \ \ 0,
\end{equation}
where $\omega_i$ are local one-forms on
$\Hecke_\Sigma(G,c_1,c_2,\orbit)$, $[\cdot,\cdot]$ is the Lie bracket
of vector fields, and $cp(1,2,3)$ indicates terms obtained by
cyclic permutations of $1$, $2$ and $3$ from the first two terms. 
The proof is identical to that of Theorem 6.1 in \cite{Po}. 
One uses the set-up of \cite{Bo1,Bo2}. 
A lengthy cohomological calculation reduces 
the equality (\ref{eq-jacobi-identity})  to the 
vanishing of the tri-linear tensor
\begin{equation}
\label{eq-vanishing-of-tri-linear-tensor}
\langle x,[y,z]\rangle, \ \ \ x,y,z \in \Delta_\varphi, 
\end{equation}
where the Lie-bracket is the restriction of the one on $P_1\gg\oplus P_2\gg$
and the pairing $\langle\bullet,\bullet\rangle$ on $P_1\gg\oplus P_2\gg$
is $(\kappa,\kappa)$. 
The vanishing (\ref{eq-vanishing-of-tri-linear-tensor})
follows immediately from the definition of
$\Delta_\varphi$. 
\EndProof

\section{The Classical Dynamical Yang-Baxter equations}
\label{sec-dynamical-yang-baxter}

We show that the Poisson 
bracket of Etingof and Varchenko \cite{etingof-varchenko}, 
defined from the dynamical $r$-matrices, corresponds to the 
one 
on the moduli space of pairs $M_\Sigma(G,c)$, which was 
defined in Theorem \ref{thm-infinite-dimensional-Poisson-moduli} 
above. The comparison is proven first for simply connected 
groups (Theorem 
\ref{thm-comparison-with-cdybe-simply-connected-case}). 
We then consider extensions to complex reductive groups. 
We write our   elliptic curve as
\[
\Sigma \ \ := \ \ \bbc/(\bbz + \tau\bbz) \ \ = \ \ \bbc^\times/(q^\bbz),
\]
where $q=e^{2\pi i\tau}$.
We will use a description of  principal $G$-bundles
on $\Sigma$ in terms of flat connections (see, e.g., 
Atiyah and Bott \cite {atiyah-bott}); alternately, this amounts to a 
description of the bundle in terms of {\em automorphy factors}.
Indeed, the pull back of a $G$-bundle $P$ to  $\bbc^\times$ is trivial.
$P$ can be described as a quotient of $\bbc^\times \times G$, by a 
diagonal 
action of $\bbz$. The generators of $\bbz^2$ act on $G$ via 
an element of $G$, the  automorphy factor.
It can be normalized, up to a discrete choice. This language is suitable 
for 
the comparison, because the factor $V$, in the Poisson groupoid 
(\ref{eq-groupoid}), parametrizes an open set of the principal bundles 
together with a choice 
of 
a factor of automorphy. 
A meromorphic section $\varphi$, of the adjoint group bundle $P(G)$,
pulls back to a $G$-valued function on $\bbc^\times$, hence an element of 
the
loop group. 
This will enable us to describe 
an open subset in the moduli space $M_\Sigma(G,c)$, as a symplectic 
(or rather Poisson) reduction of the groupoid (\ref{eq-groupoid})
of Etingof and Varchenko.

\subsection{Some elliptic functions}

The elliptic Jacobi theta-function $\theta_1(z)$, for the curve $\Sigma$ 
with modulus $\tau$,  satisfies the relations
\[
\theta_1(z+1) = -\theta_1(z),\quad \theta_1(z+\tau)=
\theta_1(z) e^{-2\pi i z} e^{-\pi i (\tau+1)}.
\]
It has a single zero on each fundamental domain, at the translates of
the origin. If we set, following Etingof and Varchenko 
\cite{etingof-varchenko},
\[
\sigma_w(z) = {\theta_1(w-z) \theta_1'(0)\over
\theta_1(z) 
\theta_1(w)},
\]
we have the periodicity relations
\[
\sigma_w(z+1) = \sigma_w(z),\quad \sigma_w(z+\tau) =
\sigma_w(z) 
e^{2\pi i w}.
\]
This function has a single pole on each fundamental domain, with residue
1, at the translates of the origin, and a single zero at the translates 
of the point $w$. If $L_w$ is the line bundle corresponding to the divisor 
$w-0, w\in \bbc$, $\sigma_w$ represents a section of  $L_w$ in 
the trivialisation associated to the automorphy factors $1, e^{2\pi i w}$.

Let us cover the elliptic curve by an open disk $U_+$ around the origin 
and by $U_- = \Sigma-\{0\}$.
If the bundle $L_{ -w}$ is non trivial, then $H^1(\Sigma,L_{ -w})=0$, and 
so its sections 
${\cal L}_{-w}(U_+\cap U_-)$, over the 
punctured neighbourhood $(U_+\cap U_-)$ of the origin, split as 
\begin{equation}
\label{eq-splitting-two-summands}
{\cal L}_{-w}
(U_+\cap U_-) = {\cal L}_{-w}(U_+)\oplus {\cal L}_{-w}(U_-).
\end{equation}
Let $f$ represent a section of 
${\cal L}_{-w}(U_+\cap U_-)$, 
in the trivialisation associated to the automorphy factors
$1, e^{-2\pi i w}$. 
We have, for $z'\in U_+$,
\[
\oint \sigma_w(z-z') f(z) dz \ \ = \ \
\left\{ 
\begin{array}{ccc}
f(z') & {\rm if} & f\in {\cal L}_{-w}(U_+),
\\
0    & {\rm if} & f\in {\cal L}_{-w}(U_-),
\end{array}
\right.
\]
where the contour is around the boundary of $U_+$.
For example, if $f\in {\cal L}_{-w}(U_-)$, then 
$\sigma_w(z-z') f(z)$ represents a section of the trivial bundle,
which is holomorphic on the complement of $U_+$, and thus has residue $0$ 
at $z'$.
In short, $\sigma_w(z-z')$ is  the integral kernel of the projection to
${\cal L}_{-w}(U_+)$.

If one now considers in turn the function
\[
\rho(z) = {\theta_1'(z)\over \theta_1(z)},
\]
we find that it satisfies the periodicity relations 
\[
\rho(z+1) = \rho(z),\quad \rho(z+\tau) = \rho(z) - 2 \pi i. 
\]
In addition, it is an odd function.
The function $\rho$ plays the same role as $\sigma_w$, but for the trivial 
bundle. One can split the space of functions ${\cal O}(U_+\cap U_-)$ as 
\begin{equation}
\label{eq-three-summands}
{\cal O}(U_+\cap U_-) \ \ = \ \  {\cal O}(U_+ ) \ \ \oplus  \ \ 
\bbc{\rho(z)} \ \ 
\oplus \ \ \widetilde{\cal O}(U_-).
\end{equation}
${\cal O}(U_+ )$ and ${\cal O}(U_- )$  represent the spaces of holomorphic 
functions on $U_+$ and $U_-$ respectively. We exclude the constants in 
${\cal O}(U_- )$ by setting
$\widetilde{\cal O}(U_-)$ to be the
subspace of functions in ${\cal O}(U_-)$
satisfying
\[
\int_\epsilon^{1+\epsilon} f dz = 0.
\]
Let $P_+$, $P_0$, and $P_-$ denote the three projections onto the factors 
of
(\ref{eq-three-summands}).
One finds:
\[
\oint \rho(z-z') f(z)dz \ \ =  \ \  
\left\{
\begin{array}{ccc}
f(z') & {\rm if} & f\in  {\cal O}(U_+),
\\
0 &    {\rm if} & f(z) = {\rho(z)},
\\ 
0 &    {\rm if} &  f\in  \tilde{\cal O}(U_-).
\end{array}
\right.
\]
Consequently, $\rho(z-z') $ is the integral kernel of the projection 
$P_+$.

\subsection{$G$-bundles for a semi-simple simply connected group}

Now let us move to a complex reductive  Lie algebra $\gg$, with Cartan 
subalgebra $\gt$, corresponding to a 
simply connected group $G$ and Cartan subgroup $T$. 
Let $x_i, i = 1,..,r$ be an orthonormal basis of $\gt$, 
and let $e_\alpha$ be a basis of the $\alpha$-th root space
such that $e_\alpha$ and $e_{-\alpha}$ are dual to each other under 
the Killing form. Let $\Phi$ denote the set of roots.

Given $\lambda\in \gt$, we can consider the bundle 
$P_\lambda$ with automorphy factors $1$, $exp(2\pi i \lambda)$. 
It is a $G$-bundle, with a reduction of its structure group to 
$T$. Any $W$-orbit results in an isomorphic $G$-bundle. 
This construction exhibits $\gt$ as a branched cover of the moduli space
$Bun_\Sigma(G)$ 
of S-equivalence classes of semi-stable $G$-bundles. 
\begin{equation}
\label{eq-branched-covering-of-moduli-of-bundles}
\gt \ \ \ \longrightarrow \ \ \ Bun_\Sigma(G).
\end{equation}
The Galois group is the affine Weyl group, 
an extension of the Weyl group by the tensor product
$Ch(T)^*\otimes_\bbz H_1(\Sigma,\bbz)$ of the co-character lattice of $T$
with the period lattice of the elliptic curve. 
We consider $Ch(T)^*\otimes_\bbz H_1(\Sigma,\bbz)$ as a full lattice in
the complex Lie algebra $\gt$. 
The quotient of $\gt$ by this lattice is the group 
$Ch(T)^*\otimes_\bbz \Sigma$,
which is a cartesian product of $r$ copies of the elliptic curve.
The branched covering map 
(\ref{eq-branched-covering-of-moduli-of-bundles}) 
is the composition of the covering map 
$\gt\rightarrow Ch(T)^*\otimes_\bbz \Sigma$ and 
the quotient of $Ch(T)^*\otimes_\bbz \Sigma$
by the classical Weyl group $W$.
The ramification locus of
(\ref{eq-branched-covering-of-moduli-of-bundles})
consists of the translates by $Ch(T)^*\otimes_\bbz H_1(\Sigma,\bbz)$ 
of all the walls in $\gt$. 
The complement of the branch locus in $Bun_\Sigma(G)$, 
consists of the regular and semi-simple bundles \cite{FM}. The
automorphism
group of these bundles is isomorphic to the Cartan. 

We will assume that $\lambda$ does not belong to the branch
locus. The adjoint Lie algebra bundle 
$P_\lambda(\gg)$ splits into a trivial summand $P_\lambda(\gt)$ and a sum 
of non-trivial line bundles $L_{<\alpha,\lambda>}$, corresponding to the 
$\alpha$-th root space. 
The splittings (\ref{eq-splitting-two-summands}) and 
(\ref{eq-three-summands}), 
of the space of sections of a line bundle over $U_+\cap U_-$, 
extend to an analogous splitting for $ad (P_\lambda )$. 
\begin{equation}
\label{eq-splitting-of-adjoint-bundle}
P_\lambda(\gg)(U_+) \ \ \oplus \ \ [\gt\otimes {\rm span}_\bbc\{\rho(z)\}] 
\ \ \oplus \ \ \widetilde{P_\lambda(\gg)(U_-)}.
\end{equation}
If we set
\[
r(\lambda, z)= \rho (z)\sum_{i=1}^rx_i\otimes x_i +
\sum_{\alpha\in \Delta} \sigma_{-<\alpha,\lambda>}(z) e_\alpha\otimes 
e_{-\alpha},
\]
then $r(\lambda,z-z')$ is the kernel of the projection operator $P_+$ onto
$P_\lambda(\gg)(U_+)$. This function is precisely the $r$-matrix of
Felder; see \cite{etingof-varchenko}, formula (4.7). 
It is an elliptic solution of the CDYBE (\ref{eq-cdybe}).

Let $P_-$  and $P_0$ denote the projections to the corresponding 
summands of the splitting (\ref{eq-splitting-of-adjoint-bundle}).
Set 
\begin{equation}
\label{eq-R-operator}
R_\lambda \ \ :=  \ \ (P_+ - P_0- P_-)/2 \ \ = \ \  P_+  - \bbi/2.
\end{equation}
We think of $R_\lambda$ as an endomorphism of the loop algebra 
$L\gg$, which is determined by a $G$-bundle $P_\lambda$ together with a 
choice 
of an automophic factor (and hence a trivialization on $U_+$).
We will denote $R_\lambda$ by $R$, when it does not lead to any ambiguity.
We note that $R_\lambda$ is skew adjoint, with respect to the pairing 
defined using the Killing form and integration around the origin.



\subsubsection {The Sklyanin Poisson structure and its reduction}

Let $V$ be the complement in $\gt^*$, of the ramification locus 
of (\ref{eq-branched-covering-of-moduli-of-bundles}), where we identify
$\gt$ with $\gt^*$ via the Killing form. 
The bundles parametrized by $V$ are non-special, 
in the sense that none of the 
$L_{<\alpha,\lambda>}$, in the decomposition of the adjoint bundle, 
is a trivial line-bundle. 
The operators $R_\lambda$, given in (\ref{eq-R-operator}),  
are used by Etingof and Varchenko \cite{etingof-varchenko} 
to define a Poisson structure on $V\times LG\times V$. 
Let $f$ and $g$ be functions on $LG$. Given elements $a$ and $b$ of $\gt$,
we denote by $a_1$ and $b_1$ the corresponding linear function on
the first factor $V$, while $a_2$ and $b_2$ are the linear functions
on the second factor $V$. 
The Poisson structure on $V\times LG\times V$ is defined 
at $(v_1, \varphi, v_2)$ by the following relations. 

\begin{eqnarray} 
\nonumber
\{a_i, b_j \} &=& 0 \ \ \ \ i,j = 1,2,
\\
\nonumber
\{a_1, f\} &=& X^R_a(f), 
\\
\nonumber
\{a_2, f\} &=&X^L_a(f),
\\
\label{eq-Sklyanin-Poisson-structure}
\{f,g \} &=& <R_{v_1}(Df), Dg>- <R_{v_2}(D'f),D'g>.
\end{eqnarray}

\noindent
Here $X^L_a, X^R_a$ denote the left and right-invariant vector fields 
on $LG$ corresponding to $a$. 
We denote by $Df$ the left differential in $L\gg$ of 
$f$, and $D'f$ the right differential:
in terms of the Maurer-Cartan forms $\theta = d\varphi\cdot \varphi^{-1}$ 
and 
$\theta' =  \varphi^{-1}d \varphi$, thought of as maps $TLG\rightarrow 
L\gg$, 
we have 
\begin{eqnarray*}
Df &=& (\theta^{-1})^* (df),\cr
D'f &=& ({\theta'}^{-1})^* (df).
\end{eqnarray*}
In matricial terms,
\[
< Df, \dot  \varphi  \varphi^{-1}> \ \ \ = \ \ \ 
<df,\dot  \varphi> \ \ \ = \ \ \  <D'f,\varphi^{-1}\dot \varphi>,
\]
so that we get the identifications
\[
Df =  L_\varphi df \ \ \ \mbox{and} \ \ \  
D'f = R_\varphi df.
\]
At an element $ \varphi$ of $LG$, 
\[
Df = Ad_\varphi(D'f).
\]

There  are left and right actions of $T$ on $V\times LG\times V$, 
given by 
\[
L_t(v_1, g, v_2) = 
(v_1, tg, v_2),\quad R_t(v_1, g,v_2) = (v_1, gt^{-1}, v_2),
\]
with moment maps 
\[
(v_1, g, v_2)\mapsto v_1, \quad (v_1, g, v_2)\mapsto -v_2
\]
respectively. The reduction by the conjugation action of $T$ is then  the 
space 
\[
\Delta(V)\times LG/\simeq,
\]
where $\Delta(V)\subset V\times V$ is the diagonal, 
and $\simeq$ is the equivalence relation 
$(v, g, v)\simeq (v,tgt^{-1},v)$, $t\in T$.

\subsubsection{Relation to the algebro-geometric structure}

As the 
bundles parametrised by $V$ all have $T$ as their global group of 
automorphisms, our reduced phase space can be identified as a Galois 
covering of an open set 
of the space of isomorphism classes 
\[
{\cal M} = \{ (P, \varphi)| P{\rm\ an\ } T-\mbox{bundle over} \ \Sigma, 
\ \ \ \varphi\in \Gamma(U_+\cap U_-, P(G)) \} \ \ / \ \  
\sim. 
\]
The Galois group is again the affine Weyl group $\widetilde{W}$.

Now choose a faithful representation $\rho$ of $G$ and an effective
divisor $D$ on $\Sigma$. 
A Zariski open subset ${\cal M}_D$, of the moduli space $M(G,c,\rho,D)$
from section \ref{sec-an-anti-symmetric-tensor-D}, embeds in 
${\cal M}$. ${\cal M}_D$ is the subset of ${\cal M}$ of simple pairs
$(P,\varphi)$, where $\varphi$ is meromorphic with poles only at $D$. It 
fibers over $V/\widetilde{W}$, and we take the fiber product
\begin{equation}
\widetilde {\cal M}_D = {\cal M}_D\times_{V/\widetilde{W}} V
\end {equation}

The tangent space of ${\cal M}_D$ is given by the first hyper-cohomology 
of the complex (\ref{eq-tangent-complex-of-M-D}). Its cotangent bundle is
given in terms of the complex (\ref{eq-cotangent-complex-M-D}). 
The algebro-geometric Poisson structure $\Psi$ is given by 
(\ref{eq-Psi-M-D}).
We have a diagram 
\[
\begin{array}{ccc}
ad(P,\varphi,\rho)^*(-D) & {\buildrel{-ad_\varphi}\over {\longrightarrow}} 
& P\gg
\\
-L_\varphi \downarrow \hspace{3ex} & & \hspace{4ex}\downarrow L_\varphi
\\
P\gg & {\buildrel{ad_\varphi}\over {\longrightarrow}} &
ad(P,\varphi, \rho)(D).
\end{array}
\]
The top and left homomorphisms are defined in terms of the 
generic isomorphism 
(\ref{eq-the-dual-of-ad-phi-rho-is-the-one-for-phi-inverse}). 
For bundles of degree zero, the tangent space fits into an exact sequence
\[ 
H^0(\Sigma, ad(P,\varphi,\rho)(D))/ 
ad_\varphi(H^0(\Sigma,P\gg)) \ \ \rightarrow \ \ 
T{\cal M}_D \ \ \rightarrow \ \ 
H^1(\Sigma, P\gg).
\]
The cotangent space, dually, fits into
\begin{equation}
\label{eq-exact-sequence-of-hypercohomology-of-cotangent-complex}
H^0(\Sigma, P\gg) \ \ \rightarrow \ \ T^*{\cal M}_D
\rightarrow \ \ 
[ker: H^1(\Sigma, ad(P,\varphi,\rho)^*(-D))\rightarrow
H^1(\Sigma, P\gg)].
\end{equation}
$H^0(\Sigma, P\gg)$ is isomorphic to $\gh$ and $H^1(\Sigma, P\gg)$
to $\gh^*$.
Let us represent the hypercohomology on the level of cocycles. To do this 
we use our covering $U_+, U_-$, and assume for simplicity that the 
support of the divisor $D$ does not intersect $U_+$. A cocycle 
representing an element of the cotangent space is given by 
\begin {equation}(a_\pm, b_+, b_-),\end {equation}
where $a_\pm$  is a section of $ad(P,\varphi,\rho)(-D)$ over $U_+\cap
U_-$,
$b_+, b_-$ are sections of $P\gg$ over  $U_+, U_-$ respectively, 
with $-ad_\varphi(a_\pm) = b_+-b_-$ on the overlap. 
Similarly, represent an element of the tangent space by a cocycle 
\[
(\alpha_\pm, \beta_+, \beta_-),
\]
where $\alpha_\pm$  is a section of $P\gg$ over $U_+\cap U_-$,
$\beta_+, \beta_-$ are sections of $ad(P,\varphi,\rho)(D)$ over  $U_+,
U_-$
respectively, 
with $ad_\varphi(\alpha_\pm) =\beta_+-\beta_-$ on the overlap. 
The pairing on the level of  cocycles is given by 
\[
(\alpha_\pm, \beta_+, \beta_-), (a_\pm, b_+, b_-) =
<\alpha_\pm, 
(b_++b_-)>
+<a_\pm, (\beta_+ +\beta_-)>.
\]
(This formula is more or less forced by the constraint that the pairing 
must give zero on coboundaries).
Now let $(a_\pm, b_+, b_-)$
represent  the differential $df$ of a function on ${\cal M}_D$ at the 
point $(P,\varphi)$.
Apply the projections  to $L_\varphi(a_\pm), R_\varphi(a_\pm)$ to get 
\begin{eqnarray*}
L_\varphi(a_\pm) &=& \mu_+ +\mu_0 +\mu_-,\\
R_\varphi(a_\pm) &=&
\nu_+ +\nu_0 +\nu_-.
\end{eqnarray*}
Recall that the differential $ad_\varphi$ of the tangent complex was 
defined as $L_\varphi-R_\varphi$. The differential of the cotangent
complex
can be interpreted as $R_\varphi-L_\varphi$. 
Hence, 
\[
R_\varphi(a_\pm)-L_\varphi(a_\pm)=b_+-b_-.
\]
  From the exactness of 
(\ref{eq-exact-sequence-of-hypercohomology-of-cotangent-complex}), 
one has $\mu_0 =\nu_0$, and
\[
b_+ = -\mu_+ +\nu_+ + b,\quad b_- = \mu_--\nu_- + b,
\]
where $b$ is a constant, global section of $P\gg$, lying in $P(\gt)$. 

We need to determine the form of these different cocycles, for the 
differentials of the various 
functions considered in the definition 
(\ref{eq-Sklyanin-Poisson-structure}) of the Poisson bracket
on our phase space. 
One is fairly straightforward: with respect to the 
fibration of ${\cal M}_D$ over $V$, the differentials of 
functions lifted from $V$ are of the form $(0,b,b)$. 

One also has an inclusion $\widetilde{\cal M}_D\rightarrow V\times 
LG/\simeq$, given
as a 
subset of pairs $(T, \varphi)$ with $T= (T_1, T_2)$ representing the
automorphy factors 
of the bundle $P$, $\varphi $ a meromorphic function $\bbc\rightarrow G$,
holomorphic on $U_+$,  with  poles only at 
$D$ and 
$Ad_{T_i}(\varphi (z)) = \varphi(z+w_i)$ for the periods $w_i$ of the 
curve. We assume for simplicity  that $D$ does not intersect $U_+$
On the level of tangent spaces, the inclusion is given by 
\[
(\alpha_\pm, \beta_+, \beta_-) \mapsto  (\alpha_\pm,
\beta_+),
\]
with $ad_\varphi(\alpha_\pm) = \beta_+ - \beta_-$. 
Dually, one has 
\[ 
(a_\pm, b) \mapsto (a_\pm, b_+ +b, b_- +b),
\]
with $-ad_\varphi(a_\pm) = b_+-b_-$. The elements  $b_+,b_-$ are
determined
by comparing the pairing on bigger space $V\times LG/\simeq$:
\[ 
2<\beta_\pm, b> + 2<a_\pm, \alpha_+>
\]
and that on the smaller space ${\cal M}_D$:
\[
\begin{array}{l}
<\beta_\pm, b_+ +b_- +2b> + <a_\pm,\alpha_++\alpha_->
\\
= \ \ <\beta_\pm, b_+ +b_- +2b> + <a_\pm, 2\alpha_+- ad_\varphi(\beta_\pm)>
\\
= \ \  2<\beta_\pm, b> + 2<a_\pm, \alpha_+> + 
<\beta_\pm, b_+ +b_- + ad_\varphi(a_\pm)>
\\
= \ \  2<\beta_\pm, b> + 2<a_\pm, \alpha_+> + <\beta_\pm, 2b_->,
\end{array}
\]
 whence $b_-$ pairs to zero with any cocycle in $\gt\cdot \rho(z)$,
 forcing the equalities 
$b_- = P_-(ad_\varphi(a_\pm))$ and $b_+ = -P_+ (ad_\varphi(a_\pm))$. 
 In summary, the differential a function on $LG/\simeq$ is represented by 
a cocycle $a_\pm$, such that $ad_\varphi(a_\pm)$ is a coboundary. Mapping 
to the forms on ${\cal M}_D$, the cocycle lifts 
to 
$(a_\pm, -P_+ (ad_\varphi(a_\pm)),P_-(ad_\varphi(a_\pm)))$, giving 
\[
b_+ = -\mu_+ +\nu_+, \quad b_- = \mu_--\nu_- .
\]

Choosing another function $g$ on ${\cal M}_D$, one can represent the
differential $dg$ as 
\[
dg\simeq (c_\pm, d_+, d_-)
\]
with 
\begin{eqnarray*}
L_\varphi(c_\pm) &=& \xi_+ +\xi_0 +\xi_-,\cr R_\varphi(c_\pm)&=&
\zeta_+ +\zeta_0 +\zeta_-,\cr \xi_0 &=& \zeta_0
\end{eqnarray*}
giving 
\[
d_+ = -\xi_+ +\zeta_+ +d,\quad  d_- = \xi_--\zeta_- +d.
\]

We now want to compute the Poisson brackets of $f$ and $g$. 
The Poisson tensor $\Psi$ acts by $-L_\varphi$ on 1-cocycles, and by 
$L_\varphi$ on 0-cocycles. 
We modify the cocycle $\Psi(c_\pm, d_+, d_-)= (-\xi_+ -\xi_0 -\xi_-,
 L_\varphi(-\xi_++\zeta_+ + d),  L_\varphi(\xi_--\zeta_- + d) $ by the 
coboundary 
$(\xi_++\xi_-, ad_\varphi(\xi_+), ad_\varphi(-\xi_-))$,
so that the new representative cocycle is 
\[
\Psi(dg)\simeq (-\xi_0, -R_\varphi(\xi_+)+
L_\varphi(\zeta_+) + L_\varphi(d), 
R_\varphi(\xi_-) - L_\varphi(\zeta_-)+ L_\varphi(d))
\]
The algebro-geometric Poisson bracket (\ref{eq-Psi-M-D}) is given by
\[
\begin{array}{l}
-2\{f,g\} \ = \   <df,\Psi(dg)>
\\
= \  <a_\pm, R_\varphi(-\xi_++\xi_-)+
L_\varphi(\zeta_+-\zeta_-)+ 2L_\varphi(d)> - <\xi_0, b_+ +b_->
\\
= \  <L_\varphi(a_\pm), -\xi_++\xi_-> + 
<R_\varphi(a_\pm), \zeta_+-\zeta_- + 2d> - <\xi_0, b_+ + b_->
\\
= \ <\mu_+ +\mu_0 +\mu_-, -\xi_++\xi_-> + 
<\nu_+ +\nu_0 +\nu_-, \zeta_+-\zeta_- + 2d>
\\ 
\ \ \  - <\xi_0,-\mu_+ +\nu_+ + \mu_--\nu_- +2b>
\\
= \ <\mu_+ +\mu_0 +\mu_-, -\xi_++\xi_-> +
 <\nu_+ +\nu_0 +\nu_-, \zeta_+-\zeta_-> + 
\\ 
\ \ \ <\mu_+  + \mu_0  +\mu_-, \xi_0> - <\nu_+ +\nu_0 + \nu_-,\zeta_0>
+ <\mu_0, 2d> -<\xi_0, 2b>
\\
= \  <\mu_+ +\mu_0 +\mu_-, -\xi_++\xi_0 +\xi_-> +
 <\nu_+ +\nu_0 +\nu_-, \zeta_+ - \zeta_0 -\zeta_->
\\ 
\ \ \ + <\mu_0, 2d> -<\xi_0, 2b> 
\end{array}
\]
where we use the fact that the $+$-summand in our decomposition is 
isotropic, the $-,0$-summands are orthogonal and isotropic, and the fact 
that global sections are orthogonal to the $+,-$ summands.

Now suppose that $f,g$ are lifted from $LG/\simeq$; we have then that the 
terms $b, d$ vanish, and $L_\varphi(a_\pm) = Df$, 
$R_\varphi(a_\pm)= D'f$, while 
$\xi_+ = P_+(Dg)$, etc., giving for the Poisson bracket:
\[
\{f,g\} =<Df, R(Dg)>-<D'f, R(D'g)>, 
\]
where we recall that $R = (P_+ -P_0-P_-)/2$. 

 Next, let us consider the case when one of the functions, $f$, is lifted 
from  $LG/\simeq$, while the other, $g$, is lifted from $V$, and so is 
represented by a cocycle $(0,d,d), d\in \gt$. We then get
\[
\{f,g\} =<Df,d>
\]
This is minus the  derivative of $f$ along the 
left invariant vector field corresponding to $d$; one could take the 
right derivative also, as the quotienting out by the adjoint action means 
that the two are equivalent. Finally, when both functions are lifted from 
$V$, the formulae above show that their bracket is zero.  We have thus 
shown that our algebro-geometric
structure
satisfies the defining relations (\ref {eq-Sklyanin-Poisson-structure})
of the Sklyanin bracket. We note, that the algebro-geometric structure 
is defined on the whole moduli space $M(G,c)$. 
This moduli space includes pairs $(P,\varphi)$,
where the $G$-bundle $P$ belongs to the ramification locus of 
(\ref{eq-branched-covering-of-moduli-of-bundles}), in particular those 
which are not reducible to the torus. Summarising, we have

\begin{thm}
\label{thm-comparison-with-cdybe-simply-connected-case}
The Sklyanin Poisson structure (\ref{eq-Sklyanin-Poisson-structure})
and the algebro-geometric structure (\ref{eq-Psi-M-D}) coincide.
The Sklyanin bracket is thus invariant under the affine Weyl group,
and it extends across the walls on the quotient.  
\end{thm}


\subsection {Arbitrary connected reductive $G$}

We would like to consider briefly how these results extend to the case of 
arbitrary reductive  groups. A first remark is that they extend 
unchanged to the case of $G$ a product of a simply connected 
reductive group and a torus. Secondly, we can take an arbitrary connected 
reductive group $\bar{G}$ and fit it into an exact sequence
\[
1\rightarrow K \rightarrow  G \rightarrow \bar{G}
\rightarrow 1,
\]
with $  G$ a product $G_0\times T$ of a simply connected reductive 
group and a 
torus, and $K$ central, finite. For example, one can have 
$\bar{G} = Gl(N,\bbc),   G = Sl(N,\bbc)\times \bbc^*, Z = \bbz/N$. 
Corresponding to this sequence, we have an exact 
sequence of (non-abelian) cohomology
\begin{equation}
\label{G-sequence}
\rightarrow H^1(\Sigma, {\cal K}) = K^2\rightarrow
H^1(\Sigma,  {\cal G})
\rightarrow H^1(\Sigma, \bar{\cal G})\rightarrow 
H^2(\Sigma, {\cal K})=K\rightarrow, 
\end{equation}
where ${\cal K}, {\cal G},\bar{\cal G}$ are the sheaves of holomorphic 
maps into $K, G, \bar{G}$ and the exactness of the sequence is as a 
sequence of pointed sets. What we see is that there is an obstruction, 
essentially a first Chern class, to lifting a $\bar{G}$ bundle to a $  G$
bundle. Also, since it is a Chern class, the obstruction vanishes when 
one lifts to an $N$-fold cover $\hat\Sigma\rightarrow \Sigma$ by another 
elliptic curve, for some suitably chosen $N$ (We will choose the cover to 
be cyclic, of the form $\bbc/(N\bbz + \tau \bbz)\rightarrow 
\bbc/(\bbz + \tau \bbz).$

Following \cite {FM,BFM,atiyah-bott}, the generic $\bar{G}$-bundle 
can be described 
in terms of automorphy factors which are almost commuting pairs in 
$ G$, that is, their commutators are central. Let $g_a$, $g_b$ be 
the group elements of  $ G$ corresponding to the standard 
generators $a,b$ of the fundamental group of the curve; 
by Corollary (4.2) of \cite {BFM}, 
they can be chosen generically so that $g_b$ lies in the 
torus $ T$  of $ G$, $g_a$ lies in the normaliser $N(T)$; 
furthermore, as the bundle lifted to the covering curve $\hat \Sigma$ 
reduces to $G$, one can suppose that $(g_a)^N$ lies in the Abelian 
factor $T$, so that on the covering curve one has a reduction to the 
torus. 

These bundles on the covering curve are invariant under deck 
transformations. Over $\hat \Sigma$, the element $a$ represents the deck 
transformation for the cover, 
and $g_a$ its action on the associated bundle. The pull-back of the 
$\bar{G}$ bundle to $\hat \Sigma$ lifts to $G$, and reduces to a $\bar{T}$ 
bundle $P_{\bar{T}}$. The invariance under deck transformation gets
translated into an
 invariance 
\[
 Ad^{g_a}(P_{\bar{T}}) =  \rho^*( P_{\bar{T}})
\]
where $\rho:\hat \Sigma \rightarrow \hat \Sigma$ 
is the deck transformation, and on the lift to 
a $  T$-bundle $  P_T$, referring to (\ref{G-sequence}), \
\[
 Ad^{g_a}(P_T) =  \rho^*( P_T) \otimes P_Z
\]       
where $P_Z$ is some root of unity in the Jacobian corresponding to an 
element of $H^1(\hat \Sigma, {\cal Z})$.

We can then compare the 
algebro-geometric definitions of Poison structures for the 
$\bar{G}$-bundles on $\Sigma$ and $G$-bundles on $\hat \Sigma$, 
seeing in essence that one is the restriction of 
the other to the fixed point set of the action above. Similar 
considerations hold for the Sklyanin bracket. Indeed,
on the level of infinitesimal deformations, we then have for a 
$\bar{G}$-bundle $\bar{P}$ that $H^1(\Sigma, \bar{P}\gg)=\gt^{g_a}$. 
Under the identification of the loop algebra  with sections of 
$P\gg$ on the overlap $U_+\cap U_-$, we  have a 
splitting  of the loop algebra into 
\[
P\gg(U_+)\oplus \gt^{g_a}(\rho(z)) \oplus \widetilde 
{P\gg(U_-)}.
\]  

 Etingof and Schiffman 
\cite {etingof-schiffman} consider precisely the Lie-Poisson groupoids on 
such fixed point sets, and it is clear that their solutions to the CDYBE 
correspond to the splitting above.

In any case, the equivalence of the algebro-geometric form with the 
Sklyanin one provides an alternate proof of the Jacobi identity of  
the former, first, for the simply connected case, and, then, by 
restriction to a fixed point set, in the general case.

\section{Systems related to partial flag varieties}
\label{sec-products-of-partial-flag-varieties}

A particularly interesting list, among moduli spaces of type $M(G,c,\orbit)$,
was introduced in Example 
\ref{example-partial-flag-varieties-of-maximal-parabolics}.
These moduli spaces are related to products of flag varieties of 
maximal parabolic subgroups of a simple group $G$. We study in this section 
the abelian varieties, that occur as fibers of the characteristic 
polynomial map, and their spectral curves. We are particularly interested
in such examples of minimal dimension. So, we consider pairs $(P,\varphi)$, 
where the meromorphic section $\varphi$ has singularities only at two 
points $p_1$ and $p_2$ in the elliptic curve $\Sigma$.
We choose $\orbit_{p_1}$ and $\orbit_{p_2}$ to be $W$-orbits of the 
co-characters $a$ and $-a$, where $a$ is dual to a particular choice of 
a simple root.

According to the list in Example 
\ref{example-partial-flag-varieties-of-maximal-parabolics},
the group $G$ must be either a classical group, or exceptional of type $E_6$
or $E_7$, 
in order for it to have a maximal parabolic of the type we are interested.
These are the maximal parabolic subgroups, whose partial flag variety 
embeds as an orbit of the loop grassmannian. 
The type A examples are discussed in section
\ref{sec-type-A}. A deformation of the Elliptic Calogero-Moser system
arises. The three examples of type $D_n$ are discussed in section 
\ref{sec-type-D}.

\subsection{Type $A$ examples}
\label{sec-type-A}
Section \ref{sec-products-of-grassmannians} exhibits moduli spaces of type 
$M(PGL(n),c,\orbit)$, which are birational to 
$[\bbz/n\bbz\times \bbz/n\bbz]$-quotients of products of Grassmannians. 
Their spectral curves are described in Section \ref{sec-calogero-moser}.
Certain symplectic leaves of the generalized Hitchin system 
admit deformations to symplectic leaves $M(G,c,\orbit)$
of the generalized Sklyanin systems. We discuss such examples
in Section \ref{sec-calogero-moser}. In particular, Lemma 
\ref{lemma-example-deforms-to-the-calogero-moser-of-type-A-n} exhibits a
deformation of the type $A_n$ Calogero-Moser system via a one-parameter family
of Sklyanin systems. 

\subsubsection{Products of Grassmannians}
\label{sec-products-of-grassmannians}

When the group $G$ is $PGL(n)$, certain moduli spaces $M(PGL(n),c,\orbit)$ 
are birational to $[\bbz/n\bbz\times \bbz/n\bbz]$-quotients of 
products of grassmannians (Lemma 
\ref{lemma-symplectic-leaf-is-birational-to-product-of-flag-varieties}). 
We choose $c$ to be a generator of $\pi_1(PGL(n))$ and denote by 
$P$ the unique stable bundle. 
Consequently, 
a Zariski open subset of $M(PGL(n),c,\orbit)$ is the $\Aut(P)$-quotient, of 
the set of meromorphic sections of $P(G)$ with singularities in $\orbit$. 
Of particular interest is the case $G(k,n)\times G(n-k,n)$, 
corresponding to partial flag varieties associated to 
two opposite maximal parabolic subgroups. The Zariski open subset of 
$G(k,n)\times G(n-k,n)$, of transversal pairs, 
is naturally isomorphic to the twisted cotangent bundle of $G(k,n)$. 
The twisted cotangent bundle of $G(k,n)$ is isomorphic to a 
semi-simple coadjoint orbit $\CoadjointOrbit_{ss}$ of $\gsl_n$
(as in Theorem \ref{thm-generalized-hitchin}, 
we reserve the notation $\CoadjointOrbit$ for coadjoint orbits). 
The orbit $\CoadjointOrbit_{ss}$ is the coadjoint orbit of idempotents
of rank $k$, translated by a multiple of the identity endomorphism to 
annihilate the trace. 
Let $\CoadjointOrbit:=\sum_{p\in \Sigma}\CoadjointOrbit_p$
such that $\CoadjointOrbit_p=0$ if $p\neq p_0$ and
$\CoadjointOrbit_{p_0}=\CoadjointOrbit_{ss}$. Connected components of 
the Hitchin system $\Higgs_\Sigma(PGL(n),c,\CoadjointOrbit)$ 
are birational to a $[\bbz/n\bbz\times \bbz/n\bbz]$-quotient of the 
orbit $\CoadjointOrbit_{ss}$. 
Lemma \ref{lemma-symplectic-leaf-is-birational-to-product-of-flag-varieties} 
implies that, for a suitable $\orbit$, the moduli 
$M(PGL(n),c,\orbit)$ is birational to 
$\Higgs_\Sigma(PGL(n),c,\CoadjointOrbit)$. 
The data $\orbit$ depends on a point on the elliptic curve (in addition to 
the origin), while the data $\CoadjointOrbit$ does not.
We get a deformation of $M(PGL(n),c,\orbit)$ 
which is supported, birationally, on a fixed space.
The Hitchin system $\Higgs_\Sigma(PGL(n),c,\CoadjointOrbit)$
is the limit. We will carry out the limit calculation in a special case
(Lemma \ref{lemma-example-deforms-to-the-calogero-moser-of-type-A-n}).
The general case is similar.

\begin{example}
\label{example-PGL-n-moduli-is-disconnected}
{\rm
The group $\Aut(P)$, of global automorphisms of a stable rigid 
$PGL(n)$-bundle $P$, is isomorphic to $\bbz/n\bbz\times \bbz/n\bbz$. 
%
More canonically, we have an isomorphism
\begin{eqnarray*}
\Aut(P) & \longrightarrow & \Hom[\pi_1(\Sigma),\pi_1(PGL(n))]
\\
\varphi & \mapsto & \tau_\varphi,
\end{eqnarray*}
sending a  global holomorphic section of $P(PGL(n))$, to
its topological invariant $\tau_\varphi$ given in 
(\ref{eq-topological-invariant}). 
This isomorphism is the connecting homomorphism of the long exact sequence
(\ref{eq-long-exact-seq-of-analytic-sheaf-coho}) 
in cohomology. 
The stable $PGL(n)$-bundle $P$ of type $c$ admits a lift to 
a stable vector bundle $E$ of degree $d$, where $d\equiv c (\Mod \ n)$. 
The associated bundle $P(SL(n))$ in 
(\ref{eq-long-exact-seq-of-analytic-sheaf-coho}) is precisely the bundle 
$\Aut_1(E)$, of automorphisms of $E$ with determinant $1$. 
Since $E$ is simple, only the identity automorphism of $P$ lifts to 
a global section of $P(SL(n))$. The injectivity of the 
connecting homomorphism follows. Surjectivity follows, since the homomorphism 
$H^1(\Sigma,\pi_1(PGL(n)))\rightarrow H^1(\Sigma,P(SL(n)))$ in 
(\ref{eq-long-exact-seq-of-analytic-sheaf-coho}) is the constant
map; tensorization of the stable vector bundle $E$, by a line bundle of order 
$n$, is isomorphic to $E$. 

}
\end{example}

We introduce next the moduli spaces, which are related to products of 
Grassmannians. 

\begin{example}
\label{example-product-of-grassmannians}
{\rm
If $G=PGL(n)$ and  the simple root $\alpha$ 
in Example \ref{example-partial-flag-varieties-of-maximal-parabolics} 
is $\epsilon_{k}-\epsilon_{k+1}$, then $G/P_\alpha$ is 
the grassmannian $G(k,n)$. 
The orbit $\orbit_{\alpha^*}$ is represented by 
the cocharacter ${\rm diag}(t,\dots,t,1,\dots,1)$
of $GL(n)$, with simple zeros at the first $k$ diagonal entries. 
The orbit $\orbit_{-\alpha^*}$ 
is represented by the cocharacter ${\rm diag}(1,\dots,1,t,\dots,t)$,
with simple zeros at the last $n-k$ entries.

Let $c\in \pi_1(PGL(n))\cong (\bbz/n\bbz)$ be an additive generator. 
Fix a topological type $\tau:\pi_1(\Sigma^0)\rightarrow \pi_1(PGL(n))$ 
as in (\ref{eq-topological-invariant}).
We have a Zariski open subset $M^0_\tau$ in $M(G,c,\orbit)$,
consisting of pairs $(P,\varphi)$, of topological type $\tau$, 
where $P$ is the rigid stable bundle. 
Fix two points $p_+$, $p_-$ in $\Sigma$ and let 
$\orbit=\sum_{p\in\Sigma} \orbit_p$, with $\orbit_p=PGL[[t]]$ if 
$p\not\in \{p_+,p_-\}$, $\orbit_{p_+}=\orbit_{\alpha^*}$ and
$\orbit_{p_-}=\orbit_{-\alpha^*}$. 
Formula 
(\ref{eq-sum-over-roots-in-opposite-unipotent-radical})
implies, that the dimension of $M(G,c,\orbit)$ is $2k(n-k)$. 
Choose a lift $E$ of the stable $PGL(n)$-bundle $P$ to a stable vector 
bundle. Let 
\[
ev \ : \ M^0_\tau \ \ \rightarrow \ \ [G(k,E_{p_+})\times G(n-k,E_{p_-})]/\Aut(P)
\]
be the morphism sending the meromorphic section $\varphi$ of 
$P(G)$ to its pair of kernels in the fibers over $p_+$ and $p_-$. 
The pair of kernels is well defined only up to the diagonal 
action of $\Aut(P)$. 
Note that $M^0_\tau$ admits an \'{e}tale symplectic morphism onto 
a Zariski open subset of $\Hecke(G,c,c,\orbit)$ 
(Section \ref{sec-dimensions-of-Hecke-and-M} and Theorem 
\ref{thm-Poisson-structure-on-Hecke}).
The image of $ev$ is isomorphic to a Zariski open subset 
of $\Hecke(G,c,c,\orbit)$ 
(see Claim \ref{claim-Hecke-with-fixed-P-2} for the relationship between 
$\Hecke(G,c,c,\orbit)$ and the product of Grassmannians). 
}
\end{example}

More generally, let $S\subset \Sigma$ be a finite subset and
$\orbit=\sum_{p\in \Sigma}\orbit_p$, 
such that $\orbit_p=PGL(n)[[t]]$, if $p\not\in S$, and 
$\orbit_p=\orbit_{\alpha_{k(p)}^*}$ for some simple root 
$\alpha_{k(p)}=\epsilon_{k(p)}-\epsilon_{k(p)+1}$. 
Assume that 
\[
\sum_{p\in S} k(p) \equiv 0 \ \  (\Mod \ n). 
\]
(Compare with Remark \ref{rem-topological-non-emptiness-condition}). 
Let 
\[
ev \ : \ M^0_\tau \ \ \rightarrow \ \ 
\left[\prod_{p\in S} G(k(p),E_p)\right]/\Aut(P)
\]
be the morphism sending the meromorphic section $\varphi$ of 
$P(G)$ to its kernels in the fibers over $S$.

\begin{lemma}
\label{lemma-symplectic-leaf-is-birational-to-product-of-flag-varieties}
The morphism $ev$ is injective and dominant. 
\end{lemma}

\noindent
{\bf Proof:} 
Fix a stable $PGL(n)$-bundle $P$ of topological type $c$. 
The group $\Aut(P)$, of global sections of $P(G)$, is 
isomorphic to $\Hom[\pi_1(\Sigma),\pi_1(PGL(n))]$
(see Example \ref{example-PGL-n-moduli-is-disconnected}). 
Let $q:\widetilde{M}^0\rightarrow M^0$ be the 
$\Aut(P)$-Galois cover obtained by taking the global sections of $P(G)$
representing isomorphism classes in $M^0$. 
We get the commutative diagram:
\[
\begin{array}{ccc}
\widetilde{M}^0 & \LongRightArrowOf{\widetilde{ev}} & 
\prod_{p\in S} G(k(p),E_p)
\\
q \ \downarrow \ \hspace{1ex} & & \downarrow
\\
M^0 & \LongRightArrowOf{ev} & 
\left[\prod_{p\in S} G(k(p),E_p)\right]/\Aut(P).
\end{array}
\] 
The morphism $\widetilde{ev}$ is $[\Aut(P)\times\Aut(P)]$-equivariant. A pair 
$(f,g)\in [\Aut(P)\times\Aut(P)]$ takes $\varphi$ to 
$f\circ\varphi\circ g^{-1}$. The action on the product of Grassmannian 
factors through the projection $(f,g)\mapsto f$. 
We will see, that the group $\Aut(P)\times\Aut(P)$ acts 
transitively on the set of connected components of $\widetilde{M}^0$,
by changing the topological type $\tau$. 
Note, that each topological type is invariant under
the diagonal action of $\Aut(P)$ (via conjugation).

The injectivity of $\widetilde{ev}$, on the component determined by 
a fixed topological type $\tau$, is straightforward. If $\varphi$ and $\eta$ 
are two sections of in $M^0_\tau$ and 
$\widetilde{ev}(\varphi)=\widetilde{ev}(\eta)$,
then $\eta\circ \varphi^{-1}$ is a holomorphic section of $P(PGL(n))$ 
of trivial topological type. Hence, $\eta\circ \varphi^{-1}=id$.
The fact, that $\eta\circ \varphi^{-1}$ is holomorphic, 
can be seen more easily,
once we lift $\eta$ and $\varphi$ to meromorphic automorphisms 
$\tilde{\eta}$ and $\tilde{\varphi}$ of the 
pull back $f^*P(SL(n))$ over the branched covering 
$f:D_\tau\rightarrow \Sigma$ determined by the topological invariant
(see section \ref{sec-spectral-curves-and-group-isogenies}). 
Then $\tilde{\eta}$ and $\tilde{\varphi}$ are meromorphic automorphisms 
of the pull back $f^*E$ of a stable vector bundle $E$. Moreover, 
$\tilde{\eta}^*$ and $\tilde{\varphi}^*$ have the same image subsheaf
in $f^*E(S)$. Thus, $\tilde{\eta}\circ \tilde{\varphi}^{-1}$
is a holomorphic automorphism of $f^*E$.

Given a point  $x:=(K_p)_{p\in S}$ in $\prod_{p\in S} G(k(p),E_p)$,
we denote by $Elm(x)$ the  subsheaf of $E(S)$ of meromorphic 
sections of $E$ with poles along $(K_p)_{p\in S}$. Applying 
Lemma \ref{lemma-generic-elementary-transformation-is-semistable} 
$dn$-times, we get that $Elm(x)$ is stable for a generic point $x$.
$Elm(x)$ is a vector bundle of rank $n$ and determinant line bundle 
$\det(E)(\sum_{p\in S}k(p)\cdot p)$. 
Thus, if $Elm(x)$ is stable, then
it must be isomorphic to $E\otimes R$, for every $n$-th root $R$ of
$\StructureSheaf{\Sigma}\left(\sum_{p\in S}k(p)\cdot p\right)$
(see \cite{atiyah}). 
The composition $E\otimes R\cong Elm(x) \hookrightarrow E(S)$
represents an element $\varphi$ in $\widetilde{M}^0$. 
Composing with an element of $\Aut(P)$, we can translate 
$\varphi$ to an element of $\widetilde{M}^0$ with 
the topological invariant $\tau$ chosen. 
We conclude that $\widetilde{ev}$ is dominant. 
\EndProof


\begin{lemma}
\label{lemma-generic-elementary-transformation-is-semistable}
Let $E$ be a rank $n$ semi-stable vector bundle over an elliptic curve 
$\Sigma$ and $p$ a point in $\Sigma$. Assume that $E$ is not isomorphic to
the tensor product of a stable vector bundle with a trivial vector bundle of 
rank $>1$. 
The space $\bbp(\restricted{E}{p})$, of lines in the fiber, parametrizes 
subsheaves $E'$ of $E(p)$, containing $E$, with quotient $E'/E$ of length 
$1$. A Zariski dense open subset of 
$\bbp(\restricted{E}{p})$ parametrizes {\em semistable} 
subsheaves $E'$. Moreover, the generic such $E'$ is not isomorphic to
the tensor product of a stable vector bundle with a trivial vector bundle of 
rank $>1$. 
\end{lemma}

\noindent
{\bf Proof:} 
Let $\ell$ be a line in the fiber of $E$ over $p$. Denote by
$E(\ell)$ the subsheaf of $E(p)$ of sections with poles  along $\ell$.
We prove first that $E(\ell)$ is semi-stable, for a generic $\ell$.
Assume that $F$ is a  subbundle of $E$, of rank $r<n$, such that $F_p$ 
contains $\ell$ and $F(\ell)$ destabilizes $E(\ell)$. 
We may assume that $\ell$ is not contained in 
the fibers of any sub-bundle of $E$ with slope equal to the slope of $E$ 
(here we use the assumption on $E$). 
We have the inequalities 
\[
\frac{\deg(F)}{r} \ <  \  \frac{\deg(E)}{n} \ \ \ \ \ \ \mbox{and} \ \ \ \ \ \ 
\frac{\deg(F)+1}{r} \ > \ \frac{\deg(E)+1}{n}. 
\]
Combining both inequalities, we get
\[
\deg\SheafHom(F,E) \ \ < \ \ n-r.
\]
Semi-stability of $E$ implies that $\Hom(E,F)=0$. Hence, 
$\dim \Hom(F,E)=\deg\SheafHom(F,E)< n-r$.
Now, the dimension
of the family $H^1(\SheafEnd(F))$, of infinitesimal deformations of $F$,
is equal to $\dim H^0(\SheafEnd(F))$. The latter space is contained in 
$\Hom(F,E)$. The dimension of the family $B$ of sub-sheaves
of $E$, which are deformations of $F$, is
$\dim\Hom(F,E)-\dim\End(F)+\dim H^1(\SheafEnd(F))=\dim\Hom(F,E)$. 
We conclude that the dimension of $B$ is $< n-r$. Hence, the
family of lines $\ell$, in the fibers of the deformations of $F$,
has dimension $< n-1$.
Consequently, $E(\ell)$ is semi-stable, for a generic $\ell$. 

Suppose $E(\ell)$ is  isomorphic to the direct sum of $k$ copies of
the stable bundle $V$. Then 
\[
\dim\Hom(V^*,E^*)=\deg\SheafHom(V^*,E^*)=\rank(V).
\]
Thus, the image of $H^0(\SheafHom(V^*,E^*))\otimes_{\bbc} V^*$ in $E^*$ 
is equal to the image of $E(\ell)^*$. 
It follows that every embedding 
$V^*\hookrightarrow E^*$ maps the fiber $V_p^*$ into the hyperplane in 
$E^*$ annihilating $\ell$. 
Hence, if there is a line $\ell$ in $E_p$, for which $E(\ell)$
is isomorphic to $V^{\oplus k}$, then $\ell$ is unique. 
Since $\det E(\ell)$ is $\det(E)(p)$, then there are only finitely 
many such $V$.
\EndProof

\subsubsection{Symplectic leaves of minimal dimension deforming 
the Calogero-Moser 
system}
\label{sec-calogero-moser}

We study in this section the spectral curves of the $PGL(n)$ 
systems of Section \ref{sec-products-of-grassmannians}. 
We work with $GL(n)$ systems corresponding to a lift. We 
lift the singularity data $\overline{\orbit}$ for $PGL(n)$ to a
singularity data $\orbit$ of $GL(n)$. This means that a co-character of
$\orbit_p$, for each point $p\in \Sigma$, maps to a co-character of
$\overline{\orbit}_p$. Note that the lift may introduce new singularities.
If $\orbit_p$ is a new singularity, then its co-character has image 
in the center of $GL(n)$. The $PGL(n)$-systems are obtained from the 
$GL(n)$-systems by symplectic reduction: 
One fixes the determinant line bundle of the stable bundle
and takes the quotient by the $\bbc^\times$-action. 

\begin{example}
\label{example-dimension-of-symplectic-leaf-in-GL-case-and-simple-zeros-poles}
{\rm
Let $G=GL(n)$ and assume that the divisor $D$ in Example 
\ref{example-dimension-of-leaf-in-GL-case}
is reduced and $\varphi$ satisfies:
\[
E(D)(-p) \ \subset \ \varphi(E) \ \subset \ E(D), 
\]
in a neighborhood of every point $p\in \Sigma$. 
Denote by $r_p$ the rank of 
$\varphi_p:E_p\rightarrow E(D)_p$. In terms of the weights of the standard 
representation, the assumption on $\varphi$ translates to
\[
\Ord_p(\epsilon_i,\orbit_p) \ \ \in \ \ \{-1, 0, 1\}.
\]
Moreover, if $p\in D$, then the multiplicity of $1$ in 
$\{\Ord_p(\epsilon_i,\orbit_p)\}_{i=1}^n$ is $r_p$ and of $0$ is $n-r_p$.
If $p\not\in D$, then the multiplicity of $0$ is $r_p$ and of $-1$
is $n-r_p$. 
Condition (\ref{eq-linear-equivalence-of-zero-of-det-to-poles-of-det})
asserts that
\begin{equation}
\label{eq-down-to-earth-linear-equivalence-of-zero-and-poles-of-det}
\sum_{p\in D}r_p\cdot p \ \ \ \sim \ \ \
\sum_{p\in \Sigma\setminus D}(n-r_p)\cdot p.
\end{equation}
If $p$ is a point in $D$, then $\Ord_p(\bullet,\orbit_p)$ is represented by
$-(\lambda_{r_p},\bullet)$, where $\lambda_i=\sum_{j=1}^i\epsilon_i$. 
If $p\in \Sigma\setminus D$ and $r_p<n$, then 
$\Ord_p(\bullet,\orbit_p)$ is represented by $-(\lambda_{n-r_p},\bullet)$. 
Formula (\ref{eq-dimension-in-GL-case}) becomes
\begin{equation}
\label{eq-dimension-of-symplectic-leaf-in-GL-case-and-simple-zeros-poles}
\dim M(GL(n),c,\orbit) \ = \ 2 + \sum_{p\in \Sigma} r_p(n-r_p).
\end{equation}
(Compare with Lemma 
\ref{lemma-symplectic-leaf-is-birational-to-product-of-flag-varieties}). 

Let us analyze the singularities of
the spectral curve $C$ of $\varphi$ in the total space of the 
line bundle $\StructureSheaf{\Sigma}(D)$. Denote the zero section by 
$\sigma_0$. Since $\varphi$ is generically invertible, the 
scheme-theoretic intersection $C\cap \sigma_0$ is the divisor
$\sum_{p\in \Sigma} (n-r_p)\cdot \sigma_0(p)$ on $\sigma_0$. 
It follows that the multiplicity of the point $\sigma(p)$, as 
a point on the curve $C$, is at most $n-r_p$. 
We claim that the multiplicity is precisely $n-r_p$. 
Indeed, if $s$ is a local section of $\StructureSheaf{\Sigma}(D)$
with a simple zero at $p$, then the germ $C+s$, of the spectral curve
of $\varphi+s\bbi$, is the translation of $C$ by $s$. 
Since the rank of $\varphi+s\bbi$ at $p$ is $r_p$, $C+s$ intersects $\sigma_0$
with multiplicity at least $(n-r_p)$. Hence, $C$ intersects $s$ with
multiplicity $\geq (n-r_p)$. 
}
\end{example}

We consider next the spectral curves of Example 
\ref{example-product-of-grassmannians} in the special case, where
the flag variety is $\bbp^{n-1}$. 
In Example \ref{example-GL-calogero-mozer-deformed} 
we describe the spectral curves. In Lemma 
\ref{lemma-example-deforms-to-the-calogero-moser-of-type-A-n}
we show that the Calogero-Moser system is a limit of this 
integrable system. 

\begin{example}
\label{example-GL-calogero-mozer-deformed}
{\rm
Let $G=GL(n)$. 
Consider the case where $\orbit_p=GL[[t]]$ for all
$p\not\in \{p_0,p_1,p_2\}$, the polar divisor is $D=p_0$, 
and the ranks $r_i$ at $p_i$ are $r_0=n$, $r_1=1$, and $r_2=n-1$. 
If $n=2$, we assume further that $p_1\neq p_2$.
Then 
\[
\dim M(G,c,\orbit) \ = \ 2n.
\]
Since $n\cdot p_0$ and $(n-1)p_1+p_2$ are rationally equivalent, 
then either $\{p_0\}$ and $\{p_1,p_2\}$ are disjoint, 
or $p_0=p_2$ and $(p_1-p_2)$ has order 
$n-1$. We get a one-parameter family of symplectic leaves $M(G,c,\orbit)$,
depending on the choice of $p_1$ (and we do not allow $p_1-p_0$ to be of
order $n$ because that implies $p_1=p_2$).  
We claim that $M(G,c,\orbit)$ is non-empty, and the generic spectral curve has 
geometric genus $n$. Indeed, consider the $n$-th root of the 
section of ${\cal O}(np_0)$ coming from the non-zero section
of the line bundle 
${\cal O}(np_0-(n-1)p_1-p_2)$ (which is the trivial line-bundle). 
The $n$-th root is a spectral curve  $\overline{C}$ embedded in 
${\cal O}(p_0)$. Then $\overline{C}$ has a cyclic (smooth) 
ramification point over $p_2$ 
and a singularity over $p_1$ of type $y^n=z^{n-1}$. 
The normalization $C$ of $\overline{C}$ is obtained by a single blow-up. 
$C$ is a cyclic $n$-sheeted cover $\pi:C\rightarrow \Sigma$ 
branched over $p_1$ and $p_2$. Hence, the genus of $C$ is $n$. 
Given a line bundle $L$ over $C$,
we get the Higgs bundle $(E,\varphi)$, which is stable because $\overline{C}$
is reduced and irreducible. In particular, $M(G,c,\orbit)$ is non-empty. More
generally, consider any spectral curve $\overline{C}$ 
with characteristic polynomial
$y^n+a_1y^{n-1}+ \cdots + a_n$, whose $i$-th coefficient 
$a_i\in H^0({\cal O}(ip_0))$, $1\leq i<n$, comes from a section of
${\cal O}(ip_0-(i-1)p_1)$ and $a_n$ comes from a non-zero section 
of ${\cal O}(np_0-(n-1)p_1-p_2)$. Note that each $a_i$ is a section of
a line-bundle with a one-dimensional space of global sections. 
A Zariski dense open subset of this 
linear system consists of reduced and irreducible spectral curves. If
$a_1\neq 0$, such a curve has one smooth branch over $p_2$ and a singularity 
over $p_2$ of multiplicity $n-1$ at the zero point in the fiber of 
${\cal O}(p_0)$. A line-bundle $L$, over the blow-up of $\overline{C}$ 
at that singularity, corresponds to a Higgs pair $(E,\varphi)$ in $M(G,c,\orbit)$. 
}
\end{example}

The elliptic Calogero-Moser system for $GL(n)$ is the one obtained from
$\Higgs_\Sigma(GL(n),c,\orbit')$, where $\orbit'$
is the coadjoint orbit of $\mbox{diag}(1,1,\cdots, 1,1-n)$
\cite{donagi-sw,donagi-witten,krichever,treibich-verdier-elliptic-solitons,
treibich-verdier-krichever-variety}. 
(Strictly speaking, one should consider the
$PGL(n)$ Hitchin system obtained from the $GL(n)$ system by fixing the 
determinant line-bundle and taking the quotient by the $\bbc^\times$-action). 

\begin{lemma}
\label{lemma-example-deforms-to-the-calogero-moser-of-type-A-n}
The Elliptic Calogero-Moser system $\Higgs_\Sigma(GL(n),c,\orbit')$ 
can be obtained as the limit of $M(G,c,\orbit)$ in Example
\ref{example-GL-calogero-mozer-deformed} 
when $p_0=p_1=p_2$. 
\end{lemma}

A related deformation of the Elliptic Calogero-Moser system was considered in 
\cite{audin,braden-marshakov-mironov-morozov}.

\noindent
{\bf Proof:}
We follow a geometric procedure for taking the ``derivative''
of a group-type integrable system to obtain a Lie-algebra-type integrable 
system. 
%
%
%
We can take the derivative in a geometric sense, as follows: both the 
Hitchin system and this group valued one for $GL(n)$ have their Poisson 
geometry encoded in  their associated families of spectral curves and 
sheaves supported over these curves, via the Poisson surface in which the 
curves embed. In both cases under consideration here, the Hitchin case, 
($t=0$), and the group valued case, ($t\ne 0$), the surface is the same, that 
is the ruled surface 
$X:=\bbp[{\cal O}_\Sigma(p_0)\oplus {\cal O}_\Sigma]$ over $\Sigma$. 
The Poisson structure $\psi_t$ will vary, 
along with the symplectic leaf one is 
considering; the latter is given geometrically, as we have seen, by fixing the 
intersection of the spectral curves with the zero divisor of the Poisson 
structure.  These two points of intersection will be labelled by 
$\tilde{p}_1(t)$ and $\tilde{p}_2(t)$. 

For $t\neq 0$, we choose 
a Poisson structure $\psi_t$ with a degeneracy locus 
consisting of two sections $\sigma_0(t)$ and $\sigma_\infty$
of the $\bbp^1$ bundle $X\rightarrow \Sigma$. 
The two points $\tilde{p}_1(t)$ and $\tilde{p}_2(t)$ lie on the section 
$\sigma_0(t)$ over the two points $p_1(t)$ and $p_2(t)$. 
When $t=0$, the degeneracy locus $2\sigma_\infty+f_{p_0}$
of the Poisson structure $\psi_0$  
consists of the section $\sigma_\infty$ with multiplicity $2$ 
and the fiber $f_{p_0}$ of $X$ over $p_0$. 
The degeneration $\psi_t$ is chosen so that
the limits of the two points 
$\tilde{p}_1(t)$ and $\tilde{p}_2(t)$ are two distinct points 
$\tilde{p}_1(0)$ and $\tilde{p}_2(0)$ on the fiber $f_{p_0}$. 

We can construct such a degeneration as follows. 
Denote by $x_\infty$ the point of intersection of 
the fiber $f_{p_0}$ with the section $\sigma_\infty$. 
Let $x_0\in f_{p_0}$ be base point of the linear system
$\linsys{\StructureSheaf{X}(f_{p_0}+\sigma_\infty)}$. 
Let $\tilde{p}_1 \ : \ T\hookrightarrow X$ be the embedding of 
a smooth curve into $X$. Assume that $T$ intersects
the fiber $f_{p_o}$ transversally and does not pass through
the points $x_0$ and $x_\infty$. 
The linear system $\linsys{\Wedge{2}T_X(-\sigma_\infty)}$ is one-dimensional
and the universal curve
\[
U \ \subset \ \linsys{\Wedge{2}T_X(-\sigma_\infty)} \times X
\]
is isomorphic to the blow-up $\beta:U\rightarrow X$ of $X$ at $x_0$. 
We denote 
the proper transform of $T$ in $U$ by $T$ as well. Let
$\pi_1:U\rightarrow \linsys{\Wedge{2}T_X(-\sigma_\infty)}$ 
be the natural projection and $\pi_2:U\rightarrow \Sigma$ the composition
$\pi\circ\beta$. Set $p_1:=\pi_2\circ \tilde{p}_1$.
Let $\mu:\Sigma\rightarrow \Sigma$ be multiplication by $-(n-1)$,
regarding $p_0$ as the origin, and define
\[
p_2 \ \  = \ \ \mu \circ p_1. 
\]
Given $t\in T\setminus [T\cap (\sigma_\infty\cup f_{p_0})]$, 
the fiber of $U$ over
$\pi_1\circ \tilde{p}_1(t)\in \linsys{\Wedge{2}T_X(-\sigma_\infty)}$
is a section $\sigma_0(t):\Sigma\rightarrow X$. 
Set $\tilde{p}_2(t):=\sigma_0(t)(p_2)$. We get a morphism
\[
\tilde{p}_2 \ : \ T \ \rightarrow \ U,
\]
since $U$ is projective and $T$ is smooth. By construction, we have
$
\pi_1\circ \tilde{p}_2 \ \ = \ \  \pi_1\circ \tilde{p}_1.
$
We can identify the image $\tilde{p}_2(T)$ in $U$ explicitly. 
There is a natural isomorphism 
$\eta: \mu^*(\StructureSheaf{\Sigma}(p_0))  \IsomRightArrow 
\StructureSheaf{\Sigma}(\mu^{-1}(p_0))$. 
Composing the inclusion $\StructureSheaf{\Sigma}(p_0)\hookrightarrow
\StructureSheaf{\Sigma}(\mu^{-1}(p_0))$ with $\eta^{-1}$ we get a 
lift of $\mu$ to a morphism of the total spaces of line-bundles:
\[
\tilde{\mu} \ : \ \StructureSheaf{\Sigma}(p_0) \rightarrow 
\StructureSheaf{\Sigma}(p_0).
\] 
The composition of a constant section 
$s:\Sigma\rightarrow \StructureSheaf{\Sigma}(p_0)$ with $\tilde{\mu}$ is 
$(s\circ\mu)$. Hence, 
the morphism $\tilde{p}_2$ is the restriction of $\tilde{\mu}$ to $T$. 
The fiber of $\StructureSheaf{\Sigma}(p_0)$ over $p_0$ is naturally 
identified with the tangent line $T_{p_0}\Sigma$ and the morphism
$\tilde{\mu}$ restricts to the fiber as the differential $d_{p_0}\mu$. 
The latter is multiplication by $(1-n)$. Hence, if $\tilde{p}_1(t_0)$ is a
point of $T\cap f_{p_0}$, then $\tilde{p}_2(t_0)=(1-n)\cdot t_0$. 
Consequently, the following limit calculation for triples holds
\[
\lim_{t\rightarrow t_0}
(\sigma_0(t), \ \tilde{p}_1(t), \ \tilde{p}_2(t))
\ \ \ = \ \ \ 
(\sigma_\infty+\hat{f}_{p_0}, \  \tilde{p}_1(t_0), \  (1-n)\tilde{p}_1(t_0)).
\]


For every $t$, the spectral 
curves intersect the degeneracy locus of the Poisson structure at 
$\tilde{p}_1(t)$ with multiplicity $n-1$ and at $\tilde{p}_2(t)$
with multiplicity $1$. 
When $t=t_0$, we get the spectral curves of the Calogero-Moser system
\cite{donagi-sw,krichever}.
\EndProof


\subsection{Type $D$ examples}
\label{sec-type-D}

We use the technique of section \ref{sec-symplectic-surfaces} to study 
the fibers of the characteristic polynomial map (\ref{eq-char}). 
We study the three examples of symplectic leaves $M(PO(2n),c,\orbit)$,
whose dimension is equal to the dimension of a product of 
partial flag varieties (see Example 
\ref{example-partial-flag-varieties-of-maximal-parabolics}).
Our aim is to describe their spectral curves. 

Let $\epsilon_1, \dots, \epsilon_n; -\epsilon_1, \dots, -\epsilon_n$
be the weights of the standard $2n$-dimensional representation of
$SO(2n)$. Denote by 
$e_1, \dots, e_n$ the basis of co-characters of 
$SO(2n)$ dual to $\epsilon_1, \dots, \epsilon_n$. 
The basis $\Delta$ of simple roots of $\gso(2n)$ consists of
$\alpha_1=\epsilon_1-\epsilon_2$,
\dots, $\alpha_{n-1}=\epsilon_{n-1}-\epsilon_n$, and 
$\alpha_n=\epsilon_{n-1}+\epsilon_n$. 
Hence, the first element $\alpha_1^*$,
in the basis of co-characters of $PO(2n)$ dual to $\Delta$, is $e_1$ and the 
last two elements are
\[
\alpha_{n-1}^* \ = \ \left[\frac{1}{2}\sum_{i=1}^{n} e_i\right]-e_n \ \ \ 
\mbox{and}  \ \ \ 
\alpha_n^* \ = \ \frac{1}{2}\sum_{i=1}^n e_i.
\]
The flag variety of the maximal parabolic subgroup corresponding to 
$\alpha_1^*$ is the quadric in $\bbp^{2n-1}$ of isotropic lines. The 
flag varieties corresponding to $\alpha_{n-1}^*$ and 
$\alpha_n^*$ are the two components of the maximal isotropic 
Grassmannian. Each component has dimension $\frac{n(n-1)}{2}$.
Note that $\alpha_{n-1}^*$ and $\alpha_n^*$ are {\em not} 
co-characters of $SO(2n)$.

We will use the following notation.
The dominant weight $\lambda$ of the standard representation of $SO(2n)$ 
is $\epsilon_1$. Let $W_\lambda$ be the stabilizer of $\lambda$ in
the Weyl group and $\widetilde{W}_\lambda$ the stabilizer of the set
$\{\lambda,-\lambda\}$. 
The quotient $\widetilde{W}_\lambda/W_\lambda$ has order $2$.

\subsubsection{The quadric in $\bbp^{2n-1}$}

We describe the spectral curves for the $(4n-4)$-dimensional moduli space 
$M(SO(2n),c,\orbit)$, where $\orbit_p=G[[t]]$, for $p\not\in \{p_1,p_2\}$, and 
$\orbit_{p_i}$, $i=1,2$, correspond to the co-character $\alpha_1^*$. 

We use the set-up of section \ref{sec-symplectic-surfaces}, with 
$G=\bar{G}=SO(2n)$ and $\lambda$ is the dominant weight of the standard 
$2n$-dimensional representation. Set $S:=\Sigma\times \bbp^1$. 
Let $\sigma_0$, $\sigma_\infty$, $\sigma_1$ and $\sigma_{-1}$ 
be the sections of $S$ corresponding to the points $0$, $1$, $-1$, $\infty$
of $\bbp^1$. 
The surface $S_\lambda(\orbit)$, given in (\ref{eq-S-lambda-bar-orbit}), 
is contained in the blow-up of $S$ along the 
points on $\sigma_0$ and $\sigma_\infty$ over $p_1$ and $p_2$. 
$S_\lambda(\orbit)$ is the complement of the proper transform of the two 
fibers of $S$ over $p_1$ and $p_2$. Let $\iota$ be the element of order $2$ in
$\widetilde{W}_\lambda/W_\lambda$. The involution $\iota$ acts on $S$
interchanging the two sections $\{\sigma_0,\sigma_\infty\}$
and fixing $\sigma_1$ and $\sigma_{-1}$. Denote by 
$\sigma_{0,\infty}$, $\bar{\sigma}_1$  and $\bar{\sigma}_{-1}$ the three 
sections of $S/\iota$. 

Let $C$ be the spectral curve in $X(\orbit)$ of a generic pair in 
$M(SO(2n),c,\orbit)$. 
The generic fiber of (\ref{eq-char}) is the prym of 
(the resolution $C/W_\lambda$ of) 
a curve $C_\lambda$ in $S$ with the following properties:
1) $C_\lambda$ is $\iota$-invariant, 
2) $C_\lambda$ has nodes along $\sigma_1$ and $\sigma_{-1}$, and
3) $C_\lambda$ meets the sections $\sigma_0$ and $\sigma_\infty$
over $p_1$ and $p_2$. 
We identify first the family of quotient curves $C_\lambda/\iota$. 
Consider the linear system
$\linsys{n\sigma_{0,\infty}+f_{p_1}+f_{p_2}}$, 
where $f_p$ is the fiber over $p\in \Sigma$.
A curve $\Gamma$ in this linear system satisfies
$\Gamma\cdot \Gamma=4n$ and $\omega_{S/\iota}$ restricts to $\Gamma$ with degree $-4$.
Thus, $\omega_\Gamma$ has degree $4n-4$ and 
the arithmetic genus of $\Gamma$ is $2n-1$. 
Set $\tilde{p_i}=\sigma_{0,\infty}(p_i)$. 
We have a subvariety $B$, of the linear subsystem 
$\linsys{\StructureSheaf{S/\iota}(n\sigma_{0,\infty}+f_{p_1}+f_{p_2})
(-\tilde{p}_1-\tilde{p}_2)}$, 
consisting of curves passing through $\tilde{p_1}$ and $\tilde{p_2}$ and
intersecting $\bar{\sigma}_1$ and $\bar{\sigma}_{-1}$ tangentially. 
The three divisors $\Gamma\cap \sigma_{0,\infty}$, $\Gamma\cap \bar{\sigma}_1$ and 
$\Gamma\cap \bar{\sigma}_{-1}$ are linearly equivalent on $\Gamma$. 
Hence, their norms in $\Sigma$ are linearly equivalent. 
We conclude that $B$ consists of $16$ components $B_{(x_1,x_{-1})}$,
each consisting of curves intersecting $\bar{\sigma}_1$ tangentially at 
$\bar{\sigma}_1(x_1)$, where $x_1$ is a point on $\Sigma$ satisfying the
linear equivalence $2x_1\sim p_1+p_2$. Similarly,
the curves intersect $\bar{\sigma}_{-1}$ tangentially 
at $\bar{\sigma}_{-1}(x_{-1})$, where 
$2x_{-1}\sim p_1+p_2$.
Let $C_\lambda$ be the inverse image in $S$ of a curve $\Gamma\subset S/\iota$ 
in the subsystem $B_{(x_1,x_{-1})}$. If $\Gamma$ is
smooth at $\bar{\sigma}_1(x_1)$ and $\bar{\sigma}_{-1}(x_{-1})$, 
then $C_\lambda$  has a node over each of the two points. 
Its resolution $C/W_\lambda$ is an unramified double cover $C/W_\lambda\rightarrow \Gamma$. 
$\Prym(C/W_\lambda,\Gamma)$ has dimension $2n-2$. 

We show next, that only the four diagonal components of $B$ contain 
spectral curves of pairs in $M(SO(2n),c,\orbit)$. These are the components
$B_{(x_1,x_{-1})}$, where $x_1=x_{-1}$. 
We will use the set-up of section 
\ref{sec-spectral-curves-and-group-isogenies},
with the isogeny $Spin(2n)\rightarrow SO(2n)$. 
The fundamental group of $SO(2n)$ is $\bbz/2\bbz$. We have the 
topological invariant $\tau:\Sigma^0\rightarrow \pi_1(SO(2n))$, 
given in (\ref{eq-topological-invariant}), of a pair in 
$M(SO(2n),c,\orbit)$. The co-character 
$\alpha_1^*$ maps to the non-trivial class in $\pi_1(SO(2n))$. 
Hence, the homomorphism $\tau$ is surjective. The topological invariant
$\tau$ determines a double cover $f:D_\tau\rightarrow \Sigma$ 
branched over $p_1$ and $p_2$. The pushforward of the structure sheaf
decomposes $f_*\StructureSheaf{{D_\tau}}=
\StructureSheaf{\Sigma}\oplus \StructureSheaf{\Sigma}(-p_0)$, 
for a unique point $p_0$. 
The fiber over $p_0$ satisfies the linear equivalence 
$f^{-1}(p_0)\sim \bar{p}_1+\bar{p}_2$, where $\bar{p}_1$ and $\bar{p}_2$
are the two ramification points of $f$. 
Recall, that a pair $(P,\bar{\varphi})$ in 
$M(SO(2n),c,\orbit)$ lifts to a pair $(f^*P,\varphi)$ over $D_\tau$,
where $\varphi$ is a meromorphic section of $(f^*P)(Spin(2n))$. 
The singularity data $(f^{-1}\orbit)_{\bar{p}_i}$, at each of the 
ramification points, consists of co-characters in the $W$-orbit of $2e_1$. 
Let $s^+$ and $s^-$ be the traces of the two half spin
representations of $Spin(2n)$. We get the rational functions
$s^+(\varphi)$ and $s^-(\varphi)$ on $D_\tau$. 
The highest weights of the two half spin
representations are $\frac{1}{2}(\epsilon_1+\cdots +\epsilon_n)$ and 
$\frac{1}{2}(\epsilon_1+\cdots +\epsilon_{n-1}-\epsilon_n)$. 
We conclude, that both $s^+(\varphi)$ and $s^-(\varphi)$ 
have simple poles at $\bar{p}_1$ and $\bar{p}_2$. 
Their difference $s^+(\varphi)-s^-(\varphi)$ has simple zeroes
along the fiber $f^{-1}(x_1)$ and polar divisor $\bar{p}_1+\bar{p}_2$.
We get the linear equivalence on $D_\tau$ 
\[
f^{-1}(x_1) \ \ \sim \ \ \bar{p}_1+\bar{p}_2.
\]
Considering the sum $s^+(\varphi)+s^-(\varphi)$, 
we get the linear equivalence 
\[
f^{-1}(x_{-1}) \ \ \sim \ \ \bar{p}_1+\bar{p}_2.
\]
We conclude that both $x_1$ and $x_{-1}$ must be equal to $p_0$.

\subsubsection{The two components of the maximal isotropic Grassmannian}
\label{sec-max-isotropic-grassmannian}
Let us describe the spectral curves for the $(n^2-n)$-dimensional moduli space 
 $M(PO(2n),c,\orbit)$, where $\orbit_p=PO(2n)[[t]]$, 
for $p\not\in \{p_1,p_2\}$, and 
$\orbit_{p_i}$, $i=1,2$, correspond to the co-character $\alpha_n^*$. 
(The case, where both singularities correspond to $\alpha_{n-1}^*$, 
is analogous). 

The co-character $\alpha_n^*$ of $PO(2n)$ does not come from $SO(2n)$. 
In section \ref{sec-symplectic-surfaces} we developed a general technique, 
which enables us to study the quotients $C/W_\lambda$, 
of the cameral covers of pairs in $M(PO(2n),c,\orbit)$, where 
$\lambda$ is the dominant weight of the standard $2n$-dimensional 
representation of $SO(2n)$. We follow the technique of section 
\ref{sec-symplectic-surfaces}, with $G=SO(2n)$ and $\bar{G}=PO(2n)$. 
The center $K$ of $SO(2n)$ has order $2$. 
The Galois $K$-cover $f:D_\tau\rightarrow \Sigma$ in (\ref{eq-D-orbit})
is a double cover branched over $p_1$ and $p_2$. 

Recall, that $\lambda$ determines three surfaces,
which are denoted in Section
\ref{sec-symplectic-surfaces} by 
$S_\lambda(\orbit)$, 
$S_\lambda(\orbit,\tau)$, and 
$S_\lambda(f^{-1}\orbit)$
(see (\ref{eq-S-lambda-bar-orbit}), (\ref{eq-S-lambda-orbit}), and 
(\ref{eq-S-lambda-f-1-orbit})). The first two are symplectic surfaces
(Lemma \ref{lemma-symplectic-covering}). 
We will describe the generic spectral curve
$C_\lambda$ in $S_\lambda(\orbit,\tau)$ and its double cover 
$\widetilde{C}_\lambda$ in $S_\lambda(f^{-1}\orbit)$.
Denote by $D_\tau^0$ the complement of the two ramification points.
The surface $S_\lambda(f^{-1}\orbit)$ is a partial compactification 
of $D_\tau^0\times \bbc^\times$.  
The group $K\times K\times (\widetilde{W}_\lambda/W_\lambda)$
acts on $S_\lambda(f^{-1}\orbit)$, extending
the natural action on $D_\tau^0\times \bbc^\times$. 
The first factor of $K$ is
the Galois group of $D_\tau\rightarrow\Sigma$. The second factor of $K$ acts
on $\bbc^\times$ via the restriction of the character $\lambda$.
The involution of the group $\widetilde{W}_\lambda/W_\lambda$
acts on the $\bbc^\times$ factor as the inversion $t\mapsto \frac{1}{t}$. 
The symplectic surface $S_\lambda(\orbit,\tau)$ is the quotient of 
$S_\lambda(f^{-1}\orbit)$ by the diagonal $K$-action and 
$S_\lambda(\orbit)$ is the quotient by the $K\times K$-action. 

Let us obtain first an explicit description of the surface 
$S_\lambda(f^{-1}\orbit)$. 
The $W$-orbit of $\alpha_n^*$ consists of co-characters  
$\frac{1}{2}\sum_{i=1}^n\pm e_i$, 
such that the number of negative coefficients is even.
The character $\lambda$ takes this orbit,
of rational co-characters of $SO(2n)$, 
to the set $\{\bar{\alpha},-\bar{\alpha}\}$,
where $\bar{\alpha}=\frac{1}{2}\alpha$ is half the generator $\alpha$ of
the co-character lattice of $T/\ker\lambda$. 
The local toric model of $S_\lambda(f^{-1}\orbit)$, 
over a neighborhood of a ramification point in $D_\tau$, is
described as follows. The lattice is $\Span_\bbz\{e,\alpha\}$. 
Let $\bar{\sigma}$ be the rational cone generated by 
$\{\bar{e}+\bar{\alpha}, \bar{e}-\bar{\alpha}\}$,
where $\bar{e}=\frac{1}{2}e$ (see (\ref{eq-fan-of-symplectic-surface})). 
The local model is the complement, 
of the singular point, in the toric surface determined by the rational cone
$\bar{\sigma}$.
The local model, of the morphism to $D_\tau$, is the projection
(modulo $\alpha$) to the toric curve with cone spanned by $e$ in 
$\Span_\bbz\{e\}$. 
Note, that the cone $\bar{\sigma}$ is generated also by 
$\{e+\alpha, e-\alpha\}$. 
It is easy to see, that the surface $S_\lambda(f^{-1}\orbit)$ 
is obtained from $D_\tau\times \bbp^1$ as follows. 
Let $\tilde{p}_1$ and $\tilde{p}_2$ be the two ramification points
of $D_\tau$. 
Let $Y$ be the blow-up of $D_\tau\times \bbp^1$ 
at the points in $D_\tau\times\{0,\infty\}$ over 
$\tilde{p}_1$ and $\tilde{p}_2$. 
Denote the four exceptional divisors by 
$E_{1,0}$, $E_{1,\infty}$, $E_{2,0}$, and $E_{2,\infty}$.
The surface $S_\lambda(f^{-1}\orbit)$ is a Zariski open subset of $Y$;
the complement of the proper transforms  
of $D_\tau\times \{0,\infty\}$ and 
$\{\tilde{p}_1,\tilde{p}_2\}\times \bbp^1$. 

The surface $S_\lambda(\orbit,\tau)$ is the quotient of 
$S_\lambda(f^{-1}\orbit)$ under the diagonal action of $K$.
Observe, that the fixed locus in $Y$, under the diagonal involution 
of $K$, consists of the four exceptional divisors. 

$H^{1,1}(Y,\bbz)$ is generated by the four exceptional
divisors, the class $\phi$ of a fiber over $D_\tau$, 
and the class $z$ of the total transform of $D_\tau\times \{0\}$. 
The intersection pairing is determined by $\phi\cdot z=1$, $\phi\cdot\phi=0$,
$z\cdot z=0$, each exceptional divisor has self-intersection $-1$,
and each exceptional divisor has zero intersection with $\phi$ and $z$.
The canonical class of $Y$ is cohomologous to
\[
2\phi-2z+E,
\]
where $E$ is the sum of the four exceptional divisors.

Recall, that a $PO(2n)$-pair $(P,\bar{\varphi})$ over $\Sigma$ lifts to 
a pair $(f^*P,\varphi)$ over $D_\tau$, where $\varphi$ is a meromorphic 
section of $f^*P(SO(2n))$. 
Let $C_\varphi$ be the spectral cover of $(f^*P,\varphi)$
in $X(f^{-1}\orbit)$ 
(see section \ref{sec-spectral-curves-and-group-isogenies}).
The (compact) spectral curve $\widetilde{C}_\lambda$ of $(f^*P,\varphi)$ 
embeds in $S_\lambda(f^{-1}\orbit)$ as a $2n$-sheeted cover of 
$D_\tau$. The curve $\widetilde{C}_\lambda$ is a birational image of 
the quotient $C_\varphi/W_\lambda$, 
acquiring additional nodes, which will be discussed below. 
The morphism (\ref{eq-S-lambda-f-1-orbit}) is $K\times K$ equivariant. Hence, 
$\widetilde{C}_\lambda$ is invariant with respect to the 
diagonal action of $K$.

The cohomology class $\gamma$ of $\widetilde{C}_\lambda$
in $Y$ is determined by the intersection data
\begin{eqnarray}
\nonumber
\gamma\cdot \phi &= & 2n
\\
\label{eq-gamma-dot-z}
\gamma\cdot z &=& 2n
\\
\label{eq-gamma-dot-E}
\gamma\cdot E_{i,j} & = & n.
\end{eqnarray}
The first equation indicates, that $\widetilde{C}_\lambda$ is a $2n$-sheeted
branched cover of $D_\tau$.
Equation (\ref{eq-gamma-dot-z}) indicates, that the image of
$\widetilde{C}_\lambda$ in $D_\tau\times \bbp^1$ intersects the 
zero section with multiplicity $n$ over each of the two 
ramification points $\tilde{p}_1$
and $\tilde{p}_2$ (and the same holds for the infinity section). 
This follows from the fact, that the singularity type 
$f^{-1}(\orbit_{p_i})$, of $\varphi$ at $\tilde{p}_i$, 
is the $W$-orbit of the co-character $2\alpha_n^*$ of $SO(2n)$. 
Equation (\ref{eq-gamma-dot-E}) expresses the fact, that the proper transforms
of the zero and infinity sections are disjoint from $\widetilde{C}_\lambda$.
We get the equality
\[
\gamma \ = \ 
2n\phi + 2nz -n E.
\]
The curve $\widetilde{C}_\lambda$ belongs, in fact, to the linear system
$\pi^*\StructureSheaf{D_\tau}(n\tilde{p}_1+n\tilde{p}_2)\otimes
\StructureSheaf{Y}(2nz-nE)$, where $\pi:Y\rightarrow D_\tau$ is
the natural morphism. 
Consequently, $\gamma\cdot \gamma=4n^2$ and the intersection of 
the canonical class of $Y$ with $\gamma$ is $4n$. 
We get, that the arithmetic genus of $\widetilde{C}_\lambda$ 
is $2(n^2+n)+1$. 

We discuss next the additional nodes acquired via the map 
$C_\varphi/W_\lambda\rightarrow \widetilde{C}_\lambda$. 
Let $T$ be the maximal torus of $SO(n)$ and $x\in T$ an element, 
whose weights $\epsilon_i(x)=x_i$ satisfy: $x_1=-1$,
$\{x_2, \dots, x_n\}$ are distinct, and $x_i^2\neq 1$, for $i\geq 2$. 
The Weyl group orbit of $x$ consists of ``signed'' permutations
$w(x)_i=x_{w(i)}^{\pm 1}$, with an even number of inversions. 
Consequently, the $\widetilde{W}_\lambda$-orbit of $x$ 
consists of two distinct $W_\lambda$-orbits 
(the two orbits are distinguished by the parity of the number of
inversions in the set $\{x_2, \dots, x_n\}$). 
We conclude, that the generic quotient $C_\varphi/W_\lambda$ 
has two distinct points over each of the intersection points, of 
$\widetilde{C}_\lambda$, with the transform of $D_\tau\times \{-1\}$.
A similar discussion, of the case $\lambda(x)=1$, 
leads to the following conclusion. 
The curve $\widetilde{C}_\lambda$, of a generic cameral cover, 
meets the fixed locus of $\widetilde{W}_\lambda/W_\lambda$ in 
$S_\lambda(f^{-1}\orbit)$, along nodes of $\widetilde{C}_\lambda$. 
Hence, $\widetilde{C}_\lambda$ has $\frac{1}{2}(\gamma\cdot z)=n$ 
nodes along each of these two $\widetilde{W}_\lambda/W_\lambda$-invariant 
sections. The generic $C_\varphi/W_\lambda$ has geometric genus 
$2n^2+1$.

The diagonal $K$-action acts on $C_\varphi/W_\lambda$ with $4n$ 
fixed points over the four exceptional divisors in $Y$. 
These are ramification points of 
the quotient $q:C_\varphi/W_\lambda\rightarrow C/W_\lambda$. 
Hence, $C_\varphi/[W_\lambda\times K]$ has geometric 
genus $n^2-n+1$. 

The curve $C_\varphi/[W_\lambda\times K]$ maps to the spectral curve $C_\lambda$ in $S_\lambda(\orbit,\tau)$ with $n$ nodes. 
Lemma \ref{lemma-symplectic-covering} implies, that the fixed locus in 
$S_\lambda(\orbit,\tau)$, of the $\widetilde{W}_\lambda/W_\lambda$-action, 
is precisely the $K$-quotient of $[D_\tau^0\times \{1,-1\}]$. 
We conclude, that $\widetilde{W}_\lambda/W_\lambda$ acts on the generic curve 
$C_\varphi/[W_\lambda\times K]$ without fixed points. 
The prym variety of the $\widetilde{W}_\lambda/W_\lambda$ 
action has dimension $\frac{n^2-n}{2}$. 
This is precisely half the dimension of the moduli space
$M(PO(2n),c,\orbit)$.







\begin{thebibliography}{B-N-R}
\bibitem[AB]{atiyah-bott} M. F. Atiyah and R. Bott: 
{\it The Yang-Mills equations over Riemann surfaces.\/}
Philos. Trans. Roy. Soc. London Ser. A 308 (1983), no. 1505, 523-615. 

\bibitem[AK]{altman-kleiman} A. Altman and S. Kleiman: 
{\it Compactifying the Picard Scheme,\/ }
Adv. in Math. {\bf 35}, 50-112 (1980)

\bibitem[Ar]{artamkin} I. V. Artamkin:
{\it On deformations of sheaves,\/}
Math. USSR Izvestiya Vol. 32(1989), No. 3 663-668

\bibitem[At]{atiyah} M. F. Atiyah: 
{\it Vector bundles over an elliptic curve},
Proc. Lond. Math. Soc {\bf 7} (1957), 414--452.

\bibitem[Au]{audin} M. Audin: 
{\it Lectures on gauge theory and integrable systems.\/} 
in Gauge theory and symplectic geometry (Montreal, PQ, 1995), 1-48, 
NATO Adv. Sci. Inst. Ser. C Math. Phys. Sci., 488, 
Kluwer Acad. Publ., Dordrecht, 1997. 

\bibitem[BD1]{beilinson-drinfeld} A. Beilinson, V. G. Drinfeld: 
{\em Quantization of Hitchin's integrable system
and Hecke eigenshieves} Preprint. 


\bibitem[BD2]{belavin-drinfeld} A. A. Belavin and  V. G. Drinfeld:
{\it Solutions of the classical Yang-Baxter equations for simple Lie 
algebras.\/}
Funct. Anal. and its appl., {\bf 16 } (1982), 159-180

\bibitem[BFM]{BFM} A. Borel, R. Friedman, J. W. Morgan:
{\em Almost commuting elements in compact Lie groups.\/}
math.GR/9907007

\bibitem[Ber]{bernard} D. Bernard:
{\em On the Wess-Zumino-Witten models on the torus.} 
Nucl. Phys. {\bf B303}, 77-93 (1988)

\bibitem[BR]{biswas-ramanan} 
I. Biswas and S. Ramanan:
{\it 
An infinitesimal study of the moduli of Hitchin pairs}.
J. London Math. Soc. (2) {\bf 49} (1994), no. 2, 219--231. 

\bibitem[Bo1]{Bo1}
F. Bottacin: 
{\it Symplectic geometry on moduli spaces of stable pairs,}
Ann. Sci. Ecole Norm. Sup. (4) {\bf 28} (1995), no. 4, 391-433.

\bibitem[Bo2]{Bo2}  F. Bottacin: 
{\it Poisson structures on moduli spaces of sheaves over Poisson
surfaces}, Invent. Math. {\bf 121} (1995), no. 2, 421-436. 

\bibitem[BMMM]{braden-marshakov-mironov-morozov}
H. W. Braden, A.Marshakov, A.Mironov, A.Morozov: 
{\it The Ruijsenaars-Schneider Model in the Context of Seiberg-Witten 
Theory,\/} Nucl.Phys. B558 (1999) 371-390


\bibitem[Br]{bridgeland} T. Bridgeland:
{\it Equivalences of triangulated categories and Fourier-Mukai transforms.}
Bull. London Math. Soc. 31 (1999), no. 1, 25--34. 

\bibitem[Do1]{donagi}  R. Donagi: {\it Spectral covers},  in: 
Current topics in complex algebraic geometry (Berkeley, CA 1992/93), 
Math. Sci. Res. Inst. Publ. {\bf 28} (1995), 65-86. 

\bibitem[Do2]{donagi-sw}
R. Donagi: {\em Seiberg-Witten integrable systems.\/}
Algebraic geometry---Santa Cruz 1995, 3--43, Proc. Sympos. Pure Math., 62, 
Part 2, Amer. Math. Soc., Providence, RI, 1997. 
(alg-geom eprint no. 9705008)

\bibitem[DG]{donagi-gaitsgory} R. Donagi and D. Gaitsgory:
{\em The gerbe of Higgs bundles,\/}
Transform. Groups 7 (2002), no. 2, 109--153. 


\bibitem[DM]{DM} R. Donagi and E. Markman: {\it Spectral curves,
algebraically completely integrable Hamiltonian systems, and
moduli of bundles}, 
Springer Lecture Notes in Math. Vol 1620, pages 1-119 (1996).

\bibitem[Dr1]{drinfeld-icm}  V. G. Drinfeld: 
{\it Quantum groups}, Proc. Int. Congr. Math. Berkeley, 1986, 1, pp. 798-820.

\bibitem[Dr2]{drinfeld-poisson-groups} V. G. Drinfeld: 
{\em On Poisson homogeneous spaces of Poisson Lie-groups}, 
Theoret. and Math. Phys. 95 (1993), no. 2, 524-525 


\bibitem[DW]{donagi-witten}  R. Donagi and E. Witten:
{\it Supersymmetric Yang-Mills theory and
integrable systems}, Nucl.Phys B460 (1996), 299. HepTh, 9510101.

\bibitem[ES]{etingof-schiffman} P. Etingof and O. Schiffmann: 
{\em Twisted traces of intertwiners for 
Kac-Moody algebras and classical dynamical R-matrices
corresponding to generalized Belavin-Drinfeld triples.\/}
Math. Research Letters {\bf 6}, 593-612 (1999)

\bibitem[EV]{etingof-varchenko} P. Etingof and A. Varchenko:
{\it Geometry and classification of solutions of the classical
Dynamical Yang-Baxter Equation,\/}
Commun. Math. Phys. 192, 77-120 (1998)

\bibitem [Fa]{faltings}   G.~Faltings: 
{\it Stable $G$-bundles and projective
connections},
 J. Alg. Geom.  {\bf 2}(1993), 507-568.  

\bibitem [Fe]{felder} G. Felder:
{\em Conformal field theory and integrable systems associated to elliptic 
curves,\/}
Proceedings of the International Congress of Mathematicians, Vol. 1, 2 
(Z\"{u}rich, 1994), 1247--1255, 
Birkh\"{a}user, Basel, 1995. 

\bibitem[FM]{FM}
R. Friedman and J. W. Morgan:
{\it Holomorphic principal bundles over elliptic curves},
math.AG/9811130


\bibitem[Hi1]{hitchin-self-duality} N. J. Hitchin:
{\em The self-duality equations on a Riemann surface},
Proc. Lond. Math. Soc. {\bf 55} (1987) 59-126.

\bibitem[Hi2]{hitchin-integrable-system} N. J. Hitchin:
{\em Stable bundles and integrable systems.\/}
Duke Math. J. 54, No 1  91-114 (1987)

\bibitem[H]{hurtubise-surfaces}
J. Hurtubise: 
{\it Integrable systems and algebraic surfaces},  
Duke Math. J. {\bf 83} (1996), no.~1, 19--50.



\bibitem[HM1]{HM-pryms} 
J. Hurtubise and E. Markman: 
{\it Rank $2$ integrable systems of Prym varieties},
  Adv. Theo. Math. Phys.,{\bf 2} (1998), 633-695. 

\bibitem[HM2]{HM-sklyanin} 
J. Hurtubise and E. Markman: 
{\it Surfaces and the Sklyanin bracket. \/}
Comm. Math. Phys. 230 (2002), no. 3, 485--502. 

\bibitem[Hum1]{humphreys}
J. E. Humphreys: 
{\it Introduction to Lie algebras and representation theory},
Springer-Verlag (1972).

\bibitem[Hum2]{humphreys-groups}
J. E. Humphreys: 
{\it Linear Algebraic Groups},
Springer-Verlag GTM 21 (1991).

\bibitem[Hum3]{humphreys-conjugacy}
J. E. Humphreys: {\it Conjugacy classes in semisimple algebraic groups.}
Mathematical Surveys and Monographs, 43. 
AMS 1995. 

\bibitem[Kap]{kapranov} M. Kapranov: 
{\it Hypergeometric functions on reductive groups,\/} 
Integrable systems and algebraic geometry (Kobe/Kyoto, 1997), 236--281, 
World Sci. Publishing, River Edge, NJ, 1998. 

\bibitem[Kaw]{kawamata-T1-lifting} Y. Kawamata 
{\it Unobstructed deformations - a remark on a paper of Z. Ran.\/}
J. Alg. Geom. 1 (1992) 183-190. 
{\it Erratum on ``Unobstructed deformations'',\/}
J. Alg. Geom. 6 (1997) 803-804

\bibitem[Kr]{krichever} I. M. Krichever: 
{\em Elliptic solutions of the K-P equation and integrable systems of 
particles.\/} Funct. Anal. 14, (1980), p. 45-54

\bibitem[Lan]{lange} H. Lange: {\em Universal families of Extensions.}
J. of Alg. 83, 101-112 (1983)


\bibitem[Li]{li} L. Li: Letter, October 2001.

\bibitem[Lo]{looijenga} E. Looijenga:   
{\it Root systems and elliptic curves}, Inv. 
Math. 38(1976),17-32 and {\it Invariant Theory for generalized root 
systems}, Inv. Math. 61(1980),1-32

\bibitem[Lu]{lusztig} G. Lusztig: 
{\em Singularities, character formulas, and a $q$-analog
of weight multiplicities} analyse et topologie sur les espaces singuliers, 
Ast\'{e}risque 101-102, 1982, pp. 208-229.  

\bibitem[M1]{M}
E. Markman:   
{\it Spectral curves and integrable systems}, 
Compositio Math. {\bf 93} (1994), 255-290.


\bibitem[M2]{markman-sw}
E. Markman:    
{\it Integrable Systems, Algebraic Geometry and Seiberg-Witten Theory.
} 
In "Integrability: the Seiberg-Witten and Whitham
equations", eds H.W. Braden and I.M. Krichever, 23-41.
Amsterdam: Gordon and Breach Science Publishers. 1999
 

\bibitem[Mu1]{mukai-symplectic} S. Mukai: 
{\em Symplectic structure of the moduli
space of sheaves on an abelian or K3 surface},
Invent. math. 77,101-116 (1984)

\bibitem[Mu2]{mukai-survey} S. Mukai: 
{\it Moduli of vector bundles on $K3$ surfaces and symplectic manifolds.\/}
(Japanese) Sugaku Expositions {1} (1988), no. 2, 139--174. 
Sugaku 39 (1987), no. 3, 216--235. 

\bibitem[Mu3]{mukai-fourier} S. Mukai: 
{\it Duality between $D(X)$ and $D(\hat X)$ with its application to Picard
sheaves}. Nagoya Math. J. 81 (1981), 153--175. 

\bibitem[N]{nitsure} N. Nitsure: 
{\em Moduli space of semistable pairs on a curve.\/}
Proc. London Math. Soc. (3) 62 (1991) 275-300

\bibitem[O]{orlov} D. O. Orlov: 
{\it Equivalences of derived categories and K3 surfaces, \/}
Algebraic geometry, 7.
J. Math. Sci. (New York) 84 (1997), no. 5, 1361-1381. 

\bibitem [Po]{Po} A. Polishchuk: 
{\it Poisson structures and birational morphisms associated with
bundles on elliptic curves}, Internat. Math. Res. Notices (1998), no. 
{\bf 13}, 683--703. 

\bibitem[R1]{ziv-ran} 
Z. Ran:  
{\it On the local geometry of moduli spaces of locally free sheaves. \/}
Moduli of vector bundles (Sanda, 1994; Kyoto, 1994), 213--219, 
Lecture Notes in Pure and Appl. Math., 179, Dekker, New York, 1996.

\bibitem[R2]{ziv-ran-T1-lifting} Z. Ran: 
{\it Deformations of manifolds with torsion or negative canonical bundle.\/} 
J. Alg. Geom. 1 (1992), no. 2, 279--291. 

\bibitem[Ra]{ramanathan} A. Ramanathan, 
{\it Stable Principal Bundles on a Compact Riemann Surface,\/}
Math. Ann. {\bf 213}, 129--152 (1975).

\bibitem[RS]{reiman-semenov} 
A. G. Reiman and M.A. Semenov-Tian-Shansky: 
{\it Integrable Systems II}, chap.2, in ``Dynamical Systems VII'', 
Encyclopaedia of Mathematical Sciences, vol 16., 
V.I. Arnold and S.P.Novikov, eds., Springer-Verlag, Berlin, 1994.

\bibitem[Sco]{scognamillo}  R. Scognamillo: 
{\it An elementary 
approach to the abelianization of the Hitchin system for arbitrary 
reductive groups.}

\bibitem[Sim]{simpson} C. Simpson: 
{\it Moduli of representations of the fundamental group of a smooth
projective variety. I and II.}
Inst. Hautes \'{E}tudes Sci. Publ. Math. No. 79 (1994), 47--129 
and 
Inst. Hautes \'{E}tudes Sci. Publ. Math. No. 80 (1994), 5--79. 


\bibitem[Sk]{sklyanin}
E. K. Sklyanin:
{\it On the complete integrability  of the Landau-Lifschitz equation},
LOMI preprint E-3-79, (1979).

\bibitem[TV1]{treibich-verdier-elliptic-solitons} 
A. Treibich, J.-L. Verdier:
{\it Solitons elliptiques.\/}
The Grothendieck Festschrift, Vol. III, 437-480,
Birkhauser Boston, 1990.


\bibitem[TV2]{treibich-verdier-krichever-variety} 
A. Treibich, J.-L. Verdier:
{\it  Varietes de Kritchever des solitons elliptiques de KP.\/}
Proceedings of the Indo-French Conference on Geometry
(Bombay, 1989), 187-232, 1993.

\bibitem[Ty]{tyurin} A. Tyurin: 
{\em Symplectic structures on the varieties of moduli of
vector bundles on algebraic surfaces with $p_{g} > 0$. \/}
Math. USSR Izvestiya Vol. 33(1989), No. 1


\end{thebibliography}
\end{document}